\documentclass[12pt,twoside]{article}
\usepackage{amsmath,amssymb,mathrsfs,graphicx,yhmath,fancyhdr,indentfirst,newlfont,amscd,amsthm,floatflt,bbm,subfigure,tikz-cd,comment,float,xcolor}
\usepackage[small,sc]{caption} 
\usepackage[title]{appendix}
\usepackage[hidelinks]{hyperref}

\newtheorem{theorem}{Theorem}[section]
\newtheorem{lemma}[theorem]{Lemma}
\newtheorem{proposition}[theorem]{Proposition}
\newtheorem{corollary}[theorem]{Corollary}
\newtheorem{definition}[theorem]{Definition}

\newtheorem{rmrk}[theorem]{Remark}

\newtheorem*{theorem*}{Theorem}
\newtheorem*{corollary*}{Corollary}
\newtheorem*{notation*}{Notation}

\DeclareMathAlphabet{\mathbfit}{OML}{cmm}{b}{it}

\makeatletter
\@addtoreset{equation}{section}
\makeatother

\newenvironment{remark}
{\begin{rmrk} \em}
	{\end{rmrk}}

\newcommand{\R} {\mathbb{R}}
\newcommand{\C} {\mathbb{C}}

\newcommand{\Z} {\mathbb{Z}}
\newcommand{\N} {\mathbb{N}}
\renewcommand{\qed} {\hfill {\small Q.E.D.} \par\medskip}

\renewcommand{\proof} {\noindent \textsc{Proof.} }
\newcommand{\proofof}[1] {\noindent \textsc{Proof of {#1}.} }

\renewcommand{\top}{\text{top}}

\DeclareMathOperator{\de}{\partial}

\renewcommand{\emptyset} {\varnothing}

\newcommand{\norm}[1]{\Vert#1\Vert}

\newcommand{\id}{\textnormal{id}}

\newcommand{\ess}{\textnormal{ess}}

\newcommand{\diam}{\textnormal{diam}}
\newcommand{\T}{\mathbb{T}}

\definecolor{miorosso}{RGB}{211, 47, 47}
\definecolor{mioverde}{RGB}{56, 140, 131}
\definecolor{mioblu}{RGB}{21, 101, 192}
\definecolor{mioarancio}{RGB}{224, 116, 21}

\begin{document}
	
	\title{\textbf{A cohomological approach to Ruelle-Pollicott resonances and speed of mixing of Anosov diffeomorphisms}}
	
	\author{
		\scshape
		Daniele Galli\thanks{
			Universität Zürich, Institut für Mathematik, Winterthurerstrasse 190,
			CH-8057 Zürich, Switzerland.
			E-mail: \href{mailto:daniele.galli@math.uzh.ch}{\texttt{daniele.galli@math.uzh.ch}}.
		}\\
        Final version for\\
        \textit{Communications in Mathematical Physics}}
	
	\date{May 2026}
	
	\maketitle
	%\updateinfo
	
	\begin{abstract}
		We investigate Ruelle-Pollicott resonances of smooth Anosov diffeomorphisms, acting on manifolds of every dimension, with respect to the measure of maximal entropy. We highlight a profound connection between resonances and eigenvalues of the action induced by the dynamics on de Rham cohomology. In particular, resonances appear as eigenvalues of a quasi-compact transfer operator acting on suitable anisotropic spaces  of currents. After defining the anisotropic Banach spaces, we introduce the anisotropic de Rham cohomology and we show that it is isomorphic to the \emph{standard} de Rham cohomology. The relation between resonances and cohomological eigenvalues is deduced from a comparison of spectra. We finally exploit these results to get information about the Ruelle-Pollicott asymptotics of the correlation function and to establish a cohomological bound for the speed of mixing of Anosov diffeomorphisms.
	\end{abstract}
	
	\section{Introduction}
	The statistical properties of Anosov diffeomorphisms, whose definition is recalled in Section \ref{section_Setup_results}, represent one of the main topics in the field of (uniformly) hyperbolic dynamics, because they model many chaotic phenomena. Most of the ergodic properties of Anosov diffeomorphism were established, starting from the 60's, with the important contribution of D.V. Anosov \cite{anosov}, S. Smale \cite{smale_differentiable}, R. Bowen \cite{bowen}, D. Ruelle \cite{ruelle_thermodynamic} and Ya. G. Sinai \cite{sinai1,sinai2}. In particular, it is well known that a topologically transitive Anosov diffeomorphism $f\colon M\rightarrow M$ of class $C^2$, on a compact, connected Riemannian manifold, admits a unique measure of maximal entropy\footnote{The measure of maximal entropy is also called Bowen-Margulis measure, after R. Bowen \cite{bowen2,bowen3} and G. Margulis \cite{margulis}, who gave two different, but equivalent, constructions of this invariant measure. That is why we use the notation $\mu_{BM}$.} $\mu_{BM}$ and the dynamics turns out to be mixing. In addition, it is a classical result that, for H\"older observables, the system is exponentially mixing, i.e.,  there exists $\sigma\in (0,1)$ such that\footnote{Here $C_{\phi,\psi}$ is a constant depending on the H\"older norm of $\phi$ and $\psi$, but it does not depend on $n$. See the notation remark at the end of the introduction.} 
	\begin{equation}
		\label{eq_decay_of_correlation_Holder}
			\left|\int_M \phi(\psi\circ f^n) d\mu_{BM}-\int_M\phi d\mu_{BM}\int_M\psi d\mu_{BM}\right|\leq C_{\phi,\psi} \sigma^n,
	\end{equation}
	for every couple of H\"older observables $\phi,\psi\in C^\alpha(M).$ 
	The \textit{classical} approach of the above-mentioned authors is based on coding the dynamics with a suitable Markov partition and studying the ergodic properties of the conjugated subshift of finite type, which in turn are consequences of the spectral properties of an appropriate transfer operator (see for instance \cite{baladi} for a complete overview of these techniques). The low regularity of the conjugacy map, that is in general only H\"older continuous\footnote{The stable/unstable manifold are merely H\"older continuous (see for instance \cite{moser}). Since the Markov partition is made of rectangles whose sides are pieces of stable/unstable manifolds, the conjugacy cannot be smoother than H\"older in general.}, prevents us from getting information stronger than \eqref{eq_decay_of_correlation_Holder}, even by taking very smooth observables. A groundbreaking  approach, introduced in \cite{blank-keller-liverani} and developed with the main contribution of C. Liverani, S. Gou\"ezel, V. Baladi, M. Tsujii et al. \cite{baladi2,baladi-tsujii,baladi-tsujii2,giulietti-liverani,giulietti-liverani-pollicott,gouezel-liverani,gouezel-liverani08, liverani,liverani-tsujii}, allowed researchers to overcome the regularity problem, by studying directly the spectral properties of the transfer operator,  without coding the system.  The core of their method consists in the definition of proper \textit{anisotropic} Banach spaces and a subsequent direct analysis of the spectrum of a specific transfer operator acting on it. 
	 Indeed, by using \textit{modern} techniques,  S. Gou\"ezel and C. Liverani \cite{gouezel-liverani08}, as well as V. Baladi and M. Tsujii \cite{baladi-tsujii,baladi-tsujii2}, proved independently  the following  result, that we rewrite in the setting of this paper.
	\begin{theorem}
		Let $f\in C^\infty$ be a topologically transitive Anosov diffeomorphism on a compact, connected, $C^\infty$ Riemannian manifold  $M$. Then there exists a unique measure of maximal entropy $\mu_{BM}$. In addition, for every $\epsilon>0$, there exists a finite set $\Xi=\{\xi_1,\dots, \xi_{n_{\epsilon}}\}\subset\C$, with $\epsilon<|\xi_i|<1,$ such that, for every $\xi_i,$ there are a finite number $N_i\in\Z^+$ of finite-rank bilinear forms $\{c_{\xi_i,k}(\cdot,\cdot)\}_{k=1}^{N_{i}}$ for which
		\begin{equation}
			\label{eq_ruelle_asymptotics}
			\int_M\phi(\psi\circ f^n) d\mu_{BM}=\int_M\phi d\mu_{BM}\int_M\psi d\mu_{BM}+\sum_{i=1}^{n_\epsilon}\sum_{k=1}^{N_i}\xi_i^n n^k c_{\xi_i,k}(\phi,\psi)+o(\epsilon^n),
		\end{equation}
		for any $\phi,\psi\in C^\infty(M).$
	\end{theorem}
	The expansion \eqref{eq_ruelle_asymptotics} is generally called a Ruelle-Pollicott asymptotics, after D.Ruelle \cite{ruelle3,ruelle4} and M.Pollicott \cite{pollicott}, and the complex numbers $\{1, \xi_1,\dots,\xi_{n_\epsilon}\}$ are called Ruelle-Pollicott resonances. We also underline that the  above theorem only proves the existence of the asymptotics of the correlation function, but it does not give any information about the actual existence of Ruelle-Pollicott resonances and about their location. In fact, even if the existence of the Ruelle-Pollicott resonances is established for many dynamical systems, both in the discrete-time case \cite{gouezel-liverani08,baladi-tsujii,baladi-tsujii2,faure_gouezel_lanneau} and in the continuous-time case \cite{tsujii,baladi-demers-liverani,dyatlov-guillarmou}, there are only few examples that admit a complete localization of the resonances.
	For instance, it is well known that, for hyperbolic automorphisms of tori, $\Xi=\emptyset$ and the set of Ruelle-Pollicott resonances is reduced to $\{1\}.$ This fact can be easily proved by  Fourier analysis. On the other hand, there exist Anosov diffeomorphisms which admit nontrivial Ruelle-Pollicott resonances for the SRB-measure \cite{adam,slipantschuk,pollicott2}, or for which there is an estimate on the number of resonances \cite{jezequel}. Similarly, there are examples of continuous-time hyperbolic for which is known that the Ruelle-Pollicott spectrum belongs to vertical bands of the complex plane \cite{faure-tsujii,datchev-dyatlov-zworsky,cekic-guillarmou}, or to the union of the real axis with countably many vertical lines in the case of geodesic flow on constant negative curvature hyperbolic manifolds \cite{dyatlov-faure-guillarmou}. It is important to underline that examples admitting a precise location of the Ruelle-Pollicott spectrum are, in general, linear (or affine) hyperbolic dynamical systems \cite{butterley-kiamari-liverani,faure_gouezel_lanneau}  or group actions on homogeneous spaces \cite{dyatlov-faure-guillarmou}. On the other hand, in this paper we consider nonlinear Anosov diffeomorphisms and we can find a correspondence between certain Ruelle-Pollicott resonances and the eigenvalues of the linear map induced by $f$ on the cohomology.
	Moreover, we treat the measure of maximal entropy, for which even less results are known about the existence and location of Ruelle-Pollicott resonances, and a complete asymptotics like \eqref{eq_ruelle_asymptotics} is proved only in the trivial case of hyperbolic automorphisms on tori, where the Bowen-Margulis measure coincides with the SRB-measure, which is the Lebesgue measure. Some new ideas to face this problem arose from \cite[Section 5]{butterley-kiamari-liverani},  where the authors proved the following theorem for transitive Anosov diffeomorphisms of the 2-torus. 	
	\begin{theorem}
		\label{thm_liverani}
		Let $f$ be a $C^\infty$ topologically transitive Anosov diffeomorphism of $\T^2.$
		There exist $r\in\N$, $C>0$ and $k\in(0,1)$ such that, for any couple of observables  $\phi,\psi\in C^\infty(\mathbb{T}^2),$
		$$\left|\int_{\mathbb{T}^2}\phi(\psi\circ f^n) d\mu_{BM}-\int_{\mathbb{T}^2}\phi d\mu_{BM}\int_{\mathbb{T}^2}\psi d\mu_{BM}\right|\leq C\norm{\phi}_{C^r}\norm{\psi}_{C^r}(ke^{-h_\top })^n,$$
		where $h_\top$ is the topological entropy of $f$.
	\end{theorem}   
	The above Theorem \ref{thm_liverani} shows that, for the measure of maximal entropy, there are no Ruelle-Pollicott resonances in the annulus $\{z\in\C:\ e^{-h_\top}\leq |z|\leq1\},$ except for the trivial resonance $\{1\}.$ 
	
	 As a corollary of our main result, we prove a generalization of Theorem \ref{thm_liverani} for every transitive Anosov diffeomorphism on higher-dimensional tori.  We point out that  Theorem \ref{thm_liverani} was previously proved by V. Baladi in \cite{baladi3} and later  by G. Forni in \cite{forni}. For different reasons, \cite{baladi3} and \cite{forni} are specific to the 2-torus and the techniques cannot be extended to higher dimensional manifolds.  
	
	We now explain our main results, which will be properly stated in Section \ref{section_Setup_results}. 	Let $M$ be an orientable, closed, $\dim(M)-$dimensional Riemannian manifold and $f\colon M\rightarrow M$ a $C^\infty$, topologically transitive, Anosov diffeomorphism. Let $\lambda>1$ be the least expansion of unstable vectors and $\lambda^{-1}$ be the least contraction of stable vectors as in the definition of Anosov diffeomorphisms recalled at the beginning of Section \ref{section_Setup_results}.  We recall that every smooth diffeomorphism $f$ of a compact manifold induces a linear action $f_{\#}$ on the finite-dimensional de Rham cohomology (see Section \ref{section_Setup_results} for few basic facts about de Rham cohomology). Let $\sigma(f_\#|_{H_{dR}^{d_s}(M)})$ be the  spectrum of $f_\#$ acting on the de Rham cohomology group $H_{dR}^{d_s}(M),$ where $d_s$ is the dimension of the stable subbundle of the Anosov diffeomorphism $f$. Let us consider the finite subset of cohomological eigenvalues $\{\Lambda_1=e^{h_\top},\Lambda_2,\dots,\Lambda_m\}\subseteq\sigma(f_\#|_{H_{dR}^{d_s}(M)})$ such that\footnote{Here $h_\top$ denotes the topological entropy of $f$.} $\Lambda_1>|\Lambda_2|\geq\dots\geq|\Lambda_m|>\lambda^{-1}e^{h_\top}.$ In our main Theorem \ref{thm_strong}, we prove the following partial Ruelle-Pollicott asymptotics 
		\begin{equation*}
		\begin{split}
			\left|\int_{M}\phi(\psi\circ f^n) d\mu_{BM}-\right.&\int_{M}\phi d\mu_{BM}\int_{M}\psi d\mu_{BM}-\\-&\left.\sum_{i=2}^{m}\sum_{k=0}^{N_i-1}(\Lambda_ie^{-h_\top})^nn^k c_{\Lambda_i,k}(\phi,\psi)\right|\leq C_{\phi,\psi}\lambda^{-n},
		\end{split}
	\end{equation*}
	for any couple of smooth observables $\phi,\psi\in C^\infty.$ In other terms, we obtain that the set of Ruelle-Pollicott resonances in the annulus $\{z\in\C:\ |z|\in(\lambda^{-1},1]\}$
	comes from the action of the dynamics on cohomology.
	
	A straightforward consequence of the main Theorem \ref{thm_strong} is a bound on the speed of exponential decay of correlation (Corollary \ref{cor_gen}). In fact, denoting by $\nu=\max\{|\Lambda_2|,\lambda^{-1}e^{h_\top}\},$ we get
		\begin{equation*}
		\left|\int_{M}\phi(\psi\circ f^n) d\mu_{BM}-\int_{M}\phi d\mu_{BM}\int_{M}\psi d\mu_{BM}\right|\leq C_{\phi,\psi}\nu^ne^{-nh_\top},
	\end{equation*}
		for any couple of smooth observables $\phi,\psi\in C^\infty.$  
		
		Moreover, whenever $f$ is topologically conjugated to a hyperbolic automorphism of a torus, for instance when $f$ is an Anosov diffeomorphism on a torus \cite{franks_anosov_diffeonorphisms_of_tori} or, in general, if the system is codimension 1 \cite{franks,newhouse} \footnote{Codimension 1 means that the dimension of either the stable or the unstable manifolds is 1.} , we prove $|\Lambda_2|<\lambda^{-1}.$ Accordingly, whenever $M$ is a torus,
			\begin{equation}
			\label{decay_of_correlation2}
			\left|\int_{M}\phi(\psi\circ f^n) d\mu_{BM}-\int_{M}\phi d\mu_{BM}\int_{M}\psi d\mu_{BM}\right|\leq C_{\phi,\psi}\lambda^{-n}.
		\end{equation}
		Thus, there are no Ruelle-Pollicott resonances in $\{z\in\C:\ \lambda^{-1}<|z|=1\}\setminus\{1\}$ and Corollary \ref{cor_gen} is actually equivalent to Theorem \ref{thm_strong}. This is a consequence of the analysis of the spectrum of the action of $f$ on de Rham cohomology, which is easily accessible in the particular case of Anosov diffeomorphisms on tori (see the proof of Theorem \ref{cor_main}). 
		We point out that Theorem \ref{thm_strong} could be improved to include Anosov diffeomorphisms and observables of regularity $C^r$, with $r\geq 2.$ Since this extension only adds more technicalities, but it does not require any new significant idea, we limit to consider $C^\infty$ Anosov diffeomorphisms and $C^\infty$ observables.  For more details, see Remark \ref{rem_finite_regularity}.
	
	The paper is organized as follows. In Section \ref{section_Setup_results} we set some assumptions for the dynamical system, we recall some basic concepts of de Rham cohomology  and state our main results. Section \ref{sec_anisotropic_banach_spaces_of_currents} is devoted to the construction of suitable 
	anisotropic Banach spaces of currents. In particular, in Section \ref{sec_quasi-compactnes} we prove that the pushforward operator acts quasi-compactly on our anisotropic Banach spaces. Section \ref{sec_spectrum} and \ref{sec_anisotropic_de_rham_cohomology_and_spectrum} represent the core of the paper. In fact, in Section \ref{sec_spectrum} we investigate the spectrum of the pushforward operator and construct the measure of maximal entropy. Section \ref{sec_anisotropic_de_rham_cohomology_and_spectrum} contains the connection between the eigenvalues of the transfer operator, that are related to the Ruelle-Pollicott resonances, and the eigenvalues of the action $f_\#$ on de Rham cohomology. We point out that Section \ref{sec_anisotropic_de_rham_cohomology_and_spectrum} inspired the techniques that R. Trevi\~no exploited in the recent paper \cite{trevino}. In addition, there is an analogy between our tools and the techniques of F.Faure, S. Gou\"ezel and E. Lanneau in \cite{faure_gouezel_lanneau}, where they found the complete set of Ruelle-Pollicott resonances for linear pseudo-Anosov maps on half- translation surfaces. Our paper, in particular Section \ref{sec_anisotropic_de_rham_cohomology_and_spectrum},  can be seen as a weak generalization of \cite{faure_gouezel_lanneau} to nonlinear Anosov diffeomorphisms, and it is reasonable to think that similar ideas can be used to treat the same problems for nonlinear pseudo-Anosov maps. 
	Section \ref{sec_proofs} contains the proofs of the main results, stated in Section \ref{sec_statement}. Finally, Appendix \ref{technical_results} collects the proofs of some technical lemmas that we use in several part of the paper.
	\begin{notation*}
		Throughout this paper we denote by $C$ a generic constant that can depend on the manifold, the dynamics or the atlas of $M$. We point out that $C$ may change also inside the same equation. If we want to emphasize the dependence of $C$ on a parameter $a$, we write $C_a.$ These constants as well may change at any occurrence inside a single equation.
	\end{notation*}
	\paragraph{Acknowledgments.} This paper comes out of the author's PhD thesis that is  available on his personal website. I warmly express gratitude to Marco Lenci, Giovanni Forni and Carlangelo Liverani for their suggestions and clarifications, and for supporting me during the preparation of this manuscript. I also thank Sébastien Gouëzel and Rodrigo Treviño for reporting some mistakes of a first version. I finally thank the anonymous referee for many comments and suggestions. Part of the present research was carried out at the Department of Mathematics of the University of Maryland, with the  support of the Marie Sklodowska-Curie grant agreement No 777822 (‘GHAIA’), and the PRIN Grant 2017S35EHN, MUR, Italy. It is also part of the author's activity 
	within the UMI Group \emph{DinAmicI}.
	\section{Setup and results}
	\label{section_Setup_results}
	In this section we recall the basic concept of Anosov diffeomorphisms, we set assumptions and we finally state the main results  of this paper. 
	
	Let $M$ be an orientable, closed, $\dim(M)-$dimensional Riemannian manifold, endowed with the metric $g$. Let $f\colon M\rightarrow M$ be a $C^\infty$ Anosov diffeomorphism.  In other terms, there exist a constant $\lambda>1$ and two continuous  families of cones $\mathcal{C}^s$  and $\mathcal{C}^u$ that satisfy the following two properties: 
	\begin{enumerate}
		\item 	$d_xf^{-1}\mathcal{C}^{s}_x\subseteq \text{int}(\mathcal{C}^{s}_{f^{-1}(x)})\cup\{0\},\ d_xf\mathcal{C}^{u}_x\subseteq \text{int}(\mathcal{C}^{u}_{f(x)})\cup\{0\}$, for any $x\in M$;
		\item $	\norm{d_xf^{-n}v} \geq \lambda^n\norm{v} ,\ \text{for }v\in \mathcal{C}_x^{s},\text{ and }
		\norm{d_xf^{n}v} \geq \lambda^n \norm{v} ,\ \text{for }v\in \mathcal{C}_x^{u}.$
	\end{enumerate}
	We denote by $d_s$ (resp. $d_u$) the dimension of the \textit{stable} (resp. \textit{unstable}) cone $\mathcal{C}^s$ (resp. $\mathcal{C}^u$).
	We finally assume that $f$ is topologically transitive.\footnote{In this setting, topological transitivity is equivalent to the existence of a dense orbit. It is conjectured that every Anosov diffeomorphism on a connected manifold is topologically transitive. See \cite{micena} for a recent account to this problem and \cite{micena2} for a new result relating the measure of maximal entropy and topological transitivity.}
	
	Before stating our main results, we recall some basic notions of de Rham cohomology (see \cite[Chapter 15]{lee} for a complete discussion of the topic). Let us denote by $\Omega^l(M)$ the space of $l$-differential forms on $M$. $\omega\in\Omega^l(M)$  is \textit{closed} if $d\omega=0,$ where $d\colon\Omega^l(M)\rightarrow\Omega^{l+1}(M)$ denotes the exterior derivative on differential forms. A differential form $\omega\in\Omega^l(M)$ is \textit{exact} if $\omega=du$, for some $u\in\Omega^{l-1}(M)$. As a consequence of the relation $d\circ d=0,$ one can define the $l$-de Rham cohomology group $H^l_{dR}(M)$ as the quotient of closed $l$-forms w.r.t. exact $l$-forms. Since $M$ is compact, $H^l_{dR}(M)$  is a finite-dimensional vector space.  Every diffeomorphism $f$ on $M$ induces an action $f_*$ by pushforward on the spaces of differential forms. Since $f_*$ commutes with $d$, the pushforward descends to a linear map $f_\#$ on every finite-dimensional de Rham cohomology group. We denote by $\sigma(f_\#|_{H^l_{dR}(M)})$ the spectrum of $f_\#\colon H^l_{dR}(M)\rightarrow H^l_{dR}(M).$

	\subsection{Statement of the results}
	\label{sec_statement}
	We can now state our main results.
\begin{theorem}
	\label{thm_strong}
	Let M be an orientable, closed, $\dim(M)-$dimensional Riemannian manifold and let $f\colon M\rightarrow M$ be a $C^\infty$, topologically transitive, Anosov diffeomorphism on $M$ with expanding, resp. contracting, factor $\lambda>1$, resp. $\lambda^{-1}<1$.  Let $\{\Lambda_1=e^{h_\top},\Lambda_2,\dots,\Lambda_m\}$ be the set of eigenvalues of the induced action $f_{\#}$ on $H_{dR}^{d_s}(M)$ such that $|\Lambda_i|>\lambda^{-1}e^{h_\top}$, for any $i=1,\dots,m$. Then, there exist $r\in\N$, $C>0$ and, for  any $i=2,\dots, m$, there exist $N_i\in\N$ and finite-rank bilinear forms $\{c_{\Lambda_i,k}(\cdot,\cdot)\}_{k=0}^{N_i-1},$ such that 
	\begin{equation*}
		\begin{split}
			\left|\int_{M}\phi(\psi\circ f^n) d\mu_{BM}-\right.&\int_{M}\phi d\mu_{BM}\int_{M}\psi d\mu_{BM}-\\-&\left.\sum_{i=2}^{m}\sum_{k=0}^{N_i-1}(\Lambda_ie^{-h_\top})^nn^k c_{\Lambda_i,k}(\phi,\psi)\right|\leq C\lambda^{-n} \norm{\phi}_{C^r}\norm{\psi}_{C^r},
		\end{split}
	\end{equation*}
	for every couple of observables $\phi,\psi\in C^\infty(M).$
\end{theorem}
Theorem \ref{thm_strong} shows that Ruelle-Pollicott resonances larger than $\lambda^{-1}$, as well as their multiplicities, are completely determined by the action induced by $f$ on de Rham cohomology.  A straightforward consequence of Theorem \ref{thm_strong} is the following corollary.
\begin{corollary}
	\label{cor_gen}
	Let $M$ be a compact $\dim(M)$-dimensional Riemannian manifold and let $f\colon M\rightarrow M$ be a $C^\infty$ topologically transitive Anosov diffeomorphism on $M$ with expanding, resp. contracting, factor $\lambda>1$, resp. $\lambda^{-1}<1$.  Let $\nu=\max\{|\Lambda_2|,\lambda^{-1}e^{h_\top}\}$, where $\Lambda_2$ is the second largest eigenvalue of the induced action $f_{\#}$ on the $d_s$-de Rham cohomology group $H_{dR}^{d_s}(M)$. Then, there exist $r\in\N$ and $C>0$ such that, for every couple of observables $\phi,\psi\in C^\infty(M),$
	\begin{equation*}
		\left|\int_{M}\phi(\psi\circ f^n) d\mu_{BM}-\int_{M}\phi d\mu_{BM}\int_{M}\psi d\mu_{BM}\right|\leq C\nu^ne^{-nh_\top} \norm{\phi}_{C^r}\norm{\psi}_{C^r}.
	\end{equation*}
\end{corollary}

The spectrum of $f_\#$ on de Rham cohomology can be explicitly computed when $f$ is a hyperbolic automorphism of the torus.  Since every Anosov diffeomorphism of the torus \cite{franks_anosov_diffeonorphisms_of_tori}, as well as every codimension 1 Anosov diffeomorphism\footnote{Codimension 1 means that the dimension of the stable (or the unstable) manifolds is 1.} \cite{franks,newhouse}, is topologically conjugated to a hyperbolic automorphism of the torus, we obtain the following corollary. 
\begin{corollary}
	\label{cor_main}
	Let M be an orientable, closed, $\dim(M)-$dimensional Riemannian manifold and let $f\colon M\rightarrow M$ be a $C^\infty$ transitive Anosov diffeomorphism on $M$ with expanding, resp. contracting, factor $\lambda>1$, resp. $\lambda^{-1}<1$. Assume that $f$ is topologically conjugate to a hyperbolic automorphism of the torus $\T^{\dim(M)}.$ Then, there exist $r\in\N$ and $C>0$ such that, for every couple of observables $\phi,\psi\in C^\infty(M),$
	\begin{equation*}
		\left|\int_{M}\phi(\psi\circ f^n) d\mu_{BM}-\int_{M}\phi d\mu_{BM}\int_{M}\psi d\mu_{BM}\right|\leq C\lambda^{-n} \norm{\phi}_{C^r}\norm{\psi}_{C^r}.
	\end{equation*}
\end{corollary}
As a consequence of this corollary, there are no Ruelle-Pollicott resonances in the set $\{z\in\C|\ \lambda^{-1}<|z|=1\}$, except for the trivial resonance $\{1\}$. Corollary \ref{cor_main} can be seen as an generalization of Theorem \ref{thm_liverani} to higher-dimensional tori.

\section{Anisotropic Banach spaces of currents}
\label{sec_anisotropic_banach_spaces_of_currents}
This section is devoted to the construction of a family of anisotropic Banach spaces which make the pushforward  operator $f_*$ quasi-compact, as proved in Section \ref{sec_quasi-compactnes}. We recall that, by definition, $f_*$ is quasi-compact if the essential spectral radius\footnote{Recall that there are different definition of essential spectrum. Here we adopt the following Browder's definition \cite{browder}. The essential spectrum of $f_*$ is the subset $\sigma_\ess(f_*)$ of points $\lambda$ in the spectrum $\sigma(f_*)$ such that at least one of the following holds: 1) $\text{Range}(\lambda-f_*)$ is not closed; 2)$\cup_{r\geq 1}\ker(\lambda-f_*)^r$ is infinite dimensional; 3) $\lambda$ is a limit point of $\sigma(f_*)\setminus\{\lambda\}$.} $\rho_{ess}(f_*)$ is strictly smaller than the spectral radius $\rho(f_*)$. Consequently, we can write  $f_*=C+R,$ where $C$  is a compact linear operator of finite rank and $R$ is another linear operator such that $\rho(R)<\rho(C)$. The adjective  \textit{anisotropic} means that the elements belonging to our spaces have different behaviors along stable and unstable manifolds. In particular, our objects behave as smooth differential forms in the unstable direction, while they behave as distributional differential forms, actually currents, in the stable direction (see Definition \ref{def_norm} and Lemma \ref{lemma_correnti}). The spaces we define in this section recall the ones introduced in
\cite{gouezel-liverani, giulietti-liverani-pollicott}(see also \cite{bahsoun-liverani,butterley-kiamari-liverani,butterley-liverani,giulietti-liverani,gouezel-liverani08}). In particular, Section \ref{section_leaves} is focused on the construction of a suitable set of stable leaves and it refers to \cite[Section 3]{gouezel-liverani}. Then, Section  \ref{Section_Norms_and_Banach_spaces} is a simplified version of \cite[Section 3.2]{giulietti-liverani-pollicott}, since we work with Anosov diffeomorphisms, while the cited paper treats Anosov flows. For a recent survey of these constructions we also refer the reader to \cite{demers-kiamari-liverani}. 
\subsection{Local charts and the family  of stable leaves}
\label{section_leaves}
We now introduce local coordinates on the manifold, conformed to the hyperbolic dynamics. Given a small constant  $\rho>0,$  there exists a smooth atlas $(U_i,\psi_i)_{i=1}^m,$
such that\footnote{We denote by $\mathcal{B}_t(p,r)$ the ball centered in $p$ of radius $r$ in $\mathbb{R}^t.$ When the subscript $t$ is not specified we mean $t=\dim(M)$.}
$$\overline{\mathcal{B}_{d_s}(0,3\rho)\times\mathcal{B}_{d_u}(0,3\rho)}\subseteq U_i $$ and the maps $\psi_i\colon U_i\rightarrow M$ satisfy the following properties:\footnote{One may use the exponential map to construct such charts around every point of the manifold. A finite covering of $M$ gives the atlas we need.}
\begin{enumerate}
	\item $\bigcup_{i=1}^m\psi_i(\mathcal{B}(0,\rho))=M;$
	\item  $d_0\psi_i(\mathbb{R}^{d_s}\times\{0\})=E_{p_i}^s$ and $d_0\psi_i(\{0\}\times\mathbb{R}^{d_u})=E_{p_i}^u,$ where $p_i=\psi_i(0);$   
	\item For every $x\in U_i$, let $$\zeta^s_{x,i}=\{v+w\in T_xU_i : v\in \mathbb{R}^{d_s}\times\{0\},\ w\in\{0\}\times\mathbb{R}^{d_u},\ \norm{w}\leq\norm{v}\},$$
	and let
	$$\zeta^u_{x,i}=\{v+w\in T_xU_i: v\in \mathbb{R}^{d_s}\times\{0\},\ w\in\{0\}\times\mathbb{R}^{d_u},\ \norm{v}\leq\norm{w}\}.$$
	Choosing $\rho>0$ small enough, we require that the Euclidean stable/unstable cones $\zeta_{x,i}^{s/u}$ satisfy 
	$$\mathcal{C}^s_{\psi_i(x)}\subseteq d_x\psi_i\zeta^s_{x,i},\ \ \mathcal{C}^u_{\psi_i(x)}\subseteq d_x\psi_i\zeta^u_{x,i},$$
	$$d_{\psi_i(x)}f^{-1}(\mathcal{C}_{\psi_i(x)}^s\setminus\{0\})\subseteq d_{\psi_i(x)}f^{-1}d_x\psi_i(\zeta^s_{x,i}\setminus\{0\})\subseteq\text{int}(\mathcal{C}^s_{f^{-1}\circ\psi_i(x)}),$$
	$$d_{\psi_i(x)}f(\mathcal{C}_{\psi_i(x)}^u\setminus\{0\})\subseteq d_{\psi_i(x)}fd_x\psi_i(\zeta^u_{x,i}\setminus\{0\})\subseteq\text{int}(\mathcal{C}^u_{f\circ\psi_i(x)}),$$
	where $\mathcal{C}^{s/u}=\mathcal{C}^{s/u,\alpha}$ are the stable/unstable bundles of $f$.
	As a consequence, $\zeta_{x,i}^{s/u}$ are \textit{Euclidean} invariant cones for the dynamics and satisfy the Anosov condition in local charts.
\end{enumerate}
\begin{remark} From now on, for any $r\in\N$, $U\subseteq\R^d$ and any Banach algebra\footnote{Recall that a Banach space $(\mathcal{B},\norm{\cdot})$ is  Banach algebra if it is a linear algebra and $\norm{xy}\leq\norm{x}\norm{y}$, whenever $x,y\in\mathcal{B}.$} $\mathcal{B},$ we fix a norm on $C^r(U,\mathcal{B})$ that makes it a Banach algebra.\footnote{For instance, the norms $\norm{\phi}_{C^0}=\sup_{x\in U}\norm{\phi(x)}_{\mathcal{B}},$ $\norm{\phi}_{C^r}=\sum_{i=0}^r\binom{r}{i}\sup_{|\beta|=i}\norm{\de^\beta\phi}_{C^0},$ for some $a\geq1$, make $C^r(U,\mathcal{B})$ a Banach algebra.}
We also point out that, in our setting, $r$ will be chosen as big as we want, because our Anosov map is $C^{\infty}.$ 
Finally, given a linear map $L\colon C^r(U,\mathcal{B})\rightarrow C^r(U,\mathcal{B}),$ we define the usual operator norm
$$\norm{L}_{(C^r)^\star}=\sup_{\substack{\phi\in C^r\\\norm{\phi}_{C^r}\leq1}}\norm{L\phi}_{C^r}.$$
\end{remark}

We are ready to describe the set of stable leaves that we are going to use to define our anisotropic Banach spaces. We remark that by \textit{stable leaf} we do not mean a piece of stable manifold, but instead a small piece of manifold whose tangent space belongs to the stable cone bundle $\mathcal{C}^s$. This permits to overcome the problem of regularity of foliations that, in general, are only H\"older continuous. 

We firstly consider the following set of stable graphs in $\mathbb{R}^{\dim(M)}$
\begin{equation}
	\label{eq_stable_graphs}
	\begin{split}
		\mathcal{F}=\{&F\in C^r(\mathcal{B}_{d_s}(0,2\rho);\mathbb{R}^{d_u})|\ F(0)=0,\\\
		&\norm{F}_{C^0(\mathcal{B}_{d_s}(0,2\rho))}\leq2\rho,\ \norm{dF}_{(C^r(\mathcal{B}_{d_s}(0,2\rho)))^\star}\leq 1\}.
	\end{split}
\end{equation}
Given $F\in\mathcal{F}$ and a point $x\in\mathcal{B}(0,\rho),$ let $G_{x,F}$ be the graph of $F$ in $\mathbb{R}^{\dim(M)}$ centered at $x$, namely $G_{x,F}(\mathcal{B}_{d_s}(0,2\rho))=x+(y,F(y))_{y\in\mathcal{B}_{d_s}(0,2\rho)}.$ Notice that the tangent space to the graph $G_{x,F}$ belongs to the Euclidean stable cone $\zeta^s_{x,i}$. Finally, we define the set of full admissible leaves
$$\widetilde{\Sigma}=\{\psi_i\circ G_{x,F}(\mathcal{B}_{d_s}(0,2\rho))|\ x\in\mathcal{B}(0,\rho),\  F\in\mathcal{F},\ i=1,\dots,m\},$$
and the set of admissible leaves
$$\Sigma=\{\psi_i\circ G_{x,F}(\mathcal{B}_{d_s}(0,\rho))|\ x\in\mathcal{B}(0,\rho),\  F\in\mathcal{F},\ i=1,\dots,m\}.$$
Observe that, for any admissible leaf $W\in\Sigma$, there is a full admissible leaf $\widetilde{W}\in\widetilde{\Sigma}$ containing $W$. Moreover, notice that the sets of leaves $\Sigma$ and $\widetilde{\Sigma}$ are well-defined. Indeed, the graph of $F\in\mathcal{F}$ is included in $\mathcal{B}_{d_s}(0,2\rho)\times\mathcal{B}_{d_u}(0,2\rho)$ and, since $x\in\mathcal{B}(0,\rho),$ 
$G_{x,F}(\mathcal{B}_{d_s}(0,2\rho))\subset\mathcal{B}_{d_s}(0,3\rho)\times\mathcal{B}_{d_u}(0,3\rho)\subseteq U_i,$ for all $i=1,\dots,m.$

The importance of the set $\Sigma$ is given by the following lemma. 
\begin{lemma}
	\label{lemma_foglie}
	There exist $n_0\in\N$ and $\rho>0$ small enough such that for each full admissible leaf $\widetilde{W}$, with corresponding admissible leaf $W,$ and for each $n\geq n_0,$ there are a finite number (depending only on $n$) of admissible leaves $W_1,\dots,W_l\in\Sigma$ such that 
	\begin{enumerate}
		\item $f^{-n}(W)\subseteq \cup_{i=1}^lW_i\subseteq f^{-n}(\widetilde{W})$;
		\item the leaves $W_1,\dots,W_l$ admit a uniformly finite (in $W$ and $n$) number of overlaps.;
		\item there exists a constant $C_\rho$ and a  $C^{r}$ partition of unity $\eta_1,\dots,\eta_l$ subordinated to $\{W_1,\dots,W_l\}$ on $f^{-n}(W),$ such that $\norm{\eta_i}_{C^r}\leq C_\rho,$ for any $i=1,\dots,l.$
	\end{enumerate}
\end{lemma}
We refer the reader to \cite[Lemma 3.3]{gouezel-liverani} for the proof of Lemma \ref{lemma_foglie}. We remark that the authors of \cite{gouezel-liverani} proved a slightly stronger result for what they call $\gamma$-admissible leaves, which are used  to apply perturbation arguments. Our Lemma \ref{lemma_foglie} corresponds to $\gamma=1$.

\subsection{Norms and anisotropic Banach spaces}
\label{Section_Norms_and_Banach_spaces}
Let us introduce our family of anisotropic Banach spaces. In particular, we are going to define anisotropic norms on spaces of differential forms. Then, we obtain the suitable family of anisotropic Banach spaces by completing the spaces of differential forms with respect to these norms. 

We denote by $\Omega^l(M)$, for each $l=0,\dots,\dim(M),$ the space of complex smooth differential forms on $M,$ namely the set of $C^\infty$ sections of the $l$-exterior algebra of the cotangent bundle $T^*M$ over $M,$ with values in $\C.$
Given an admissible leaf $W=\psi_i\circ G_{x,F}(\mathcal{B}_{d_s}(0,\rho))\in\Sigma$, $s\geq0$ and $l\in\{0,\dots,\dim(M)\}$, we denote by $\Gamma_0^{l,s}(W)$ the space of complex $C^s$ sections of the fiber bundle over $W$, with the fiber space  $\wedge^l(T^*M)$ and compact support. In other words, we may think the elements $\Gamma_0^{l,s}(W)$ as $l-$differential forms of class $C^s,$ defined on $W$ and vanishing in a neighborhood of $\de W.$ This is exactly the space introduced in \cite[Section 3]{giulietti-liverani-pollicott} and, in the definition of the norm, its elements have the role of ``test forms". Let $\mathcal{V}^s(W)$ be the space of $C^s$ vector fields defined in a $\dim(M)$-dimensional neighborhood $U(W)\subseteq U_i$ of $W.$ 

At some point in the proofs we will need to express forms and vector fields in local coordinates. Let us set the notation. Given the atlas $\{U_i,\psi_i\}_{i=1}^m$, let $\{\chi_i\}_{i=1}^m$ be a smooth partition of unity subordinate to the atlas, such that $\chi_i|_{\psi_i(\mathcal{B}(0,3\rho))}=1.$
We denote by $\de_{r_1},\ \de_{r_2},\dots,\ \de_{r_{\dim(M)}}$ a basis for the vector fields on $U_i$ such that $\de_{x_1}:=\psi_i^*(\de_{r_1}),\ \de_{x_2}:=\psi_i^*(\de_{r_2}),\dots,\ \de_{x_{d_s}}:=\psi_i^*(\de_{r_{d_s}})\in\mathcal{C}^s$ and $\de_{x_{d_s+1}}:=\psi_i^*(\de_{r_{d_s+1}}),\ \de_{x_{d_s+2}}:=\psi_i^*(\de_{r_{d_s+2}}),\dots,\ \de_{x_{d_s+d_u}}:=\psi_i^*(\de_{r_{d_s+d_u}})\in\mathcal{C}^u.$ Without loss of generality, we may suppose that this is an orthonormal basis of vector fields, otherwise one can apply the Gram-Schmidt procedure without essentially affecting the forthcoming arguments.
Finally, let $dx_1,\ dx_2,\ \dots,\ dx_d$ be the dual basis of $\de_{x_1},\ \de_{x_2},\dots,\ \de_{x_d}$, i.e., the corresponding basis of differential forms on $\psi_i(U_i).$
Let $\mathcal{J}_l=\{\overline{j}=(j_1,\dots,j_l)\in\{1,\dots,d\}^l|j_1<j_2<\dots<j_l\}$ the set of ordered $l$-multi-indexes. We adopt the following notation for fields and differential forms: 
$\de_{x_{\overline{i}}}:=\de_{x_{i_1}}\wedge\dots\wedge\de_{x_{i_l}},\ dx_{\overline{i}}:=dx_{i_1}\wedge\dots\wedge dx_{i_l},$ whence $dx_{\overline{i}}(\de_{x_{\overline{i}}})=l!.$

We can decompose every form $h\in\Omega^l(M)$ as $h=\sum_{i=1}^m h_i,$ where $h_i=h\chi_i\in\Omega_0^l(\psi_i(U_i)),$ i.e., $h_i$ is a smooth differential form on $\psi_i(U_i)$ with compact support. Moreover, using the local basis, we can write every $h_i$ in coordinates as  $h_i=\sum_{\overline{j}\in\mathcal{J}_l}h_i^{\overline{j}}dx_{\overline{j}}.$ We define the $C^s$ norm of $h\in\Omega^l(M)$ as
\begin{equation}
	\label{norma_forme}
	\Vert h\Vert_{C^s(M)}=\sup_{i=1,\dots,m}\sup_{\overline{j}\in\mathcal{J}_l}\Vert h_i^{\overline{j}}\Vert_{C^s(\psi_i(U_i))}=\sup_{i=1,\dots,m}\sup_{\overline{j}\in\mathcal{J}_l}\Vert h_i^{\overline{j}}\circ\psi_i^{-1}\Vert_{C^s(U_i)}\ .
\end{equation}
Similarly, given $\phi\in\Gamma_0^{s,l}(W),$ we can write $$\phi=\sum_{i=1}^m\phi\chi_i=\sum_{i=1}^m\phi_i=\sum_{i=1}^m\sum_{\overline{j}\in\mathcal{J}_l}\phi_i^{\overline{j}}dx_{\overline{j}}$$
and
$$\Vert\phi\Vert_{\Gamma_0^{l,s}(W)}=\sup_{i=1,\dots,m}\sup_{\overline{j}\in\mathcal{J}_l}\Vert \phi_i^{\overline{j}}\Vert_{C^s(\psi_i\circ G_{x,F}(\mathcal{B}_{d_s}(0,\delta)))}=\sup_{i=1,\dots,m}\sup_{\overline{j}\in\mathcal{J}_l}\Vert \phi_i^{\overline{j}}\circ\psi_i^{-1}\Vert_{C^s(G_{x,F}(\mathcal{B}_{d_s}(0,\delta)))}.$$

The last ingredient we need is a scalar product for differential forms. We point out that such a scalar product, and consequently the induced norm, depends on the metric. On the other hand, the Banach space we are establishing is independent of the metric. Let us consider the metric $g$ introduced in Section \ref{section_Setup_results} and the induced volume form $\omega_0\in\Omega^{\dim(M)}(M)$. The non-degeneracy condition of $g$ induces an isomorphism $\xi$ between smooth vector fields $\mathcal{V}(M)$ and smooth 1-forms $\Omega^1(M)$ such $g(v,\cdot)=\xi(v)(\cdot)$ for every $v\in\mathcal{V}(M).$ A pointwise scalar product between 1-forms $\omega_1,\ \omega_2\in\Omega^1(M)$ is $\langle\omega_1,\omega_2\rangle=g(\xi^{-1}(\omega_1),\xi^{-1}(\omega_2)).$ Similarly, for $\{\omega_{i,j}\}_{i=1,2;\ j=1,\dots,l}\subseteq\Omega^1(M),$ a pointwise scalar product between $l$-forms  is given by
\begin{equation*}
	\langle\omega_{1,1}\wedge\dots\wedge\omega_{1,l},\omega_{2,1}\wedge\dots\wedge\omega_{2,l}\rangle=\det\left(\begin{matrix}
		\langle\omega_{1,1},\omega_{2,1}\rangle&\dots&\langle\omega_{1,l},\omega_{2,1}\rangle\\
		\vdots&\ddots&\vdots\\
		\langle\omega_{1,1},\omega_{2,l}\rangle&\dots&\langle\omega_{1,l},\omega_{2,l}\rangle\\
	\end{matrix}\right).
\end{equation*}
Finally, the scalar product between $l$-forms is the integral of the pointwise scalar product, i.e.,
$$\left(\omega_1,\omega_2\right)=\int_{M}\langle\omega_1,\omega_2\rangle\omega_0,\hspace{1cm}\omega_1,\omega_2\in\Omega^l(M)\ .$$

The following definition introduces the anisotropic norms.
\begin{definition}
	\label{def_norm}
	Let $p,q\in\mathbb{N}$. Given $h\in\Omega^l(M)$, we define a seminorm $$\Vert h\Vert_{p,q,l}^-=\sup_{W\in\Sigma}\ \ 
	\sup_{\substack{v_1,\dots,v_p\in\mathcal{V}^{p+q}(U(W))\\\Vert v_k\Vert_{C^{p+q}(U(W))}\leq1}}\ \ \sup_{\substack{\phi\in\Gamma_0^{p+q,l}(W),\\\Vert\phi\Vert_{\Gamma_0^{p+q,l}(W)}\leq 1}}\left|\int_W\langle \phi,L_{v_1}\dots L_{v_p}h\rangle\omega_W\right|,$$
	where $L_{v_i}$ is the Lie derivative of the $l$-form $h$ w.r.t.\ the vector field $v_i$ and $\omega_W$ is the measure induced by $g$ on $W.$
	We define, for every $h\in\Omega^l(M)$, the norm 
	$$\norm{h}_{p,q,l}=\max_{t\leq p}\Vert h\Vert_{t,q,l}^-\ .$$
	Finally, we denote by $\mathcal{B}^{p,q,l}=\overline{\Omega^l(M)}^{\Vert\cdot\Vert_{p,q,l}}$ the closure of the space of $l$-forms w.r.t.\ such a norm.
\end{definition}
\begin{remark}
	Notice that the Banach spaces $\mathcal{B}^{p,q,l}$ are similar to the spaces studied in \cite[Section 5]{butterley-kiamari-liverani} for Anosov map on the 2-torus. The authors of \cite{butterley-kiamari-liverani} were in turn inspired by \cite{giulietti-liverani-pollicott}, where these spaces are used to analyze dynamical zeta functions of Anosov flows. The following Lemma \ref{lemma_correnti}, whose proof can be found in \cite[Lemma 3.10]{giulietti-liverani-pollicott}, shows that we can think these spaces  as subspaces of currents, i.e., the continuous dual space of differential forms. For an overview of currents' properties we refer the reader to de Rham's book \cite{de_rham_book}.
\end{remark}
Given $p,q\in\N,$ let us denote by $\Omega_{p+q}^l(M)$ the space of $C^{p+q}$ $l$-forms equipped with the $C^{p+q}$ norm as defined in \eqref{norma_forme}. Let $(\Omega_{p+q}^l(M))^\star$ be its dual space with the weak-$^\star$topology, i.e., the space of currents of dimension $l$, degree $\dim(M)-l$ and regularity $C^{p+q}.$
\begin{lemma}
	\label{lemma_correnti}
	The space $\mathcal{B}^{p,q,l}$ can be identified with a subspace of the space of currents of dimension $\dim(M)-l$, degree $l$ and regularity $C^{p+q}$ on the manifold $M;$ i.e., there exists an injective bounded linear operator $\iota\colon\mathcal{B}^{p,q,l}\rightarrow(\Omega_{p+q}^l(M))^\star.$
\end{lemma}
\begin{remark}
	\label{rem_properties}
	We conclude this section with a short overview of the main properties of these anisotropic Banach spaces. Firstly, the following inequalities are trivial consequences of the definition: 
	$\norm{h}_{p,q,l}\leq\norm{h}_{p+1,q-1,l}$ and $\norm{h}_{p,q,l}\leq C\norm{h}_{C^p}.$ It follows that for all $p,q\in\N$, such that $p>0$,  $\mathcal{B}^{p+1,q-1,l}\subseteq\mathcal{B}^{p,q,l}$, while the second inequality gives  $\Omega^l_p(M)\subseteq\mathcal{B}^{p,q,l}.$ Moreover, by definition, the space of $C^\infty$ differential $l$-forms, as well as the space of $C^p$ $l$-forms,  is dense in $\mathcal{B}^{p,q,l}.$ 
	Furthermore, the following inequality $$\sup_{\substack{v_1,\dots,v_p\in\mathcal{V}^{p+q}(U(W))\\\Vert v_k\Vert_{C^{p+q}(U(W))}\leq1}}\norm{L_{v_1}\dots L_{v_p}h}^-_{0,p+q,l}\leq\norm{h}^-_{p,q,l}\ $$
	is an immediate consequence of the definition, and it will be used later in this paper.  
	Finally, Lemma \ref{lemma_correnti} gives the inclusion $\mathcal{B}^{p,q,l}\subset (\Omega_{p+q}^l(M))^*$ of the anisotropic Banach spaces into the space of currents, but notice that  $\mathcal{B}^{p,q,l}$ is not dense in $(\Omega_{p+q}^l(M))^*.$
	
\end{remark}

\subsection{Quasi-compactness of the transfer operator}
\label{sec_quasi-compactnes}
This section is devoted to prove that the action of the pushforward operator $f_*$, initially defined on differential forms\footnote{Recall that the pullback of a differential form $h\in\Omega^l(M)$ by $f$ is $(f^*h)_x(v_1,\dots,v_l)=h_{f(x)}(d_xfv_1,\dots, d_xf v_l),$ where $v_1,\dots, v_l$ are smooth vector fields on $M$. Since $f$ is a diffeomorphism, the pushforward is defined as $f_*=(f^{-1})^*.$}, extends to a bounded linear operator on  $\mathcal{B}^{p,q,l}.$ Moreover, we show that $f_*\colon\mathcal{B}^{p,q,l}\rightarrow\mathcal{B}^{p,q,l}$ is quasi-compact.
Notice that $f_*$, that is the transfer operator we are dealing with in this paper, encloses  the classical Perron-Frobenius transfer operator $$\mathcal{L}h(x)=\frac{h\circ f^{-1}(x)}{|\det d_{f^{-1}(x)}f|}.$$
Indeed, for  every $\omega\in\Omega^{\dim(M)}(M)$, can write $\omega=h\omega_0,$ where $h\in C^\infty(M)$ and $\omega_0$ is the volume form induced by the metric $g.$ Then, a straightforward computation shows that 
\begin{equation}
	\label{eq_transfer_op}
	f_*\omega=(\mathcal{L}h)\omega_0,
\end{equation} i.e., the spectral properties of the Perron-Frobenius operator can be obtained by studying the action of $f_*$ on $\dim(M)$-differential forms or, more generally, on $\mathcal{B}^{p,q,\dim(M)}$.  It is well known that the spectral properties of the  Perron-Frobenius operator $\mathcal{L}$, acting on $\mathcal{B}^{p,q,0},$ give information about the statistical properties of smooth Anosov diffeomorphisms with respect to the SRB measure (see \cite{gouezel-liverani}). In fact, the space $\mathcal{B}^{p,q,0}$ contains the unique invariant density for the SRB measure, that is a distribution $h$ such that $\mathcal{L}h=h$.  In view of \eqref{eq_transfer_op}, the same holds for $f_*\colon\mathcal{B}^{p,q,\dim(M)}\rightarrow \mathcal{B}^{p,q,\dim(M)},$ with the difference that $\mathcal{B}^{p,q,\dim(M)}$ contains a current $\omega$ (of degree $\dim(M)$ and dimension 0) such that $f_*\omega=\omega$, and it is the SRB measure.   
In this paper, we are interested in the statistical properties of the measure of maximal entropy and it turns out that both $\mathcal{B}^{p,q,0}$ and $\mathcal{B}^{p,q,\dim(M)}$ are not the right spaces, because the Bowen-Margulis measure does not belong to them. On the other hand, we recall that the measure of maximal entropy is the unique invariant measure $\mu$ that can be decomposed as the tensor product of the marginals along the stable and the unstable foliation, i.e., $\mu=\mu^s\otimes\mu^u$, with the following scaling properties\footnote{Notice that the pushforward of a measure $m$ by $f$ is defined as $f_*m(B)=m(f^{-1}(B)),$ for any measurable set $B$. Hence, $f_*$ acts on the measure as $f^{-1}$ and for this reason the stable measure $\mu^s$ is expanded by $e^{h_\top}$, while the unstable measure is contracted by $e^{-h_\top}$. }: $f_*\mu^s=e^{h_\top}\mu^s$ and $f_*\mu^u=e^{-h_{\top}}\mu^u$. 
The previous argument suggests that $\mu^s$ should belong to a Banach space containing $d_s$-differential forms. Indeed, we are going to show in the subsequent sections that $\mu^s$ can be found in  $\mathcal{B}^{p,q,d_s}$ as the unique eigenvector of $f_*$ for the simple maximal eigenvalue $e^{h_\top}.$ By duality, $\mu^u$ belongs to the dual space of $\mathcal{B}^{p,q,d_s}$ and satisfies the eigenvalue equation $(f^{-1})_*\mu^u=e^{h_\top}\mu^u$ (see Section \ref{sec_peripheral_spectrum} for details).

The standard procedure to prove quasi-compactness of linear operators is based on the following classical result of Hennion \cite{hennion}. From now on, let us denote by $\rho(\mathcal{T})$ (resp. $\rho_{ess}(\mathcal{T})$) the spectral radius (resp. the essential spectral radius) of a bounded linear operator $\mathcal{T}.$
\begin{theorem}[\cite{hennion} Corollary 1]
	\label{thm_hennion}
	Let $\mathcal{T}$ be a bounded linear operator on the Banach space $(\mathcal{B},\Vert\cdot\Vert).$ Assume that there exists another norm $|\cdot|$ on $\mathcal{B}$ such that 
	\begin{enumerate}
		\item[(1)] the immersion $i\colon(\mathcal{B},\Vert\cdot\Vert)\rightarrow(\mathcal{B},|\cdot|)$ is bounded and compact;
		\item[(2)] for all $n\in\mathbb{N}$ there exist $R_n,\ r_n>0,$ with $\liminf_{n\rightarrow+\infty}r_n^{1/n}<\rho(T)$ such that 
		$$\Vert \mathcal{T}^nf\Vert\leq r_n\Vert f\Vert+R_n|f|;$$
	\end{enumerate}
	then the operator, acting on $(\mathcal{B},\norm{\cdot}),$ is quasi-compact with essential spectral radius $\rho_{ess}(\mathcal{T})\leq \liminf_{n\rightarrow+\infty}r_n^{1/n}.$ 
\end{theorem}
Next  lemma  gives property \textit{(1)} of Theorem \ref{thm_hennion} in our setting.
\begin{lemma}
	\label{inclusione_compatta}
	The inclusion $\iota\colon\mathcal{B}^{p,q,l}\rightarrow\mathcal{B}^{p-1,q+1,l}$ is compact for each $p\in\Z^+,$ $q\in\N$ and for all $l=1,\dots,\dim(M).$ 
\end{lemma}
The above Lemma \ref{inclusione_compatta} is proved in \cite[Lemma 2.1]{gouezel-liverani} for the anisotropic Banach space of distributions $\mathcal{B}^{p,q,0}.$ The extension to $\mathcal{B}^{p,q,l}$ for $l\geq 1$ follows exactly the same strategy and it does not require any new idea.

Property \textit{(2)} of Hennion's Theorem \ref{thm_hennion} is called a Lasota-Yorke type inequality. The following theorem states that such inequalities hold for the pushforward operator when acting on $\mathcal{B}^{p,q,l}$.
\begin{theorem}
	\label{thm_lasota-yorke}
	$f_*$ acts as a bounded linear operator on the spaces $\mathcal{B}^{p,q,l}$. In particular, for $p\geq0$  and for $0<l<\dim(M),$
	\small
	\begin{align}
		&\label{Lasota-Yorke_1}
		\Vert  f_*^nh\Vert _{p,q,0}\leq C\Vert  h\Vert _{p,q,0};\\
		&\label{Lasota-Yorke_2}
		\Vert  f_*^nh\Vert _{p,q,l}\leq C\lambda^{-|d_s-l|n}e^{h_\top n}\Vert  h\Vert _{p,q,l};\\
		&\label{Lasota-Yorke_5}\Vert  f_*^nh\Vert _{p,q,\dim(M)}\leq C\Vert  h\Vert _{p,q,\dim(M)}.
	\end{align}
	\normalsize
	Moreover, for $p>0$ and for $0<l<\dim(M),$
	\small
	\begin{align}
		&\label{Lasota-Yorke_3}
		\Vert  f_*^nh\Vert _{p,q,0}\leq C\lambda^{-n\min\{p,q\}}\Vert  h\Vert _{p,q,0}+C\Vert h\Vert_{p-1,q+1,0};\\
		&\label{Lasota-Yorke_4}
		\Vert  f_*^nh\Vert _{p,q,l}\leq C\lambda^{-n(\min\{p,q\}+|d_s-l|)}e^{h_\top n}\Vert h\Vert _{p,q,l}+C \lambda^{-n|d_s-l|}e^{nh_\top }\Vert h\Vert_{p-1,q+1,l};\\
		&\label{Lasota-Yorke_6}\Vert  f_*^nh\Vert _{p,q,\dim(M)}\leq C\lambda^{-n\min\{p,q\}}\Vert  h\Vert _{p,q,\dim(M)}+C\Vert h\Vert_{p-1,q+1,\dim(M)}.
	\end{align}
	\normalsize
\end{theorem}
\begin{remark}
	Notice that the proof \eqref{Lasota-Yorke_2} and \eqref{Lasota-Yorke_4} also work for $l=0$ and $l=\dim(M)$, but it gives weaker bounds than the other four inequalities. In fact, by using Ledrappier-Young entropy formula \cite[Theorem D.3.1]{brown}, 
	it is easy to check that $e^{h_\top }\lambda^{-{d_s}}\geq 1.$ The reasons behind this mismatch will be clear at the end of the proof, but, just to give an idea, the factor $\lambda^{-|d_s-l|}e^{h_\top }$ comes up from the estimate of the expansion of $l$-dimensional subspaces of the tangent bundle, under the action of the differential. Hence, there is no expansion on 0-dimensional subspaces when $l=0$ and, by duality \eqref{eq_transfer_op}, when $l=\dim(M)$. On the other hand, the bounds of \eqref{Lasota-Yorke_2} and \eqref{Lasota-Yorke_4} work for any other $l=1,\dots,\dim(M)-1$ and the factor $\lambda^{-|d_s-l|}e^{h_\top }$ is the best estimate of $l$-dimensional subspaces' growth.
	
	The expert reader could also notice a possible discrepancy with the classical Thermodynamic Formalism of uniformly hyperbolic maps (see \cite{bowen,ruelle_thermodynamic,sinai1,sinai2}). In fact, for $l=0$, the pushforward operator $f_*h= h\circ f^{-1}$ can be interpreted as a weighted transfer operator with potential $\phi=0$. Accordingly, the spectral radius should be $\exp(P(f,0))=\exp(h_\top),$ where $P(f,\phi)$ denotes the topological pressure. On the other hand, in this paper we prove that the spectral radius is 1. Similarly, the potential for the SRB measure in the classical setting is\footnote{Note that $E_u$ is in general only H\"older continuous, hence the same is true for $\phi$. This generates an important limitation, because forces to apply the transfer operator to low regular functions.} $\phi(x)=-\log|\det(df|_{E^u})|,$ while in  our setting the potential is $\widetilde{\phi}(x)=-\log|\det(df)|.$ To explain this mismatch we recall that the classical approach is based on the coding of the Anosov diffeomorphisms and the construction of a semiconjugacy with a two-sided subshift of finite type. An additional reduction allows to a one-sided subshift of finite time, which drops the dependence on the stable directions and allows to study the transfer operator on a space of H\"older functions. On the contrary, if we avoid to code the system, it is no longer true that the spectral radius is $\exp(P(f,\widetilde{\phi})),$ because of the contraction along the stable, but we can work on spaces of higher regularity, since $\widetilde{\phi}(x)=-\log|\det(df)|$ is smooth. So, why is the spectral radius of $\mathcal{L}\psi=\psi\circ f^{-1}|\det(df)|^{-1}$ on $\mathcal{B}^{p,q,0}$ bounded by 1? An heuristic motivation is given by the following argument. We can see $|\det(df)|^{-1}=|\det(df)|_{E^s}|^{-1}|\det(df)|_{E^u}|^{-1}.$ Then $\psi\circ f^{-1}$ behaves well under derivative in the unstable direction and $|\det(df)|_{E^u}|^{-1}$ is bounded by 1. The problematic term in the computation of the norm is $|\det(df)|_{E^s}|^{-1}$, which is large, but this is compensated by the negative (distributional) regularity $-q$ along the stable direction. 
	The same is true, when we look at the pushforward $f_* h=h\circ f^{-1}$ on $\mathcal{B}^{p,q,0}.$ Both the derivatives along the unstable and the ``antiderivatives" along the stable of behave well on $h\circ f^{-1}$ in a bounded way, hence the spectral radius is bounded by 1. 
\end{remark}

\proof
Let us prove \eqref{Lasota-Yorke_1} for $p=0$.
Let $h\in C^\infty(M)$, $W\in\Sigma,$ $\phi\in \Gamma_0^{q,0}(W)=C^q_0(W)$. The following integral can be split, using the partition of unity $(W_i,\eta_i)_{i\in\mathcal{I}}$ of Lemma \ref{lemma_foglie}, as 
$$\int_W \phi f_*^nh=\int_W \phi (h\circ f^{-n})=\int_{f^{-n}W}(\phi\circ f^n) h \lambda_n^s=\sum_{i}\int_{W_i}(\phi\circ f^n) h\lambda_n^s\eta_i,$$
where $\lambda_n^s(x)$ is the Jacobian of the change of variables.
As a consequence of the Distortion Lemma \ref{lemma_determinant}, the $C^q$-norm of $\lambda_n^s$ is bounded by $C|f^n(W_i)|,$ where $|f^n(W_i)|$ is the $\omega_{\widetilde{W}}$-measure of $f^n(W_i)$. Moreover, the $C^q$-norm of $\phi\circ f^n$ is bounded by the $C^q$-norm of $\phi,$ because the composition with $f^n$ reduces norms along stable manifolds. Thus, since we can assume that the $C^q$-norm of $\eta_i$ is uniformly bounded, the norm of $\phi\circ f^n \lambda_n^s\eta_i$ is bounded by $C|f^n(W_i)|\Vert \phi\Vert _{C_0^q(W)}$ , hence 
$$\left|\int_W \phi f_*^nh\right|\leq C\sum_{i}|f^n(W_i)|\Vert \phi\Vert _{C_0^q(W)}\Vert h\Vert_{0,q,l}.$$
By Lemma \ref{lemma_foglie}-\textit{2.}, we can choose an appropriate covering $\{W_i\}_{i\in\mathcal{I}}$ of $f^{-n}(W)$, with a uniformly bounded (in $W$ and $n$) number of overlaps. Accordingly, the sum in $i$ is bounded by a constant and, by density, we get \eqref{Lasota-Yorke_1} with $p=0$. 

Let us tackle together \eqref{Lasota-Yorke_1} for $p\geq0$ and  \eqref{Lasota-Yorke_3}. We firstly show that, for $h\in\Omega^l(M),$ 
\begin{equation}
	\label{eq_Lasota-Yorke_3}
	\norm{f_*^nh}^-_{t,q,0}\leq\begin{cases}
		\displaystyle
		C\lambda^{-nq}\norm{h}_{p,q,0}+C_n\norm{h}_{p-1,q+1,0}\ \ \text{if }t<p\\
		\displaystyle
		C\lambda^{-n\min\{p,q\}}\norm{h}_{p,q,0}+C_n\norm{h}_{p-1,q+1,0}\ \ \text{if }t=p
	\end{cases}
\end{equation}
so that
\begin{equation}
	\label{eq_dependent_n}
	\Vert f_*^nh\Vert_{p,q,0}=\max_{t\leq p}\Vert f_*^nh\Vert_{t,q,0}^-\leq C\lambda^{-n\min\{p,q\}}\Vert h\Vert_{p,q,0}+C_n\Vert h\Vert_{p-1,q+1,0}.
\end{equation}
The inequality \eqref{eq_dependent_n} does not prove \eqref{Lasota-Yorke_3}, as well as \eqref{Lasota-Yorke_1} for $p\geq0,$ because the coefficient of the weak norm does depend on $n$. 
We proceed by induction on $p$ to remove the dependence on that $n$. We have already shown \eqref{Lasota-Yorke_1} for $p=0.$ Assume that \eqref{Lasota-Yorke_1} holds for $p-1,$ then we prove \eqref{Lasota-Yorke_3} and \eqref{Lasota-Yorke_1} for $p$.
Let $\widetilde{\lambda}>\lambda$ be a constant for which the Anosov property continues to be true. Then \eqref{eq_dependent_n} implies 
$$\norm{ f_*^nh}_{p,q,0}\leq \widetilde{C}\widetilde{\lambda}^{-n\min\{p,q\}}\norm{h}_{p,q,0}+\widetilde{C}_n\norm{h}_{p-1,q+1,0}.$$
Let $N\in\Z^+$ be a positive integer such that $\widetilde{C}\widetilde{\lambda}^{-N\min\{p,q\}}\leq\lambda^{-N\min\{p,q\}}.$ Then, for every $n\in\Z^+,$ we write $n=Q_nN+R_n$ with $0\leq R_n<N$ and $Q_n\in\N.$ Consequently,
\begin{equation*}
	\begin{split}
		\norm{f_*^n h}_{p,q,0}\leq& \widetilde{C}\widetilde{\lambda}^{-N\min\{p,q\}}\norm{f_*^{n-N}h}_{p,q,0}+\widetilde{C}_N\norm{f_*^{n-N}h}_{p-1,q+1,0}\leq\\
		\leq&\widetilde{C}\widetilde{\lambda}^{-N\min\{p,q\}}\norm{f_*^{n-N}h}_{p,q,0}+\widetilde{C}_N\norm{h}_{p-1,q+1,0}\leq\\
		\leq&\left(\widetilde{C}\widetilde{\lambda}^{-N\min\{p,q\}}\right)^{Q_n}\norm{f_*^{R_n}h}_{p,q,0}+\widetilde{C_N}\sum_{i=0}^{Q_N-1}\left(\widetilde{C}\widetilde{\lambda}^{-N\min\{p,q\}}\right)^i\norm{h}_{p-1,q+1,0}\leq\\\leq&C\lambda^{-(NQ_n+R_n)\min\{p,q\}}\norm{h}_{p,q,0}+\widetilde{C_N}\frac{1-\lambda^{-NQ_n\min\{p,q\}}}{1-\lambda^{-N\min\{p,q\}}}\norm{h}_{p-1,q+1,0}+\\&+C_{R_n}\norm{h}_{p-1,q+1,0}\leq 
		C\lambda^{-n\min\{p,q\}}\norm{h}_{p,q,0} +C\norm{h}_{p-1,q+1,0},
	\end{split}
\end{equation*}
where we used \eqref{eq_dependent_n} to estimate $\norm{f_*^{R_n}h}_{p,q,0}$ and the inductive hypothesis to estimate $\norm{f_*^{n-iN}h}_{p-1,q+1,0}.$ This computation also gives \eqref{Lasota-Yorke_1} for $p\geq0$. By density, both \eqref{Lasota-Yorke_1} and \eqref{Lasota-Yorke_3} extend to $\mathcal{B}^{p,q,0}.$

We are left with the proof of \eqref{eq_dependent_n}.
Let $h\in\Omega^0(M)= C^\infty(M),$ $W\in\Sigma,$ $\phi\in \Gamma_0^{t+q,0}(W)=C_0^{t+q}(W),$ such that $\Vert\phi\Vert_{C^{t+q}_0(W)}\leq 1,$ and $v_1,\dots,v_t\in\mathcal{V}^{t+q}(U(W))$ with $\Vert v_i\Vert_{C^{t+q}}\leq 1.$
We compute 
$$\int_W\phi(x)L_{v_1}\dots L_{v_t}f_*^nh(x)=\int_W\phi(x)v_1\dots v_t(h\circ f^{-n})(x).$$
By linearity we can assume that $v_j=g_j\de_{x_j}$ where $\de_{x_j}$ is a coordinate vector field, hence  
\begin{equation*}
	\begin{split}
		&\int_W\phi(x)v_1\dots v_t(h\circ f^{-n})(x)=\int_W\phi(x)\prod_{j=1}^t (g_j(x)\de_{x_j})(h\circ f^{-n})(x)=\\&=\sum_i\int_{f^n(W_i)}\phi(x)\eta_i(f^{-n}(x))\prod_{j=1}^t g_j(x)\prod_{j=1}^t \de_{x_j}(h\circ f^{-n})(x)+C_n\Vert h\Vert_{p-1,q+1,0}
	\end{split}
\end{equation*}
where the last term comes out deriving at least one of the coefficients of the vector fields. With a slight abuse of notations we rewrite $\phi(x)\eta_i(f^{-n}(x))\prod_{j=1}^d g_j(x),$ which is again a test function, as $\phi(x).$  
By \cite[Lemma 6.5]{gouezel-liverani}, given a $C^{t+q}(U(W_i))$ vector field $v,$ there exist $C^{t+q}$ vector fields $w^u$ and $ w^s$, in a neighborhood $U(W_i)$ of $W_i$,  such that 
$v=w^u+w^s$, for all $x\in f^n(W_i)$ the vector $w_x^s\in T_xf^n(W_i),$ $\Vert w^s\Vert_{C^{t+q}(U_i)}\leq C_n,$ $\Vert w^s\circ f^n\Vert_{C^{t+q-1}(W_i)}\leq C$ and $\Vert d_{f^n(x)}f^{-n}w^u(f^n(x))\Vert_{C^{t+q}(f^{-n}(U(W)))}\leq C\lambda^{-n}.$
We obtain
\begin{equation*}
	\begin{split}
		\int_{f^n(W_i)}\phi(x)\prod_{j=1}^d \de_{x_j}(h\circ f^{-n})(x)=\sum_{\sigma\in\{s,u\}^t}\int_{f^n(W_i)}\phi(x)\prod_{j=1}^d w_j^{\sigma_j}(h\circ f^{-n})(x).
	\end{split}
\end{equation*}
Since $w_j^u w_k^s=w_k^s w_j^u+[w_j^u,w_k^s],$ we can swap two vector fields up to a term which is again $C_n\Vert h\Vert_{p-1,q+1,0}.$ Thus, we need to estimate terms of the form
\begin{equation*}
	\begin{split}
		&\int_{f^n(W_i)}\!\phi(x)\!\prod_{j=1}^g w_j^{s}\!\prod_{j=g+1}^dw_j^u(h\circ f^{-n})(x)=\!(-1)^g\!\int_{f^n(W_i)}\!\prod_{j=1}^g w_j^{s}\phi(x)\!\prod_{j=g+1}^dw_j^u(h\circ f^{-n})(x)
	\end{split}
\end{equation*}
where we applied an integration by parts on the $g$ derivatives along $f^n(W_i)$.
Every vector field $w^s_j$ can be written again in terms of the coordinate vector fields $w^s_j=\sum_{z=1}^db_{j,z}\de_{x_z}.$ As above, if one of the vector fields acts on one of the coefficients we get a term bounded by $C_n\Vert h\Vert_{p-1,q+1,o}.$ Therefore, we are left with the following terms
\begin{equation*}
	\begin{split}
		&\int_{f^n(W_i)}\prod_{j=1}^g \de_{x_j}\phi(x)\prod_{j=1}^gb_{j,z_j}(x)\prod_{j=g+1}^dw_j^u(h\circ f^{-n})(x)=\\&=\int_{W_i}\prod_{j=1}^g(\de_{x_j}\phi)\circ f^n(x)\prod_{j=1}^gb_{j,z_j}\circ f^n(x)\prod_{j=g+1}^d\widetilde{w}_j^u(h)\lambda_n^s(x).
	\end{split}
\end{equation*}
By the above estimates on norms, $\Vert \prod_{j=1}^gb_{j,z_j}(x)\circ f^n\Vert_{C^{t+q}}\leq C.$
We distinguish two cases 
If $t=p$ and $g=0$ it holds
$$\left|\int_{W_i}\prod_{j=1}^p\widetilde{w}^u_j(h)\phi\right|\leq C\lambda^{-np}\Vert h\Vert_{p,q,o}\Vert\phi\Vert_{C_0^{p+q}}.$$
On the other hand, if $t<p$ or $g>0,$ the integral can be estimated by $C_n\norm{h}_{p-1,q,l}\norm{\phi}_{C^{p+q}_0}$. This is not enough, because we need an estimate with $\norm{h}_{p-1,q+1,l}.$ To improve the regularity, let $\overline{\phi}=\prod_{j=1}^g\de_{x_j}\phi\in C^{q+t-g}_0(\widetilde{W}).$
We need to smoothen this function through the following standard lemma.
\begin{lemma}
	\label{lemma_Lasota_Yorke}
	Let $\alpha$ be the bigger integer smaller than $q+t-g.$ For $\epsilon>0$, there exists $\phi_\epsilon\in C^{q+t-g+1}$ such that $\Vert\phi_\epsilon\Vert_{C^{q+t-g}}\leq C\Vert\overline{\phi}\Vert_{C^{q+t-g}}$, $\Vert\phi_\epsilon\Vert_{C^{q+t-g+1}}\leq C\epsilon^{-1}\Vert\overline{\phi}\Vert_{C^{q+t-g}}$ and 
	$\Vert\phi_\epsilon-\overline{\phi}\Vert_{C^\alpha}\leq C\epsilon^{q+t-g-\alpha}\Vert\overline{\phi}\Vert_{C^{q+t-g}}$.
\end{lemma}
The proof of the above lemma easily follows convolving the function $\overline{\phi}$ with a mollifier of order $\epsilon$ and then computing the norms. In our context we fix $\epsilon$ to be $\lambda^{-(q+t-g)n/(q+t-g-\alpha)}$, so that 
$$\Vert\phi_\epsilon-\overline{\phi}\Vert_{C^\alpha}\leq C\lambda^{-(q+t-g)n}\Vert\overline{\phi}\Vert_{C^{q+t-g}}.$$
This implies that
\begin{equation*}
	\begin{split}
		&\left|\int_{W_i}\prod_{j=1}^g(\de_{x_j}\phi)\circ f^n(x)\prod_{j=1}^gb_{j,z_j}\circ f^n(x)\prod_{j=g+1}^d\widetilde{w}_j^u(h)\lambda_n^s(x)\right|\leq\\
		\leq&\left|\int_{W_i}(\overline{\phi}-\phi_\epsilon)\circ f^n(x)\prod_{j=1}^gb_{j,z_j}\circ f^n(x)\prod_{j=g+1}^d\widetilde{w}_j^u(h)\lambda_n^s(x)\right|+\\+&\left|\int_{W_i}\phi_\epsilon\circ f^n(x)\prod_{j=1}^gb_{j,z_j}\circ f^n(x)\prod_{j=g+1}^d\widetilde{w}_j^u(h)\lambda_n^s(x)\right|.
	\end{split}
\end{equation*}
The first term of the sum is bounded by
$$C\lambda^{-(q+t-g)n}\lambda^{-(t-g)n}\Vert\phi\Vert_{C^{q+t}}\Vert h\Vert_{p,q,0}\leq C\lambda^{-qn}\norm{\phi}_{C^{q+t}}\norm{h}_{p,q,0},$$
where we used that $\lambda^{-(q+t-g)}\lambda^{-(t-g)}\leq \lambda^{-q}$ when $g\leq t.$ Since the function $\phi_\epsilon$ is strictly more regular than $\overline{\phi},$ because $t-k<p,$ the second term is bounded by $C_n\Vert h\Vert_{p-1,q+1,0}.$ 
This concludes the proof of \eqref{eq_Lasota-Yorke_3}.

For \eqref{Lasota-Yorke_2} and \eqref{Lasota-Yorke_4}, we proceed as above. Hence, we firstly prove \eqref{Lasota-Yorke_2} for $p=0$, then we show the analogous of \eqref{eq_dependent_n} and, finally, we conclude by induction proving that \eqref{Lasota-Yorke_2} for $p-1$ implies \eqref{Lasota-Yorke_4} and \eqref{Lasota-Yorke_2} for $p$.  
Let us show \eqref{Lasota-Yorke_2} for $p=0$.
We consider $h\in\Omega^l(M),$ $W\in\Sigma$ and a test form $\omega\in\Gamma^{l,q}_0(W).$
We need to compute the integral 
$$\int_{W}\langle\omega,f_*^n h\rangle.$$
We fix the local bases for vector fields and forms introduced in section \ref{Section_Norms_and_Banach_spaces}. Assuming that $W\subseteq\psi_z(U_z,)$,  we can write $\omega=\omega\circ\chi_z=\sum_{\overline{j}\in\mathcal{J}_l}\omega_{\overline{j}}dx_{\overline{J}}$  on $\psi_z(\mathcal{B}(0,3\rho)).$ 
It follows that
\begin{align*}
	\int_W\langle\omega,f_*^nh\rangle=&\int_W\sum_{\overline{j}\in\mathcal{J}_l}\omega_{\overline{j}}\langle dx_{\overline{j}},f_*^nh\rangle=\sum_{\overline{j}\in\mathcal{J}_l}\int_{f^{-n}(W)}\omega_{\overline{j}}\circ f^n\langle dx_{\overline{j}},f_*^nh\rangle\circ f^n\lambda_n^s=\\=&\sum_{\overline{j},\overline{k}\in\mathcal{J}_l}\sum_i\int_{W_i}h_{\overline{k}}\omega_{\overline{j}}\circ f^n\langle dx_{\overline{j}},f_*^ndx_{\overline{k}}\rangle\circ f^n\lambda_n^s\eta_i.
\end{align*}
where $h=h\circ\chi_i=\sum_{\overline{k}\in\mathcal{J}_l}h_{\overline{k}}dx_{\overline{k}},$ on $\psi_i(\mathcal{B}(0,3\rho))\supseteq W_i.$
The $C^q$-norm of the functions $\langle dx_{\overline{j}},f_*^ndx_{\overline{k}}\rangle\circ f^n\lambda_n^s$ is bounded by $C\lambda^{-|d_s-l|n}$ (see Lemma \ref{lemma_determinant2}).
Thus,
\begin{equation*}
	\begin{split}
		&\left|\int_W\langle\omega,f_*^nh\rangle\right|\leq C\lambda^{-|d_s-l|n}|f^{-n}(W)|\Vert\omega\Vert_{\Gamma_0^{l,q}}\Vert h\Vert_{0,q,l}
	\end{split}
\end{equation*}
By classical results (see for instance \cite[Lemma C.3]{giulietti-liverani-pollicott}), the volume growth of $f^{-n}(W)$ fulfills $|f^{-n}(W)|\sim e^{h_\top n}|W|,$ where $h_\top $ is the topological entropy of $f$. By density, we get \eqref{Lasota-Yorke_2} for $p=0$. 

Next, let us prove that, for $h\in\Omega^l(M),$
\begin{equation}
	\label{eq_dependent_n2}
	\norm{f_*^nh}_{p,q,l}\leq
	C\lambda^{-n\min\{p,q\}-n|d_s-l|}e^{nh_\top }\norm{h}_{p,q,l}+ C_n\norm{h}_{p-1,q+1,l}.
\end{equation}
Let $h\in\Omega^l(M)$, $W\in\Sigma$, $\omega\in\Gamma_0^{p+q}(W)$ and let  $v_1,\dots,v_t\in\mathcal{V}^{t+q}(U(W))$ be $t$ vector fields such that  $\norm{v_i}_{C^{p+q}(U(W))}\leq1.$ As above, we can write $\omega=\sum_{\overline{j}\in\mathcal{J}_l}\omega_{\overline{j}}dx_{\overline{j}}$ and we compute 
\begin{align*}
	&\int_W\!\langle\omega,L_{v_1}\!\dots\! L_{v_p}f_*^nh\rangle\!=\!\sum_i\!\int_{f^n(W_i)}\sum_{\overline{j},\overline{k}\in\mathcal{J}_l}\omega_{\overline{j}}\langle dx_{\overline{j}},L_{v_1}\!\dots\! L_{v_p}[(h_{\overline{k}}\circ f^{-n})f_*^n(dx_{\overline{k}})]\rangle\eta_i\circ f^{-n}\\
	&=\sum_i\int_{f^n(W_i)}\sum_{\overline{j},\overline{k}\in\mathcal{J}_l}\sum_{A\subseteq\{1,\dots,p\}}\omega_{\overline{j}}\prod_{a\in A}L_{v_a}(h_{\overline{k}}\circ f^{-n})\langle dx_{\overline{j}},\prod_{a\in A^c}L_{v_a}(f_*^n(dx_{\overline{k}}))\rangle\eta_i\circ f^{-n}
\end{align*}
where  $h=\sum_{\overline{k}\in\mathcal{J}_l}h_{\overline{k}}dx_{\overline{k}}$ on $W_i$, while the product of derivatives is ordered. If $A=\{1,\dots,p\}$ the terms
$$\sum_i\int_{f^n(W_i)}\sum_{\overline{j},\overline{k}\in\mathcal{J}_l}\omega_{\overline{j}}\prod_{a=1}^pL_{v_a}(h_{\overline{k}}\circ f^{-n})\langle dx_{\overline{j}},f_*^n(dx_{\overline{k}})\rangle\eta_i\circ f^{-n}$$
can be treated putting together the proofs of \eqref{Lasota-Yorke_2}, for $p=0$, and \eqref{Lasota-Yorke_3}. Hence, it is bounded by $C\lambda^{-n\min\{p,q\}}\lambda^{-n|d_s-l|}e^{nh_\top }\Vert h\Vert_{p,q,l}+C\Vert h\Vert_{p-1,q+1,l}.$ Every other term has the form 
$$\sum_i\int_{f^{n}(W_i)}\sum_{\overline{j},\overline{k}\in\mathcal{J}_l}\omega_{\overline{j}}\prod_{a\in A}L_{v_a}(h_{\overline{k}}\circ f^{-n})\langle dx_{\overline{j}},\prod_{a\in A^c\neq\emptyset}L_{v_a}(f_*^n(dx_{\overline{k}}))\rangle.$$
Since  there are $t<p$ derivatives acting on $(h_{\overline{k}}\circ f^{-n}),$ \eqref{Lasota-Yorke_2}, for $p=0,$ and \eqref{Lasota-Yorke_3} imply that these terms are bounded by $C_n\Vert h\Vert_{p-1,q+1,l}$ and this concludes the proof of \eqref{eq_dependent_n2}. Finally, by the same inductive procedure used to prove \eqref{Lasota-Yorke_1} and \eqref{Lasota-Yorke_3}, one can prove \eqref{Lasota-Yorke_2} for $p>0$ and \eqref{Lasota-Yorke_4}. Moreover, by density, these inequalities extend to $\mathcal{B}^{p,q,l}.$ Notice that in this case we cannot expect  that coefficient in front of the weak norm in \eqref{Lasota-Yorke_4} is uniformly bounded in $n$. On the contrary, we have just proved that it cannot grow more that $\lambda^{-n|d_s-l|}e^{h_\top }.$

We are left with the proof of \eqref{Lasota-Yorke_5} and \eqref{Lasota-Yorke_6}. We have already noticed in \eqref{eq_transfer_op} the duality, induced by the volume form $\omega_0,$ between $d_s$-forms and functions. Accordingly, \eqref{Lasota-Yorke_5} and \eqref{Lasota-Yorke_6} hold true in $\Omega^{d_s}(M)$ for $f_*$ if and only if the same are satisfied by the transfer operator $\mathcal{L}$ acting on functions  $C^{\infty}(M)=\Omega^0(M).$ Lasota-Yorke inequalities for $\mathcal{L}$ were proved in \cite[Lemma 2.2]{gouezel-liverani}. By density, we conclude that \eqref{Lasota-Yorke_5} and \eqref{Lasota-Yorke_6} extend to $\mathcal{B}^{p,q,\dim(M)}.$

\qed 

\begin{remark}
	\label{rem_Lasota-Yorke}
	The inequalities \eqref{Lasota-Yorke_4} and \eqref{Lasota-Yorke_2}, for $l=d_s$, also implies the following inequality that we are going to use later (see the proof of Lemma \ref{lemma_inv_measure}). In fact, for $p>0$, $q>0$, $l=d_s$ and $\omega\in\mathcal{B}^{p,q,d_s},$ one can easily prove by induction that
	\begin{equation}
		\label{eq_Lasota-Yorke_rem}
		\norm{f_*^{pn}\omega}_{p,q,d_s}\leq Ce^{pnh_\top }\lambda^{-n}\norm{h}_{p,q,d_s}+Ce^{pnh_\top }\norm{h}_{0,p+q,d_s},
	\end{equation}
	where the $p$ in $f_*^{pn}$ represents the same parameter $p$ of the norm. Indeed, \eqref{Lasota-Yorke_4} gives \eqref{eq_Lasota-Yorke_rem} for $p=1$. Assume \eqref{eq_Lasota-Yorke_rem} true up to $p-1$. Then, by using \eqref{Lasota-Yorke_2}, \eqref{Lasota-Yorke_4} and the property $\norm{\cdot}_{p-1,q+1,d_s}\leq\norm{\cdot}_{p,q,d_s},$ we obtain
	\begin{equation*}
		\begin{split}
			\norm{f_*^{pn}h}_{p,q,d_s}\leq& Ce^{nh_\top }\lambda^{-n}\norm{f_*^{(p-1)n}h}_{p,q,d_s}+Ce^{nh_\top }\norm{f_*^{(p-1)n}h}_{p-1,q+1,d_s}\leq\\
			\leq&Ce^{pnh_\top }\lambda^{-n}\norm{h}_{p,q,d_s}+Ce^{2nh_\top }\lambda^{-n}\norm{f_*^{(p-2)n}h}_{p-1,q+1,d_s}+\\&+Ce^{pnh_\top }\norm{h}_{0,p+q,d_s}\leq\\
			\leq&Ce^{pnh_\top }\lambda^{-n}\norm{h}_{p,q,d_s}+Ce^{pnh_\top }\norm{h}_{0,p+q,d_s},
		\end{split}
	\end{equation*}
	which proves \eqref{eq_Lasota-Yorke_rem}.
\end{remark}

\begin{corollary}
	\label{cor_spectral_radius}
	The spectral radius and the essential spectral radius of the pushforward operator $f_*\colon\mathcal{B}^{p,q,l}\rightarrow\mathcal{B}^{p,q,l}$ fulfill $$\rho(f_*|_{\mathcal{B}^{p,q,l}})\leq\begin{cases}1&\text{ if }l=0\text{ or }l=d\\
		\lambda^{-|d_s-l|}e^{h_\top }&\text{ if }0<l<d
	\end{cases}$$$$\ \ \rho_{ess}(f_*|_{\mathcal{B}^{p,q,l}})\leq\begin{cases}\lambda^{-\min\{p,q\}}&\text{ if }l=0\text{ or }l=d\\
		\lambda^{-\min\{p,q\}}\lambda^{-|d_s-l|}e^{h_\top }&\text{ if }0<l<d
	\end{cases}$$
\end{corollary}
\proof
The estimates on spectral radii follow by \eqref{Lasota-Yorke_1}, \eqref{Lasota-Yorke_2} and \eqref{Lasota-Yorke_5} applying the formula $\rho(f_*)=\lim_{n\rightarrow+\infty}\norm{f_*^n}^{\frac{1}{n}}.$
The estimates on the essential spectral radii are consequences of Hennion's theorem \ref{thm_hennion} whose hypotheses are satisfied by Lemma \ref{inclusione_compatta} and Theorem \ref{thm_lasota-yorke}.
\qed

Once we have established this spectral picture a natural question may arise: how does the spectrum of $f_*$ acting on $\mathcal{B}^{p,q,l}$, denoted by $\sigma(f_*|_{\mathcal{B}^{p,q,l}}),$ depends on $p$ and $q$?
The following lemma answers this question, at least for the spectrum we are interested in.
\begin{lemma}
	\label{lemma_same_spectrum}
	Let $\mathcal{B}^{p,q,l}$ and $\mathcal{B}^{p',q',l}$ be two anisotropic Banach spaces of currents for some parameters $p,q,p',q'\in \N.$ Assume that $\rho_{ess}(f_*|_{\mathcal{B}^{p',q',l}})\leq\rho_{ess}(f_*|_{\mathcal{B}^{p,q,l}}).$ Then 
	$$\sigma(f_*|_{\mathcal{B}^{p,q,l}})\cap\{z\in\C:\ |z|>\rho_{ess}(f_*|_{\mathcal{B}^{p,q,l}})\}=\sigma(f_*|_{\mathcal{B}^{p',q',l}})\cap\{z\in\C:\ |z|>\rho_{ess}(f_*|_{\mathcal{B}^{p,q,l}})\}.$$
	Moreover, the corresponding generalized eigenspaces coincide and they are included in $\mathcal{B}^{p,q,l}\cap\mathcal{B}^{p',q',l}.$ 
\end{lemma}
\proof
This is a consequence of \cite[Lemma A.1]{baladi-tsujii2}. We recall the statement for the sake of completeness. 
\begin{lemma}{\cite[Lemma A.1]{baladi-tsujii2}} 
	\label{lemma_baladi_tsujii}
	Let $\mathcal{B}$ be a separated topological linear space and let $\mathcal{B}_1,$ $\mathcal{B}_2$ be Banach spaces such that $\mathcal{B}_i\hookrightarrow\mathcal{B}$ is continuous, for $i=1,2.$ Moreover, assume that there exists a subspace $\mathcal{B}_0\subset\mathcal{B}_1\cap\mathcal{B}_2,$ which is dense in both spaces $\mathcal{B}_1$ and $\mathcal{B}_2$. Let $L\colon\mathcal{B}\rightarrow\mathcal{B}$ be a continuous linear map that preserves $\mathcal{B}_0$, $\mathcal{B}_1$ and $\mathcal{B}_2$. Suppose that $L\colon\mathcal{B}_1\rightarrow\mathcal{B}_1$ and $L\colon\mathcal{B}_2\rightarrow\mathcal{B}_2$ are bounded linear operators with essential spectral radius strictly smaller than some $\rho>0$. Then the spectrum of $L|_{\mathcal{B}_1}$ and the spectrum of $L|_{\mathcal{B}_2}$ in $\{z\in\C:\ |z|>\rho\}$ coincide. Moreover, the corresponding generalized eigenspaces coincide and belong to $\mathcal{B}_1\cap\mathcal{B}_2$.
\end{lemma}
Lemma \ref{lemma_same_spectrum} follows from Lemma \ref{lemma_baladi_tsujii} by taking 
 $\mathcal{B}_0=\Omega^l(M),$ $\mathcal{B}_1=\mathcal{B}^{p,q,l},$ $\mathcal{B}_2=\mathcal{B}^{p',q',l}$ and $\mathcal{B}=\mathcal{B}^{\min\{p,p'\},\max\{q,q'\},l}.$ Indeed, $f_*|_{\mathcal{B}}$ is a bounded operator by Theorem \ref{thm_lasota-yorke}, $\mathcal{B}_0$ is dense in $\mathcal{B}_1$ and $\mathcal{B}_2$ by definition, $\mathcal{B}_1,\mathcal{B}_2\subseteq \mathcal{B},$ since $\mathcal{B}^{p,q,l}\subset\mathcal{B}^{p-1,q+1,l},$ and the essential spectral radius of $f_*|_{\mathcal{B}_1}$ and $f_*|_{\mathcal{B}_2}$  is strictly larger than $\rho=\rho_{\ess}(f_*|_{\mathcal{B}_1}).$
 
\qed
\section{The spectrum of the pushforward operator}
\label{sec_spectrum}
We now want to investigate the spectrum  of the pushforward operator acting on anisotropic Banach spaces $\mathcal{B}^{p,q,l}$, $0\leq l\leq\dim(M)$.  In particular, we are going to use information about the spectrum of $f_*\colon\mathcal{B}^{p,q,d_s}\rightarrow \mathcal{B}^{p,q,d_s}$, where $d_s$ is the dimension of the stable bundle. 

From now on, we always assume  $p$ and $q$ large enough, so that there exists $\nu\in(0,1)$ such that 
\begin{equation}
	\label{eq_nu}
	\max_{i=0,\dots,d_s}\max\{\lambda^{-\min\{p\pm i,q\mp i\}-|d_s-l|}e^{h_\top }, \lambda^{-\min\{p\pm i,q\mp i\}}\}<\nu<1,
\end{equation}
for any $i=0,\dots,\dim(M)$.
Corollary \ref{cor_spectral_radius} and \eqref{eq_nu} ensure that the essential spectral radius $f_*$ acting on the Banach spaces we are interested in is bounded by $\nu.$

\subsection{Peripheral spectrum and measure of maximal entropy}
\label{sec_peripheral_spectrum}
By peripheral spectrum we mean the set of eigenvalues of $f_*|_{\mathcal{B}^{p,q,d_s}}$ of maximal modulus.  Corollary \ref{cor_spectral_radius} only tells us that the spectral radius $\rho(f_*|_{\mathcal{B}^{p,q,d_s}})$ is bounded by $e^{h_\top }.$ On the other hand, the following lemma proves that this upper bound is actually attained.
\begin{lemma}
	\label{lemma_lower_bound}
	The spectral radius of $f_*|_{\mathcal{B}^{p,q,d_s}}$ is exactly $e^{h_\top }.$
\end{lemma}
Before giving the proof of Lemma \ref{lemma_lower_bound}, we introduce the following $d_s$-differential form $\omega_\Sigma\in\Omega^{d_s}(M),$ that we are going to use along this section. We need a $d_s$-form which gives positive volume to every admissible stable leaf $W\in\Sigma,$ in the sense that 
$\int_W\omega_\Sigma>0$, for any $W\in\Sigma.$ The first idea would be to consider the volume $\omega_W$ induced by $\omega_0$ on every admissible leaf $W\in\Sigma, $ but this is only a $d_s$-form on $W$ and not on $M$.  On the other hand, on every chart $(U_i,\psi_i),$ one can easily define a $d_s$-differential form $u_i$ which gives positive volume to every leaf with tangent space in the Euclidean stable cone bundle $\zeta^s.$ Consequently, $\psi_i^* u_i\in\Omega^{d_s}(\psi_i(U_i))$ gives positive volume to every $W\in\psi_i(U_i).$ Finally, by using the partition of unity $\{\chi_i\}_{i=1}^m$, we define $\omega_\Sigma=\sum_{i=1}^{m}\chi_i\psi_i^*(u_i)\in\Omega^{d_s}(M).$ \\
\proofof{Lemma \ref{lemma_lower_bound}} We have already used, in the proof of Theorem \ref{thm_lasota-yorke}, that  $$e^{h_\top }=\limsup_{n\rightarrow +\infty}\left(\sup_{W\in\Sigma}|f^{-n}(W)|\right)^{\frac{1}{n}},$$
where $|f^{-n}(W)|$ is the volume of $f^{-n}(W)$ w.r.t.\ the measure induced by the Riemannian volume $\omega_0$ on $f^{-n}(W)$ (see \cite[Lemma C.3]{giulietti-liverani-pollicott}). Let  $\omega_{\Sigma}\in\Omega^{d_s}(M)$ the differential form defined above which gives positive volume to every admissible stable leaf of $\Sigma$. Then, given any $W\in\Sigma,$ let $\{W_i\}_{i=1}^l$ be the covering of $f^{-n}(W),$ constructed in Lemma \ref{lemma_foglie}. By compactness of every $W_j,$ there exists $C_j$ such that 
$$\int_{W_j}\omega_{W_j}\leq C_j\int_{W_j}\omega_\Sigma,$$
where $\omega_{W_j}$ is the volume form induced by $\omega_0$ on $W_j.$ As a consequence, 
\begin{equation*}
	\begin{split}
		|f^{-n}(W)|&\leq \sum_{j=1}^l |W_j|\leq \sum_{j=1}^l C_j\int_{W_j}\omega_\Sigma \leq \sum_{j=1}^l C_j\int_{f^n(W_j)}f_*^n\omega_\Sigma\\&\leq C\int_{W}f_*^n\omega_\Sigma\leq C \norm{f_*^n\omega_\Sigma}_{p,q,d_s}\leq C\norm{f_*^n}_{\mathcal{B}^{p,q,d_s}\rightarrow\mathcal{B}^{p,q,d_s}},
	\end{split}
\end{equation*}
hence
\begin{equation*}
	\begin{split}
		e^{h_\top }=\limsup_{n\rightarrow +\infty}\left(\sup_{W\in\Sigma}|f^{-n}(W)|\right)^{\frac{1}{n}}\leq \limsup_{n\rightarrow+\infty}\norm{f_*^n}^{\frac{1}{n}}_{\mathcal{B}^{p,q,d_s}\rightarrow\mathcal{B}^{p,q,d_s}}=\rho(f_*|_{\mathcal{B}^{p,q,d_s}}).
	\end{split}
\end{equation*}
Corollary \ref{cor_spectral_radius} gives $\rho(f_*|_{\mathcal{B}^{p,q,d_s}})\leq e^{h_\top },$ hence we conclude that $\rho(f_*|_{\mathcal{B}^{p,q,d_s}})= e^{h_\top }.$

\qed

As a consequence of Lemma \ref{lemma_lower_bound} and the quasi-compactness of $f_*$, we can write 
\begin{equation}
	\label{eq_sviluppo}
	f_*=\sum_{i=0}^{N}(z_ie^{h_\top })\Pi_i+R,
\end{equation}
where every $z_i$ is a complex number of modulus $1$, the operator $\Pi_i$ is the finite rank projection on the eigenspace corresponding to the eigenvalue $z_i e^{h_\top }$ and $R$ is a bounded linear operator whose spectral radius is strictly smaller than $e^{h_\top }$. Moreover, $\Pi_i\circ\Pi_j=\delta_{i,j}\Pi_i$ and $\Pi_i\circ R=R\circ\Pi_i=0$. Notice that, as a consequence of \eqref{Lasota-Yorke_2}, the operator $e^{-nh_\top }f_*^n$ is bounded for all $n$, i.e.,  there cannot be Jordan blocks for eigenvalues of modulus $e^{h_\top }$. In fact, let us assume by contradiction that there is a Jordan block of dimension 2 (the extension to higher dimension is very similar) for the eigenvalue $\mu_i=z_ie^{h_\top}.$ Then, there two generalized eigenvectors $h_1, h_2\in\mathcal{B}^{p,q,d_s}$ such that $f_*h_1=\mu_i h_1$ and $f_*h_2=\mu_i h_2+h_1.$ By induction, one can easily show that $f_*^nh_2=\mu_i^nh_2+n\mu_i^{n-1}h_1.$ It follows that, for large $n\in\N,$ $$\norm{\mu_i^{-n}f_*^nh_2}_{p,q,d_s}\geq ne^{-h_\top}\norm{h_2}_{p,q,d_s}-\norm{h_1}_{p,q,d_s}\xrightarrow[n\rightarrow +\infty]{}+\infty.$$
On the other hand, \eqref{Lasota-Yorke_2} implies that $\norm{\mu_i^{-n}f_*^nh_2}_{p,q,d_s}\leq C\norm{h_2}_{p,q,d_s},$ giving a contradiction. 
\begin{remark}
	We underline that in the more general case of \cite{gouezel-liverani08}, where the authors consider equilibrium states for compact locally maximal hyperbolic sets of smooth maps, the absence of Jordan blocks requires a more involved proof and the assumption that the system is topologically mixing. Since we are assuming that $f$ is topologically transitive, which is equivalent to topological mixing for Anosov diffeomorphisms \cite[Corollary 18.3.5]{Katok}, we should have used this hypothesis at some point. In fact, the growth of stable leave $|f^{-n}(W)|\sim e^{h_\top n}|W|$ is proved in \cite[Lemma C.3]{giulietti-liverani-pollicott} assuming $f$ topologically transitive. As a consequence, we obtained a more precise Lasota-Yorke inequality in \eqref{Lasota-Yorke_2}, compared to that of \cite{gouezel-liverani08}, and an easier proof of the absence of Jordan blocks for the peripheral spectrum. 
\end{remark}

\begin{remark}
	Up to now, we have not considered orientability issues. We have only assumed that $M$ is an orientable manifold, but we have never made any assumption about the orientation of the stable/unstable foliation. In fact, up to considering a finite covering of $M$ we can always assume that these two foliations are oriented. Moreover, we can also  suppose that $f$ preserves the orientation of both foliations. Otherwise, it would be enough to consider $f^2$ in place of $f$. Accordingly, from now on, we assume to work with a diffeomorphism $f$ preserving the orientation of the oriented stable and unstable manifolds.
\end{remark}
The rest of this section is devoted to the proof of the following proposition.
\begin{proposition}
	\label{prop_measure_of_maximal_entropy}
	$e^{h_\top }$ is the unique simple maximal eigenvalue of $f_*$ acting on $\mathcal{B}^{p,q,d_s}.$ Let $\hat\omega\in\mathcal{B}^{p,q,d_s}$ be a corresponding eigenvector and let $\hat t\in(\mathcal{B}^{p,q,d_s})'$ be the dual  eigenvector such that $\hat t(\hat\omega)=1.$ 
	The continuous linear operator $\phi\mapsto\hat t(\phi\hat\omega),$ defined on $C^{p+q}$-functions, extends to a bounded linear operator on $C^0$-functions, i.e., it is a measure. 
	In particular, $\mu_{BM}(\cdot)=\hat t(\cdot\ \hat\omega)$ is a positive measure and  it is the unique measure of maximal entropy.
\end{proposition}
We point out that parts of the proof are inspired by \cite[Section 4-6]{gouezel-liverani08}, where the authors showed the existence of equilibrium states for locally compact hyperbolic set. However, the authors of \cite{gouezel-liverani08} defined anisotropic Banach spaces, starting from Grassmannian spaces, while we introduced them as the completion of differential forms spaces. Moreover, we limit to consider the measure of maximal entropy, i.e., the unique equilibrium state for the null potential. 

\proofof{Proposition \ref{prop_measure_of_maximal_entropy}}
Let us consider the $d_s$-differential form $\omega_\Sigma$ that gives positive volume to every admissible leaf $W\in\Sigma$ and that is defined before the proof of Lemma \ref{lemma_lower_bound}. We set $\hat\omega=\Pi_1\omega_\Sigma,$ where $\Pi_1$ is the eigenprojector related to the eigenvalue $e^{h_\top }$ of \eqref{eq_sviluppo}. Accordingly, $f_*\hat\omega= e^{h_\top }\hat\omega.$ Notice that, a priori, $\hat\omega$ could be null, because we do not know yet that $e^{h_\top }$ is an eigenvalue of $f_*$. On the other hand, we are going to prove that $\hat\omega$ is actually nonzero.

The next two lemmas recall \cite[Lemma 4.9, Lemma 4.10]{gouezel-liverani08} adapting them to our setting.
\begin{lemma}
	\label{lemma_abs_cont}
	Let $\omega\in\mathcal{B}^{p,q,d_s}$ be an eigenvector for the eigenvalue $z_ie^{h_\top }$ such that $|z_i|=1.$ Then $\omega$ induces a measure on every admissible leaf $W\in\Sigma.$ Moreover, every such $\omega$ is absolutely continuous with bounded density w.r.t.\ the measure defined by $\hat\omega$.
\end{lemma} 
\proofof{Lemma \ref{lemma_abs_cont}}
We firstly show that, given $W\in\Sigma$ and $\phi\in\Gamma_0^{q,d_s}(W)$, it holds
\begin{equation}
	\label{eq_estimate}
	\left|\int_W\langle\phi,\omega \rangle\omega_W\right|\leq C\norm{\phi}_{\Gamma_0^{0,d_s}}.
\end{equation}
Since $\Omega^{d_s}(M)$ is dense in $\mathcal{B}^{p,q,d_s}$, $\Pi_{z_i}$ is continuous and $\Pi_{z_i}\Omega^{d_s}(M)$ is closed, we get $\Pi_{z_i}\Omega^{d_s}(M)=\Pi_{z_i}\mathcal{B}^{p,q,d_s}.$ Thus, there exists  a smooth form $\widetilde{\omega}\in\Omega^{d_s}(M)$ such that $\Pi_{z_i}\widetilde{\omega}=\omega$ and, by \eqref{eq_sviluppo},
\begin{equation}
	\label{eq_eigenproj}
	\begin{split}
		\int_W\langle\phi,\omega\rangle\omega_W=\lim_{n\rightarrow+\infty}\frac{1}{n}\sum_{k=0}^{n-1}(z_ie^{h_\top })^{-k}\int_W\langle\phi,f_*^k\widetilde{\omega}\rangle\omega_W.
	\end{split}
\end{equation}
Therefore,
\begin{equation*}
	\begin{split}
		\left|\int_W\langle\phi,\omega\rangle\omega_W\right|&\leq\lim_{n\rightarrow+\infty}\frac{1}{n}\sum_{k=0}^{n-1}e^{-kh_\top }\left|\int_W\langle\phi,f_*^k\widetilde{\omega}\rangle\omega_W\right|\\&\leq\lim_{n\rightarrow+\infty}\frac{1}{n}\sum_{k=0}^{n-1}e^{-kh_\top }\norm{\phi}_{\Gamma_0^{0,d_s}}\norm{f_*^k\widetilde{\omega}}_{0,0,d_s}\leq C \norm{\phi}_{\Gamma_0^{0,d_s}}\norm{\widetilde{\omega}}_{0,0,d_s},
	\end{split}
\end{equation*}
where last inequalities follows by the smoothness of $\widetilde{\omega}$ and by \eqref{Lasota-Yorke_2}.

We point out that, since $\omega\in\mathcal{B}^{p,q,d_s},$ then above integral \eqref{eq_eigenproj} can always be estimated by the $C^q$-norm of $\phi$. On the contrary, \eqref{eq_estimate} implies that, for eigenvectors corresponding to maximal eigenvalues, the bound is given by the $C^0$ norm of $\phi$. 

Next, let us consider a function $\psi\in C^0_0(W),$ then $\psi\omega_W\in\Gamma_0^{0,d_s}$ and we define $$\int_W\psi\mathcal{M}_W(\omega)=\int_W\langle\psi\omega_W,\omega\rangle\omega_W.$$
As a consequence of \eqref{eq_estimate},
$$\left|\int_W\psi\mathcal{M}_W(\omega)\right|\leq C\norm{\psi}_{C^0(W)},$$
hence $\omega$ defines a measure on any admissible leaf $W\in\Sigma$. When $\omega=\hat\omega=\Pi_1\omega_\Sigma,$ the equality \eqref{eq_eigenproj} implies that $\mathcal{M}_W(\hat\omega)$ is a nonnegative measure. In addition, for every $\psi\in C^q_0(W),$
\begin{equation*}
	\begin{split}
		&\left|\int_W\psi\mathcal{M}_W(\omega)\right|=\left|\int_W\langle\psi\omega_W,\omega\rangle\omega_W\right|=\left|\int_W\langle\psi\omega_W,\lim_{n\rightarrow+\infty}\frac{1}{n}\sum_{k=0}^{n-1}(z_ie^{h_\top })^{-k}f_*^k\omega\rangle\omega_W\right|\leq\\\leq&\int_W\left|\langle\psi\omega_W,\lim_{n\rightarrow+\infty}\frac{1}{n}\sum_{k=0}^{n-1}e^{-kh_\top }f_*^k\omega\rangle\right|\omega_W=\int_W|\langle\psi\omega_W,\Pi_1\omega\rangle|\omega_W\leq \\\leq &C\int_W|\langle\psi\omega_W,\hat\omega\rangle|\omega_W\leq C\int_W|\psi|\mathcal{M}_W(\hat\omega)
	\end{split}
\end{equation*}
and, since above inequality extends to continuous functions by density, the measure $\mathcal{M}_W(\omega)$ is absolutely continuous w.r.t.\ $\mathcal{M}_W(\hat\omega)$ with bounded density. 

\qed

Let us continue the proof of Proposition \ref{prop_measure_of_maximal_entropy}.
Let $\omega\in\Omega^{d_s}$ and let $W,W'\in\Sigma$ be admissible leaves, whose intersection is again a $d_s$-dimensional submanifold. In addition, assume that the orientations of $W$ and $W'$ agree on $W\cap W'$.  Then, for every  function $\psi\in C^{q}_0(W), $  the operator $\mathcal{M}_W(\omega)$, such that 
$$\int_W\psi\mathcal{M}_W(\omega)=\int_W\langle\psi\omega_W,\omega\rangle\omega_W,$$ 
is a bounded operator on $C^q_0(W).$ Moreover, for $\psi\in C^q_0(W\cap W')$
\begin{equation*}
	\begin{split}
		&\int_W\psi\mathcal{M}_W(\omega)=\int_{W}\langle\psi\omega_W,\omega\rangle\omega_W=\int_{W}\langle\psi i_W^*\omega_0,\omega\rangle i_W^*\omega_0=\\=&\int_{W'}\langle\psi i_{W'}^*\omega_0,\omega\rangle i_{W'}^*\omega_0=\int_{W'}\langle\psi\omega_{W'},\omega\rangle\omega_{W'}=\int_{W'}\psi\mathcal{M}_{W'}(\omega),
	\end{split}
\end{equation*}
where $i_W^*$, resp.\ $i_{W'}^*,$ is the pull-back of the embedding of $W,$ resp $W'$, in $M$. Accordingly,  every $\omega\in\Omega^{d_s}$ defines an element of the dual of $C^q_0(\mathcal{S})$, where $\mathcal{S}$ is every manifold obtained by gluing elements of $\Sigma.$  By density, if $\omega\in\mathcal{B}^{p,q,l}$ is an eigenvector such that $\Pi_{z_i}\omega=\omega,$ then $\omega$ induces a measure, denoted by $\mathcal{M}(\omega)$, on every oriented stable manifold of $M$.

\begin{lemma}
	\label{lemma_injectivity}
	The function $\omega\mapsto\mathcal{M}(\omega),$ defined on the eigenspace $\Pi_{z_i}\mathcal{B}^{p,q,l}$ corresponding to the eigenvalue $z_ie^{h_\top },$ is injective. Additionally, $\hat\omega$ is nonzero.
\end{lemma}
\proofof{Lemma \ref{lemma_injectivity}}
Assume that $\mathcal{M}(\omega)=0$ for some $\omega\in\Pi_{z_i}\mathcal{B}^{p,q,d_s}.$ We firstly show that $\norm{\omega}_{0,q,d_s}=0.$ As a consequence of the compact inclusion of Lemma \ref{inclusione_compatta},  there exists a constant $C>0$ such that, for any $\epsilon>0$, there are a finite number of admissible leaves $W_1,\dots, W_k\in\Sigma$, such that for every $W\in\Sigma$ and for any $\phi\in\Gamma_0^{q,d_s},$ there is at least one $W_j$ such that
\begin{equation}
	\begin{split}
		\label{eq_epsilon_closed}
		\left|\int_W\langle\phi,\omega\rangle\omega_W-\int_{W_j}\langle\bar\phi_1,\omega\rangle\omega_{W_j}\right|\leq C\epsilon\norm{\phi}_{\Gamma_0^{q,d_s}},
	\end{split}
\end{equation}
where $\bar\phi_1=\Psi_1^*\phi$ and $(\Psi_t)_{t\in[0,1]}$ is the holonomy map from $W$ to $W_j$. Since  $W_j$ is a compact oriented $d_s$-dimensional manifold endowed with a volume form $\omega_{W_j}$, there exists $\bar\psi_1\in C_0^q(W_j)$ such that $\bar\phi_1=\bar\psi_1\omega_{W_j}.$ Consequently, if $W_j$ is contained in a stable manifold, then 
$$\int_{W_j}\langle\bar\phi_1,\omega\rangle\omega_{W_j}=\int_{W_j}\bar\psi_1\mathcal{M}(\omega)=0$$
and above inequality becomes
\begin{equation*}
	\begin{split}
		\left|\int_W\langle\phi,\omega\rangle\omega_{W}\right|\leq C\epsilon\norm{\phi}_{\Gamma_0^{q,d_s}}.
	\end{split}
\end{equation*}
More generally, there exist $n_0\in\N$ and a sequence $\{\epsilon_n\}_{n>n_0},$ going exponentially fast to zero, such that, for each  $n>n_0$ and for any full admissible leaf $\widetilde{W}\in\widetilde\Sigma,$ 
$$f^{-n}(W)\subseteq \cup_{i=1}^l W_i^{(n)}\subseteq f^{-n}(\widetilde{W}),$$
as stated in Lemma \ref{lemma_foglie}. In addition, every $W_i^{(n)}$ is $\epsilon_n$-closed to some leaf contained in a stable manifold in the sense of \eqref{eq_epsilon_closed}.
Next, given $W\in\Sigma$ and $\phi\in\Gamma_0^{q,d_s}(W),$ there exists $\psi\in C^q_0(W)$ such that $\phi=\psi\omega_W.$  We compute
\begin{equation*}
	\begin{split}
		&\left|\int_W\langle\phi,\omega\rangle\omega_W\right|=\left|e^{-nh_\top }\int_W\langle\phi,f_*^n\omega\rangle\omega_W\right|=\\=&\left|e^{-nh_\top }\int_W(-1)^{l(\dim(M)-l)}\langle\star (f^n)^*\star\phi,\omega\rangle\circ f^{-n}\det(df^{-n})\omega_W\right|\leq\\
		\leq&e^{-nh_\top }\sum_i\left|\int_{W_i^{(n)}}(-1)^{l(\dim(M)-l)}\langle\star (f^n)^*\star\phi,\omega\rangle\det(df^{-n})\circ f^n\lambda_n^s\omega_{W_i}\right|,
	\end{split}
\end{equation*}
where $\star\colon \Omega^l(M)\rightarrow \Omega^{\dim(M)-l}(M)$ is the Hodge $\star$ operator (see for instance \cite{de_Cataldo}) and $\lambda_n^s$ denotes again the Jacobian of the change of variables. 
Notice that a trivial computation gives $$\langle\star (f^n)^*\star\phi,\omega\rangle=\psi\circ f^n\langle\star (f^n)^*\star\omega_W,\omega\rangle,$$ while $$(-1)^{l(\dim(M)-l)}\star (f^n)^*\star\omega_W=(\lambda_n^s)^{-1}\det(df^n)\omega_W.$$ Indeed, by definition 
$\omega_W\wedge\star\omega_W=\omega_0$, resp.\ $\omega_{W_i^{(n)}}\wedge\star\omega_{W_i^{(n)}}=\omega_0,$ on $W$, resp.\ on $W_i^{(n)},$ hence 
\begin{equation*}
	\begin{split}
		\lambda_n^s\omega_{W_i^{(n)}}\wedge(f^n)^*\star\omega_W=(f^n)^*\omega_W\wedge (f^n)^*\star\omega_W=\det(df^n)\omega_{W_i^{(n)}}\wedge\star\omega_{W_i^{(n)}},
	\end{split}
\end{equation*}
that is $\langle\omega_{W_i^{(n)}},(-1)^{l(\dim(M)-l)}\star(f^n)^*\star\omega_W\rangle=\langle\omega_{W_i^{(n)}},(\lambda_n^s)^{-1}\det(df^n)\omega_{W_i^{(n)}}\rangle$ and this proves the second equality. 
Continuing the calculation
\begin{equation*}
	\begin{split}
	\left|\int_W\langle\phi,\omega\rangle\omega_W\right|&=e^{-nh_\top }\sum_i\left|\int_{W_i^{(n)}}\psi\circ f^n\langle\omega_{W_i^{(n)}},\omega\rangle\omega_{W_i^{(n)}}\right|\\&\leq C\epsilon_ne^{-nh_\top }\#\{W_i^{(n)}\}\norm{\psi\circ f^n}_{C^q}.
	\end{split}
\end{equation*}
Since $|e^{-nh_\top }\#\{W_i^{(n)}\}|\leq C$ and $\epsilon_n$ decays exponentially fast to zero, we obtain that $\norm{\omega}_{0,q,d_s}=0.$
Let us proceed by induction on $p$ in order to prove that $\norm{\omega}_{p,q,d_s}=0.$ Assume that the result is true up to $p-1.$ Then, using the Lasota-Yorke inequality \eqref{Lasota-Yorke_4},
\begin{equation*}
	\norm{\omega}_{p,q,d_s}=e^{-nh_\top }\norm{f_*^n\omega}_{p,q,d_s}\leq C\lambda^{-n\min\{p,q\}}\norm{\omega}_{p,q,d_s}\xrightarrow[n\rightarrow+\infty]{}0.
\end{equation*}
We conclude that $\omega=0,$ that is $\omega\mapsto\mathcal{M}(\omega)$ is injective. Finally, assume by contradiction that $\hat\omega=0,$ then, by injectivity $\mathcal{M}(\hat\omega)=0.$ Thus, for any other eigenvector $\omega$ corresponding to an eigenvalue of modulus $e^{h_\top }$, it must hold that $\mathcal{M}(\omega)=0,$ because $\mathcal{M}(\omega)$ is absolutely continuous with respect to $\mathcal{M}(\omega)$ by Lemma \ref{lemma_abs_cont}. We conclude that any such $\omega=0,$ hence the spectral radius of $f_*$ acting on $\mathcal{B}^{p,q,d_s}$ is strictly smaller than $e^{h_\top }$ and this contradicts Lemma \ref{lemma_lower_bound}.

\qed
The following step consists in proving that, assuming $f$ topologically transitive, $e^{h_\top }$ is a simple eigenvalue and it is the unique maximal eigenvalue of $f_*|_{\mathcal{B}^{p,q,d_s}}.$  Before proceeding with the proof, we need to recall some basic notions regarding continuous leafwise measure and their properties. We avoid to rewrite the proofs of results about this concept, but refer the reader to the survey \cite[Section 9]{gouezel-liverani08}, where they are proved in great generality.
\begin{definition}
	Let $X$ be a locally compact space. We assume that there exists an atlas $\{U,\phi_U\}$ such that $U\subseteq X$ is open and it is homeomorphic to $\mathcal{B}_d(0,1)\times K_U,$ for some locally compact space $K_U$, under the homeomorphism $\phi_U$. In addition, we assume that the changes of coordinates fulfill $\phi_U\circ \phi_V^{-1}(x,y)=(a(x,y),b(y)),$ i.e., they map leaves to leaves.
	A continuous leafwise measure $m$ is a family of Radon measures, each one defined on a leaf, such that, for every continuous function $\psi$ supported on the chart $(U,\phi_U),$ the integral 
	$$I_\psi(y)=\int_{\phi_U^{-1}(\mathcal{B}_d(0,1)\times\{y\})}\psi(x)\ dm(x)$$
	is a continuous function of $y$.
\end{definition}
Suppose that there exists a family of metrics on $X$, each one defined on a leaf, such that they also vary continuously with the leaf. Let $T\colon X\rightarrow X,$ be a continuous, leaves preserving homeomorphism, which uniformly expands distances on every leaf, i.e., there exists $\delta>0$ and $C>1$, such that $d_W(Tx,Ty)\geq C d_W(x,y)$, whenever $x,y$ belong to the leaf $W$ and $d_W(x,y)<\delta$.

We can now recall the result we need in order to prove the following Lemma \ref{lemma_unique_maximal}.
\begin{proposition}{\cite[Proposition 9.1, Proposition 9.4]{gouezel-liverani08}}\label{prop_leafwise_measure}\\ Let $m$ be a nonnegative continuous leafwise measure and let $m'$ be another complex continuous leafwise measure. Assume that there exists $C>0$ for which $|m'|\leq Cm$ on every leaf. Moreover, suppose that $T^*m=m$ and $T^*m'=\gamma m'$, with $|\gamma|=1.$ Finally, assume that $T$ is topologically mixing and that given any open set O of a leaf, there exists $x\in O$ with dense orbit. Then there is a $c\in\C$ such that $m'=cm,$ hence $\gamma=1$ or $m'=0.$
\end{proposition}
\begin{lemma}
	\label{lemma_unique_maximal}
	Under the assumption that $f$ is topologically transitive, $e^{h_\top }$ is the unique maximal eigenvalue of $f_*|_{\mathcal{B}^{p,q,d_s}}.$ In addition, $e^{h_\top }$ is simple.
\end{lemma}
\proofof{Lemma \ref{lemma_unique_maximal}}
We already know that $\hat\omega\neq 0$ and $f_*\hat\omega= e^{h_\top }\hat\omega,$ hence $e^{h_\top }$ is an eigenvalue of $f_*|_{\mathcal{B}^{p,q,d_s}}.$ Let $\omega\in\mathcal{B}^{p,q,d_s}$ be any other eigenvector with corresponding eigenvalue $\gamma e^{h_\top },$ such that $|\gamma|=1.$ As a consequence of Lemma \ref{lemma_abs_cont}, $\mathcal{\omega}$ is a continuous leafwise measure on each stable manifold and $|\mathcal{M}(\omega)|\leq C \mathcal{M}(\hat\omega).$   We check that the other hypotheses of Proposition \ref{prop_leafwise_measure} hold true. Firstly, an Anosov diffeomorphism is topologically transitive if and only if it is topologically mixing \cite[Theorem 5.10.3]{brin-stuck}. Next, setting $T=f^{-1}$, then $T$ is uniformly expanding on stable leaves of $f$. Moreover, defining $T^{*}\mathcal{M}(\omega)=\mathcal{M}(T^*\omega),$ we get $T^*\mathcal{M}(\hat\omega)=e^{h_\top }\mathcal{M}(\hat\omega)$ and  $T^*\mathcal{M}(\omega)=\gamma e^{h_\top }\mathcal{M}(\omega).$
It remains to prove that every open set $O$ contained in a stable manifold admits a point $x\in O$ with dense orbit. Let $x$ be a point in $O$. By topological transitivity, there exists a $y$ close to $x$, with dense orbit. By classical results (see for instance \cite[Proposition 5.9.1]{brin-stuck}), the local stable manifold centered at $x$, $W_\epsilon^s(x),$ and the local unstable manifold centered at $y$,  $W_\epsilon^u(y),$ intersect in exactly one point $z=[x,y]=W_\epsilon^s(x)\cap W_\epsilon^u(y).$ Accordingly, $z\in O$ and its orbit is dense. By Proposition \ref{prop_leafwise_measure}, we conclude that $\mathcal{M}(\omega)=c\mathcal{M}(\hat\omega),$ which in turn implies that $\omega= c\hat\omega$ and $\gamma=1.$

\qed
It remains to prove that the eigenvectors corresponding to the unique eigenvalue $e^{h_\top }$ defines a positive invariant measure and this is the measure of maximal entropy. Let $\hat\omega\in\mathcal{B}^{p,q,d_s}$ be, as above, the eigenvector for which $f_*\hat\omega=e^{h_\top }\hat\omega.$ Let $\hat t\in (\mathcal{B}^{p,q,d_s})'$ the unique element of the dual space of  $\mathcal{B}^{p,q,d_s}$ such that\footnote{We recall that the dual action of $f_*$ is the linear operator $f_*'\colon(\mathcal{B}^{p,q,d_s})'\rightarrow (\mathcal{B}^{p,q,d_s})'$ such that, for each $t\in(\mathcal{B}^{p,q,d_s})'$ and for each $\omega\in\mathcal{B}^{p,q,d_s},$ $f_*'(t)(\omega)=t(f_*\omega).$ In particular, $f_*'=(f^{-1})_*.$} $f_*'\hat t=e^{h_\top }\hat t$ and $\hat t (\hat\omega)=1.$
\begin{lemma}
	\label{lemma_inv_measure}
	The linear operator $\mu_{BM}=\hat t(\ \cdot \ \hat\omega),$ initially defined on $C^{p+q}(M)$  functions, extends to a bounded linear operator on continuous functions, i.e., it is a measure. In addition, for every $\psi\in C^0(M),$ $\mu_{BM}(\psi\circ f)=\mu_{BM}(\psi)$ and $\mu_{BM}$ is a positive probability measure.
\end{lemma}
\proofof{Lemma \ref{lemma_inv_measure}}
As above, we adapt the proofs of \cite[Lemma 6.1, Lemma 6.2]{gouezel-liverani08}.
Notice that, for every $\omega\in\mathcal{B}^{p,q,l}$  and for every function $\psi\in C^{p+q}(M),$ the product $\psi\omega\in\mathcal{B}^{p,q,l}.$ Moreover, for every $\omega\in\mathcal{B}^{p,q,d_s}$ it holds 
\begin{equation}
	\label{eq_bar_t}
	|\hat t(\omega)|\leq C\norm{\omega}_{0,p+q,d_s}.
\end{equation}
Indeed, using \eqref{eq_Lasota-Yorke_rem},
\begin{equation*}
	\begin{split}
		|\hat t(\omega)|=&e^{-pnh_\top }|(f_*')^{pn}\hat t(\omega)|=e^{-pnh_\top }|\hat t(f_*^{pn}\omega)|\leq e^{-pnh_\top }\norm{f_*^{pn}\omega}_{p,q,d_s}\leq\\\leq& C\lambda^{-n}\norm{\omega}_{p,q,d_s}+C\norm{\omega}_{0,p+q,d_s}
	\end{split}
\end{equation*}
Taking the limit for $n$ going to $\infty$ we get \eqref{eq_bar_t}. Next, by \eqref{eq_estimate}, we obtain that $\norm{\psi\hat\omega}_{0,p+q,d_s}\leq C\norm{\psi}_{C^0}.$ Thus, $$|\hat t(\psi\hat\omega)|\leq C\norm{\psi}_{C^0}$$ 
and $\hat t(\ \cdot\ \hat\omega)$ extends to a bounded operator on continuous functions.
Furthermore, 
\begin{equation*}
	\mu_{BM}(\psi\circ f)=\hat t(\psi\circ f\hat\omega)=e^{-h_\top }(f_*)'\hat t(\psi\circ f\hat\omega)=e^{-h_\top }\hat t(\psi f_*\hat\omega)=\hat t(\psi\hat\omega)=\mu_{BM}(\psi),
\end{equation*}
which proves that $\mu_{BM}$ is $f$-invariant. 

Let us prove that $\mu_{BM}$ is a positive measure. By the spectral decomposition \eqref{eq_sviluppo}, we can write 
$$\lim_{n\rightarrow+\infty}e^{-nh_\top }f_*^n\omega=\pi_1(\omega)\hat\omega,$$ 
whenever $\omega\in\mathcal{B}^{p,q,d_s},$ where $\pi_1$ is a linear form on $\mathcal{B}^{p,q,l}.$ Moreover, since $\hat t(\hat\omega)=1,$
$$\pi_1(\omega)=\pi_1(\omega)\hat t(\hat\omega)=\hat t(\pi_1(\omega)\hat\omega)=\hat t(\lim_{n\rightarrow+\infty}e^{-nh_\top }f_*^n\omega)=\lim_{n\rightarrow+\infty}e^{-nh_\top }\hat t(f_*^n\omega)=\hat t(\omega)$$
Accordingly, 
\begin{equation*}
	\lim_{n\rightarrow+\infty}e^{-nh_\top }f_*^n\omega=\hat t(\omega)\hat\omega
\end{equation*}
Given two $C^\infty$ functions $\phi,\psi\geq 0$ and a leaf $W\in\Sigma,$ we get 
\begin{equation}
	\label{eq_integral}
	\begin{split}
		0\leq\lim_{n\rightarrow+\infty}\int_W\langle\phi\omega_W,e^{-nh_\top }f_*^n(\psi\omega_\Sigma)\rangle\omega_W=\hat t(\psi\omega_\Sigma)\int_{W}\phi\mathcal{M}_W(\hat\omega)
	\end{split}
\end{equation}
Lemma \ref{lemma_abs_cont} shows that $\mathcal{M}(\omega)$ is a nonnegative and nonzero measure. Consequently, we can choose $W$ and $\phi>0$ so that last integral of \eqref{eq_integral} is strictly positive. This shows that, for every smooth $\psi\geq 0,$ $\hat t(\psi\omega_\Sigma)\geq0.$ Thus, for $\psi\geq 0,$ we get 
\begin{equation}
	\label{eq_stone}
	\begin{split}
		\mu_{BM}(\psi)=&\hat t(\psi\hat\omega)=\lim_{n\rightarrow+\infty}e^{-nh_\top }\hat t(\psi f_*^n\omega_\Sigma)=\lim_{n\rightarrow+\infty}e^{-nh_\top }\hat t(f_*^n(\psi\circ f^n\omega_\Sigma))=\\
		=&\lim_{n\rightarrow+\infty}\hat t(\phi\circ f^n\omega_\Sigma)\geq 0.
	\end{split}
\end{equation}
By the Stone-Weierstrass Theorem the inequality \eqref{eq_stone} can be extended to continuous functions. Thus, $\mu_{BM}$ is a positive measure. Finally, since $\hat t(\hat\omega)=1,$ we conclude that $\mu_{BM}$ is a probability measure. 

\qed
It remains to prove that $\mu_{BM}$ is the unique measure of maximal entropy. The following Lemma \ref{lemma_var_principle} is a streamlined version of \cite[Proposition 6.3]{gouezel-liverani08}.
\begin{lemma}
	\label{lemma_var_principle}
	Given $n\in\N,$ $x\in M$ and $\epsilon>0$, we denote by 
	$$B_n(x,\epsilon)=\{y\in M|\ d(f^{-i}(y),f^{-i}(x))<\epsilon,\text{ for }i=0,1,\dots,n-1\}$$
	the dynamical ball centered at $x$, of length $n$ and radius $\epsilon$.
	Then, there exist two constants $c_\epsilon, C_\epsilon>0$ such that 
	$$c_\epsilon e^{-nh_\top }\leq\mu_{BM}(B_n(x,\epsilon))\leq\overline{\mu_{BM}(B_n(x,\epsilon))}\leq C_\epsilon e^{-nh_\top }$$
\end{lemma}
\proofof{Lemma \ref{lemma_var_principle}}
Let $\phi\in C^{p+q}$ be a compactly supported function such that $0\leq\phi\leq 1$, $\text{supp}(\phi)\subseteq B_n(x,\epsilon)$ and $\phi|_{B_n(x,\epsilon/2)}=1.$
We show that 
$$c_\epsilon e^{-nh_\top }\leq\mu_{BM}(\phi)\leq C_\epsilon e^{-nh_\top },$$
which implies the lemma.
Let $W\in\Sigma,$ $\psi\in\Gamma_0^q(W),$ with $\norm{\psi}_{\Gamma_0^q(W)}\leq 1.$ Then, writing $\psi=\psi_0\omega_W$ for some $\psi\in C_0^q(W),$ and proceeding as in the proof of Lemma \ref{lemma_injectivity}, we get 
\begin{equation*}
	\begin{split}
		\int_W\!\!\!\langle\psi,\phi\hat\omega\rangle\omega_W=&e^{-nh_\top }\!\!\!\int_W\!\!\!\psi_0\phi\langle\omega_W,f_*^n\hat\omega\rangle\omega_W=e^{-nh_\top }\!\!\!\int_W\!\!\!\psi_0\phi\langle\omega_W,\hat\omega\rangle\circ f^{-n}(\lambda_n^s)^{-1}\omega_W=\\=&e^{-nh_\top }\sum_j\int_{W_j}\rho_j \psi_0\circ f^n\phi\circ f^n\langle\omega_W,\hat\omega\rangle\omega_{W_j}.
	\end{split}
\end{equation*}
The number of $W_j$ on which the integral is nonzero is uniformly bounded, because $\phi$ is supported in $B_n(x,\epsilon)$.
Accordingly, 
$$|\mu_{BM}(\phi)|=|\hat t(\phi\hat\omega)|\leq C\norm{\phi\hat\omega}_{0,q,d_s}\leq C \sup_{W,\psi}\left|\int_W\!\!\!\langle\psi,\phi\hat\omega\rangle\omega_W|\right|\leq C_\epsilon e^{-nh_\top }.$$

To estimate the other inequality, firstly notice that $\mathcal{M}(\hat\omega)$ gives strictly positive measure to any open piece of stable leaf. If it was not the case, there would be a ball $\mathcal{B}_{d_s}(x,\delta),$ contained in a stable manifold, such that  $\mathcal{M}(\hat\omega)(\mathcal{B}_{d_s}(x,\delta))=0.$ By invariance, $\mathcal{M}(\hat\omega)$ also assigns zero measure to $f^{-n}(\mathcal{B}_{d_s}(x,\delta)).$ By reasoning as in the proof of Lemma \ref{lemma_unique_maximal}, using the topological mixing property, $f^{-n}(\mathcal{B}_{d_s}(x,\delta))$ will meet a point $z\in W^s(x)$ with dense orbit.  Finally, since $\mathcal{M}(\hat\omega)$ is a continuous leafwise measure and since the orbit of $z$ is dense, we conclude that $\mathcal{M}(\hat\omega)=0,$ which contradicts Lemma \ref{lemma_injectivity}.

Next we prove that, if $W$ is a piece of stable manifold and $W$ contains a point $y$ with $d(x,y)<\epsilon/10$ and $d(y,\delta W)\geq\epsilon,$ then
$$\int_W\phi\mathcal{M}(\hat\omega)\geq c_\epsilon e^{-nh_\top }.$$
In fact, $f^{-n}(W)$ contains a  $d_s$-dimensional ball $\mathcal{B}$ of radius $\epsilon/10,$ which is contained in $f^{-n}(B_n(x,\epsilon/2)).$ Consequently, 
\begin{equation*}
	\begin{split}
		\int_W\phi\mathcal{M}(\hat\omega)=&\int_W\phi e^{-nh_\top }\mathcal{M}(f_*^n\hat\omega)=\\=&\int_{f^{-n}(W)}\phi\circ f^n e^{-nh_\top }\mathcal{M}(\hat\omega)\geq e^{-nh_\top }\int_{\mathcal{B}}\mathcal{M}(\hat\omega)\geq c_\epsilon e^{-nh_\top },
	\end{split}
\end{equation*}
where, in the last inequality, we used that $\mathcal{M}(\hat\omega)$ assigns positive measure to $\mathcal{B}.$ Notice that, by compactness, the constant $c_\epsilon$ does not depend on the leaf $W$.

Topological mixing also implies the following fact: for every $\delta>0$ there exists $M$, depending on $\epsilon$ and $\delta$, such that, for each $m\geq M,$ there is a constant $C$, which depends on $\epsilon,$ $\delta$ and $m$, such that, for every connected $W$ contained in a stable manifold, with $\diam(W)\geq 2\delta,$ it holds
\begin{equation}
	\label{eq_step2}
	\int_{f^{-m}(W)}\phi\mathcal{M}(\hat\omega)\geq Ce^{-mh_\top }.
\end{equation}

Finally, we prove that, for a full admissible stable leaf $\widetilde{W},$ contained in a stable manifold, there exists $C$, depending on $\epsilon$ and $\widetilde{W}$, such that, if $p$ is large enough, 
\begin{equation}
	\label{eq_step3}
	e^{-ph_\top }\int_{\widetilde{W}}f_*^p(\phi\mathcal{M}(\omega))\geq Ce^{-nh_\top }
\end{equation}
Let $L$ be a positive integer such that $M\leq L\leq p.$ Let $\{W_j\}$ be the subdivision of $f^{-p}(W)$ as described by Lemma \ref{lemma_foglie}. Then 
$$\int_{\widetilde{W}}f_*^p(\phi\mathcal{M}(\hat\omega))=\int_{f^{-p+L}(\widetilde{W})}f_*^{L}(\phi\mathcal{M}(\hat\omega))\geq C\sum_{j}\int_{B_j}f_*^L(\phi\mathcal{M}(\hat\omega))=C\sum_{j}\int_{f^{-L}(B_j)}\!\!\!\!\!\!\!\!\!\!\!\!\!\!\!\!\!\phi\mathcal{M}(\hat\omega)),$$
where $B_j$ is a $d_s$-dimensional ball of radius $2\delta$ contained in $f^L(W_j).$ To every integral on $f^{-L}(B_j)$ we can apply \eqref{eq_step2} and, since the sum  grows as $e^{ph_\top },$ we obtain that 
$$\int_{\widetilde{W}}f_*^p(\phi\mathcal{M}(\hat\omega))\geq Ce^{ph_\top }e^{-nh_\top },$$
which implies \eqref{eq_step3}.
Since $$\lim_{p\rightarrow+\infty}e^{-ph_\top }f_*^p(\phi\hat\omega)=\hat t(\phi\hat\omega),$$ considering a $W$ which satisfies \eqref{eq_step3}, we conclude that
\begin{equation*}
	\mu_{BM}(\phi)\geq\hat t(\phi\hat\omega)\geq\lim_{p\rightarrow+\infty} e^{-ph_\top }\int_{\widetilde{W}}f_*^p(\phi\mathcal{M}(\hat\omega))\geq c_\epsilon e^{-nh_\top }.
\end{equation*}

\qed
Next result is the last ingredient for the proof of Proposition \ref{prop_measure_of_maximal_entropy}.
\begin{lemma}
	\label{lemma_max_entropy}
	The measure $\mu_{BM}$ is the unique measure of maximal entropy, i.e., it is the Bowen-Margulis measure of the system.
\end{lemma}
\proofof{\ref{lemma_max_entropy}}
The Variational Principle \cite[Theorem 4.5.3]{Katok} states that
$$\sup\{h_\mu(f)|\ \mu\text{ is a }f\text{-invariant measure}\}=h_\top ,$$
where $h_\mu(f)$ denotes the $\mu$-metric entropy of $f.$
Next, notice that the spectral decomposition \eqref{eq_sviluppo} implies that $\mu_{BM}$ is mixing, hence ergodic. Thus, by the local entropy theorem \cite{brin_katok}, since $\mu_{BM}$ is a $f$-invariant, ergodic, probability measure, we obtain
\begin{equation*}
	\lim_{\epsilon\rightarrow+\infty}\limsup_{n\rightarrow+\infty}-\frac{1}{n}\log(\mu_{BM}(B_n(x,\epsilon)))\!=\!\!
	\lim_{\epsilon\rightarrow+\infty}\liminf_{n\rightarrow+\infty}-\frac{1}{n}\log(\mu_{BM}(B_n(x,\epsilon)))\!=\!h_{\mu_{BM}}(f)
\end{equation*}
By Lemma \ref{lemma_var_principle} we conclude that $h_\top =h_{\mu_{BM}}(f),$ hence $\mu_{BM}$ is a measure of maximal entropy. The proof that this $\mu_{BM}$ is the unique measure of maximal entropy easily follows repeating the proof of \cite[Theorem 20.3.7]{Katok}

\qed
This concludes the proof of Proposition \eqref{prop_measure_of_maximal_entropy}.

\qed

\section{Anisotropic de Rham cohomology and spectrum}
\label{sec_anisotropic_de_rham_cohomology_and_spectrum}
We recall that the space of $C^\infty$ (complex) differential forms $\Omega^l(M)$ with the exterior derivative $d\colon\Omega^l(M)\rightarrow\Omega^{l+1}(M)$ is a cochain complex, i.e., $d\circ d=0.$ $\omega\in\Omega^l(M)$ is closed if $d\omega=0,$  while $\omega$ is exact if there exists $u\in\Omega^{l-1}(M)$ such that $du=\omega.$ Since $d\circ d=0,$ exact forms are a vector subspace of closed forms.\footnote{Notice that, even if we consider complex valued differential forms, we are working in the category of differential vector bundles. On the other hand, since complex numbers are an algebraically closed field, we are sure to see all eigenvalues.} Accordingly, it makes sense to define the de Rham cohomology group with complex coefficients $H^l_{dR}(M,\C)=H^l_{dR}(M)$ as the quotient of closed $l$-forms w.r.t.\ exact $l$-forms.
The pushforward $f_*$ of a $C^\infty$-diffeomorphism $f$ on $M$ preserves closed and exact forms, hence it induces a linear map from the cohomology group $H^l_{dR}(M)$ to itself defined by $f_{\#}[\omega]=[f_*\omega].$

Next lemma give us the possibility to extend these ideas to our anisotropic Banach spaces. 
\begin{lemma}
	\label{lemma_ex_derivative}
	The exterior derivative extends to a continuous operator, denoted by the same letter, $d\colon\mathcal{B}^{p,q,l}\rightarrow\mathcal{B}^{p-1,q+1,l+1}$, while keeping the property that $d\circ d=0.$
\end{lemma}
\proof
Consider $h\in\Omega^l(M),$ $W\in\Sigma,$ $\phi\in\Gamma_0^{p+q-1,l+1}$ and $v_1,\dots,v_{p-1}\in\mathcal{V}^{p+q-1}(U(W)).$ If $W\subseteq\psi_i(U_i),$ then we can write, using coordinates, $h=h\circ\chi_i=\sum_{\overline{j}\in\mathcal{J}_l}h_{\overline{j}}dx_{\overline{j}}$ on $\psi_i(\mathcal{B}(0,3\rho)).$ Accordingly,
\begin{equation*}
	\begin{split}
		&\left|\int_W\langle\phi,L_{v_1}\dots L_{v_{p-1}}dh\rangle\omega_W\right|=\left|\int_W\langle\phi,L_{v_1}\dots L_{v_{p-1}}d\left(\sum_{\overline{j}\in\mathcal{J}_l}h_{\overline{j}}dx_{\overline{j}}\right)\rangle\omega_W\right|=\\
		=&\left|\sum_{\overline{j}\in\mathcal{J}_l}\sum_{s=1}^{\dim(M)}\int_W\langle\phi,L_{v_1}\dots L_{v_{p-1}}\de_{x_s}h_{\overline{j}}dx_{\overline{j}}\wedge dx_s\rangle\omega_W\right|\leq C\norm{\phi}_{\Gamma_0^{p+q-1,l+1}}\norm{h}_{p,q,l}.
	\end{split}
\end{equation*}
We conclude that $d$ extends to a bounded operator $d\colon\mathcal{B}^{p,q,l}\rightarrow\mathcal{B}^{p-1,q+1,l+1}.$ Let us prove $d\circ d=0.$ Recall that, given two differential forms $h,g\in\Omega^l(M),$ $h$ behaves as a current in the following way:
$$i(h)(g)=(h,g)=\int_M\langle h,g\rangle\omega_0=\int_Mh\wedge\star g,$$
where $\star$ denotes again the Hodge operator \cite{de_Cataldo}.
Consequently, given $h\in\mathcal{B}^{p,q,l}$ and a sequence $h_n\in\Omega^l(M)$ converging to $h$ in the $\norm{\cdot}_{p,q,l}$-norm, then $dh_n$ converges to $dh$ in $\mathcal{B}^{p-1,q+1,l+1}$ and, for each $g\in\Omega^{l+1}(M),$
$$i(dh)(g)=\lim_{n\rightarrow+\infty}i(dh_n)(g)=\lim_{n\rightarrow+\infty}i(h_n)(d^* g)=i(h)(d^* g),$$
where $d^*=(-1)^{\dim(M)(l+1)+1}\star d\star$ is the dual operator of $d$ (see \cite[Proposition 2.1.5]{de_Cataldo}). 
Accordingly, since $d^*\circ d^*=0$ on differential forms, we conclude that, for $h\in\mathcal{B}^{p,q,l},$
$i(d\circ dh)(g)=i(h)(d^*\circ d^* g)=0,$ for any $g\in\Omega^{l+2}(M), $ hence $d\circ d=0$.

\qed

We say, by analogy, that a current $\omega\in\mathcal{B}^{p,q,l}$ is \textit{closed} if $d\omega=0$, while it is \textit{exact} if there exists $u\in\mathcal{B}^{p+1,q-1,l-1}$ such that $du=\omega.$
As a consequence of Lemma \ref{lemma_ex_derivative}, we define the anisotropic De Rham cohomology $\widetilde{H}_{dR}^{p,q,l}(M)$ as the quotient of closed currents w.r.t.\ exact currents of $\mathcal{B}^{p,q,l}.$
Since $f_*\colon\mathcal{B}^{p,q,l}\rightarrow \mathcal{B}^{p,q,l}$ and $d\colon\mathcal{B}^{p+1,q-1,l-1}\rightarrow\mathcal{B}^{p,q,l}$ are both continuous linear operators, $f_*$ sends closed currents to closed currents and exact currents to exact currents. Consequently, it induces a linear map on $\widetilde{H}_{dR}^{p,q,l}(M)$ such that $f_{\#}[\omega]=[f_*\omega].$ 

Next proposition relates the spectrum of $f_*$ acting on anisotropic Banach spaces and the spectrum $f_{\#}$ on anisotropic de Rham cohomology. We just consider the spectrum of $f_*$  contained in the set $\{z\in\C\mid|z|>\nu\},$ where $\nu$ is the bound defined in \eqref{eq_nu}. Accordingly, we only take care of the discrete spectrum. Notice that this is the higher-dimensional generalization of \cite[Lemma 5.22]{butterley-kiamari-liverani}.
\begin{proposition}
	\label{prop_spectrum_1}
	$$\sigma(f_*|_{\mathcal{B}^{p,q,l}})\cap\{|z|>\nu\}\subseteq\left[\sigma(f_*|_{\mathcal{B}^{p+1,q-1,l-1}})\cup\sigma(f_\#|_{\widetilde{H}_{dR}^{p,q,l}})\cup\sigma(f_*|_{\mathcal{B}^{p-1,q+1,l+1}})\right]\cap\{|z|>\nu\}.$$
\end{proposition}
\proof
Let $\omega\in\mathcal{B}^{p,q,l}$ be an eigenvector of the pushforward operator $f_*$ corresponding to the eigenvalue $\mu,$ with $|\mu|>\nu,$ that is $f_*\omega=\mu\omega.$ If $\omega$ is not closed, then  $f_*d\omega=df_*\omega=\mu d\omega$, i.e., $d\omega\neq0$ is an eigenvector for $f_*$ in $\mathcal{B}^{p-1,q+1,l+1}.$  This proves that $\mu\in\sigma(f_*|_{\mathcal{B}^{p-1,q+1,l+1}}).$ On the other hand, if $\omega$ is closed, i.e., $d\omega=0$ we need to distinguish two cases. If $\omega$ is not exact, then it defines a nontrivial cohomology class $[\omega]\in\widetilde{H}_{dR}^{p,q,l}$ and, by definition, 
$f_\#[\omega]=[f_*\omega]=[\mu\omega]=\mu[\omega],$ hence 
$\mu\in\sigma(f_\#|_{\widetilde{H}_{dR}^l}).$
Finally, if $\omega$ is exact, then there exists a nonclosed $q\in\mathcal{B}^{p+1,q-1,l-1}$ such that $\omega=dq.$ It follows that $f_*dq=df_*q=\mu dq.$ This is not enough to conclude that $\mu\in\sigma(f_*|_{\mathcal{B}^{p+1,q-1,l-1}}),$ because it only gives that $d(f_*q-\mu q)=0,$ hence $f_*q=\mu q+v$ with $v\in\mathcal{B}^{p+1,q-1,l-1}$ closed. On the other hand,
if the operator $f_*-\mu\id$ was invertible on closed $(l-1)$-currents of $\mathcal{B}^{p+1,q-1,l-1},$ then there would exist a unique closed $u=(f_*-\mu\id)^{-1}v$ and $u$ would be closed. But then $f_* (q-u)=\mu (q-u)$ and $q-u\neq 0$, since $q$ is not closed and $u$ is closed. Thus, $\mu\in \sigma(f_*|_{\mathcal{B}^{p+1,q-1,l-1}}).$ On the contrary, if $f_*-\mu\id$ is not invertible on closed currents, then $\mu$ is again an eigenvalue of $f_*$ on $\mathcal{B}^{p+1,q-1,l-1}$, with the additional information that the related eigenvector is closed.

\qed

We can also prove that the spectrum of the action on anisotropic cohomology, outside the ball of radius $\nu$, is included into the spectrum of $f_*.$
\begin{proposition}
	\label{prop_spectrum_2}
	$$\sigma(f_\#|_{\widetilde{H}_{dR}^{p,q,l}})\cap\{|z|>\nu\}\subseteq \sigma(f_*|_{\mathcal{B}^{p,q,l}})\cap\{|z|>\nu\} $$
\end{proposition}
\proof
If $[\omega]\in \widetilde{H}_{dR}^{p,q,l}$ is an eigenvector of $f_\#$ of modulus greater than $\nu$, then $f_\#[\omega]=\mu[\omega]$ and $[\omega]\neq 0$. If $f_*\omega=\mu\omega,$ then $\mu\in\sigma(f_*|_{\mathcal{B}^{p,q,l}}).$ On the other hand,
there exists $u\in\mathcal{B}^{p+1,q-1,l-1}$ such that
$$f_*\omega=\mu\omega+du.$$
We can proceed, as above, looking for a current $\omega'=\omega+du'$, with $u'\in\mathcal{B}^{p+1,q-1,l-1},$ so that $[\omega]=[\omega'],$ and with $f_*\omega'=\mu\omega'.$ Last equality means 
$$\mu\omega+\mu du'=\mu\omega'=f_*\omega'=f_*\omega+f_*du'=\mu\omega+du+f_*du',$$
hence $(f_*-\mu \id)du'=-du.$
If $(f_*-\mu \id)$ is invertible on exact currents of $\mathcal{B}^{p,q,l}$, then the desired $u'=(f_*-\mu \id)^{-1}u$ and $\mu\in\sigma(f_*|_{\mathcal{B}^{p,q,l}}).$ Conversely, if $(f_*-\mu \id)$ is not invertible,
on exact currents of $\mathcal{B}^{p,q,l}$,  there exists $d\bar{u}\in\mathcal{B}^{p,q,l}$ such that $f_*d\bar{u}=\mu d\bar{u}$, because the spectrum is discrete in $\{|z|>\nu\}.$
\qed

In the particular case $l=d_s$ we can identify the spectrum of $f_*$ with the spectrum of the action on anisotropic cohomology out of the ball of radius $\lambda^{-1}e^{h_\top }.$
\begin{corollary}
	\label{cor_spectrum}
	$$\sigma(f_*|_{\mathcal{B}^{p,q,d_s}})\cap\{z\in\mathbb{C}: |z|>\lambda^{-1}e^{h_\top }\}=\sigma(f_{\#}|_{\widetilde{H}_{d_R}^{p,q,d_s}})\cap\{z\in\mathbb{C}: |z|>\lambda^{-1}e^{h_\top }\}$$
\end{corollary}
\proof
It follows by Proposition \ref{prop_spectrum_1}, Proposition \ref{prop_spectrum_2} and Corollary \ref{cor_spectral_radius}, because the spectral radius of $f_*$ on $\mathcal{B}^{p+1,q-1,l-1}$ and on $\mathcal{B}^{p-1,q+1,l+1}$ is bounded by $\lambda^{-1}e^{h_\top }.$

\begin{center}
	\begin{figure}[h!]
		\includegraphics[scale=0.5]{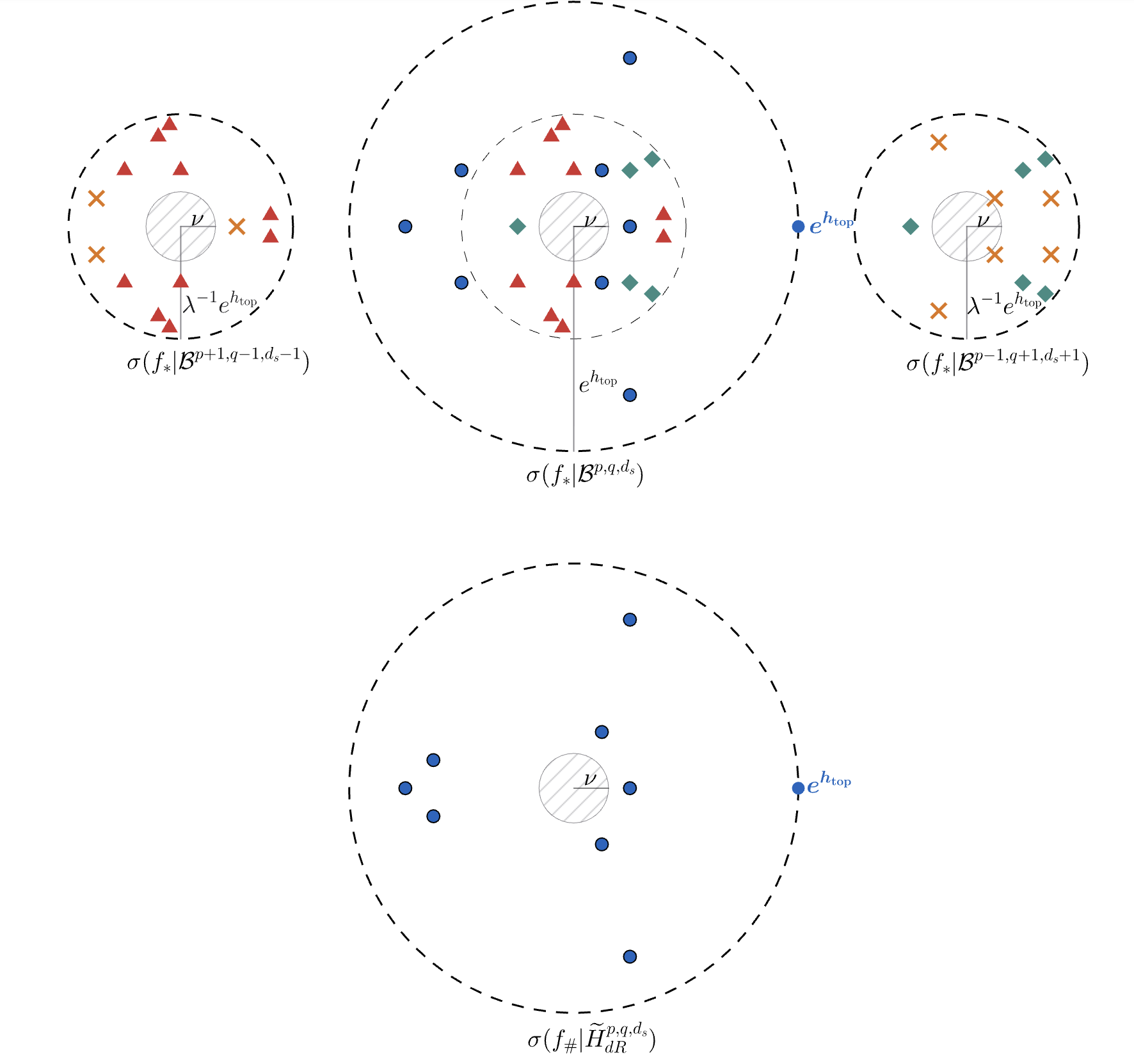}
		\caption{Example of spectra with $l=d_s.$}
			\label{pic_spectrum}
	\end{figure}
\end{center}

\qed
Figure \ref{pic_spectrum} below is a possible representation of the relation between spectra proved in Proposition \ref{prop_spectrum_1}, Proposition \ref{prop_spectrum_2} and Corollary \ref{cor_spectrum}, in the particular case $l=d_s$.The red triangles $\color{miorosso}\blacktriangle$ represent the eigenvalues of  $f_*|_{\mathcal{B}^{p+1,q-1,d_s-1}}$ that are also eigenvalues of $f_*|_{\mathcal{B}^{p,q,d_s}}.$ They correspond to eigenvectors in $\mathcal{B}^{p+1,q-1,d_s-1}$ that are closed.  The green diamonds $\color{mioverde}\blacklozenge$ are the eigenvalues of $f_*|_{\mathcal{B}^{p-1,q+1,d_s+1}}$ that are also eigenvalues of $f_*|_{\mathcal{B}^{p,q,d_s}}.$ They correspond to eigenvectors in $\mathcal{B}^{p-1,q+1,d_s+1}$ that are not exact. The blue circles $\color{mioblu}\bullet$ are the eigenvalues of $f_\#$ in the anisotropic cohomology $\widetilde{H}^{p,q,d_s}_{dR}$ that are also eigenvalues in $f_*|_{\mathcal{B}^{p,q,d_s}}.$ Notice that eigenvalues of $f_*|_{\mathcal{B}^{p,q,d_s}}$ of modulus larger than $\lambda^{-1}e^{h_\top}$ are forced to be in cohomology. Finally, the symbol $\color{mioarancio}\times$ represents eigenvalues of $f_*|_{\mathcal{B}^{p+1,q-1,d_s-1}}$ or $f_*|_{\mathcal{B}^{p-1,q+1,d_s+1}},$ that are not eigenvalues of $f_*|_{\mathcal{B}^{p,q,d_s}}.$  

\subsection{Connection with the standard de Rham cohomology}
\label{section_connection_de_Rham}
In view of Corollary \ref{cor_spectrum}, it makes sense to study the spectrum of $f_{\#}$ on the anisotropic de Rham cohomology group $\widetilde{H}_{dR}^{p,q,d_s}(M).$ The first attempt to solve this problem could be trying to show that anisotropic de Rham cohomology is actually isomorphic to the standard de Rham cohomology recalled at the beginning of Section \ref{sec_anisotropic_de_rham_cohomology_and_spectrum}.
In fact, it is well-known by classical results of G. de Rham \cite{de_rham_book} that the de Rham cohomology for currents is isomorphic to the de Rham cohomology for differential forms. On the other hand we are working with a linear subspace $\mathcal{B}^{p,q,l}$ of the dual space of $C^{p+q}$ differential forms  $(\Omega^l_{p+q}(M))^\star$(see Lemma \ref{lemma_correnti}) , which is in turn a subspace of the space of currents $(\Omega^l(M))^\star,$ i.e., the dual of $C^\infty$ differential forms. Therefore, we have fewer closed currents and fewer exact currents than the full space of currents and, a priori, there is no relation between our cohomology and the standard cohomology.

The authors of \cite[Section 5.7]{butterley-kiamari-liverani} showed that this isomorphism exists when $d_s=1$ and it is enough to prove our main result for Anosov diffeomorphisms of the 2-torus. Their strategy also works for our Anosov diffeomorphisms on higher dimensional manifolds whenever $d_s=1$, but unfortunately, the extension to other cases requires a bit of work. A motivation of the problem is given in Remark \ref{rem_lemma_Poincare}.

In order to overcome this obstacle, we firstly need to introduce an intermediate version of our anisotropic Banach spaces. 
\begin{definition}
	\label{def_new_anisotropic}
	Let $\omega\in\Omega^l(M)$ be a $C^\infty$ differential form and let $p,q\in\N.$ We define the following anisotropic norm
	$$|\omega|_{p,q,l}=\norm{\omega}_{p,q,l}+\norm{d\omega}_{p,q,l+1},$$
	where $\norm{\cdot}_{p,q,l}$ is the norm of Definition \ref{def_norm}. Let us denote by $\mathcal{C}^{p,q,l}=\overline{\Omega^{l}(M)}^{|\cdot|_{p,q,l}}$ the closure of the space of $l$-forms w.r.t.\ this norm.
\end{definition}
The following proposition collects all the properties we need about this new anisotropic Banach spaces.
\begin{proposition}
	\label{prop_new_anisotropic}
	The following properties hold 
	\begin{itemize}
		\item[a)] $\mathcal{B}^{p+j,q-j,l}\subseteq\mathcal{C}^{p,q,l}\subseteq\mathcal{B}^{p,q,l}$ for any $p,q\in\N$ and for each $j=1,\dots, q;$
		\item[b)] $\omega\in\mathcal{C}^{p,q,l}$ if and only if $\omega\in\mathcal{B}^{p,q,l}$ and $d\omega\in\mathcal{B}^{p,q,l+1};$
		\item[c)] The exterior derivative extends to a bounded linear operator $d\colon\mathcal{C}^{p,q,l}\rightarrow\mathcal{C}^{p,q,l+1}$ and $d\circ d=0.$
		\item[d)] $f_*$ extends to a bounded linear operator $f_*\colon\mathcal{C}^{p,q,l}\rightarrow\mathcal{C}^{p,q,l}.$ The spectral radius $\rho(f_*|_{\mathcal{C}^{p,q,l}})\leq\lambda^{-|d_s-l|}e^{h_\top },$ while, for $p$ and $q$ large enough, the essential spectral radius $\rho_{ess}(f_*|_{\mathcal{C}^{p,q,l}})\leq\nu,$ where $\nu$ is the same as in \eqref{eq_nu}.
	\end{itemize}
\end{proposition}
\proof
By definition $\norm{\omega}_{p,q,l}\leq\norm{\omega}_{p,q,l}+\norm{d\omega}_{p,q,l+1}=|\omega|_{p,q,l},$ hence $\mathcal{C}^{p,q,l}\subseteq\mathcal{B}^{p,q,l}.$ Next,
$|\omega|_{p,q,l}=\norm{\omega}_{p,q,l}+\norm{d\omega}_{p,q,l+1}\leq\norm{\omega}_{p,q,l}+\norm{d}\norm{\omega}_{p+1,q-1,l}\leq C\norm{\omega}_{p+1,q-1,l},$ where the first inequality follows by the continuity of $d\colon\mathcal{B}^{p+1,q-1,l}\rightarrow\mathcal{B}^{p,q,l+1}$ (see Lemma \ref{lemma_ex_derivative}), while the second one is a consequence of the inclusion $\mathcal{B}^{p+1,q-1,l}\subseteq\mathcal{B}^{p,q,l}$  (see Remark \ref{rem_properties}). Accordingly, $\mathcal{B}^{p+1,q-1,l}\subseteq\mathcal{C}^{p,q,l}\subseteq\mathcal{B}^{p,q,l}$  and, using again Remark \ref{rem_properties}, we get \emph{a)}. \emph{b)} is a trivial consequence of Definition \ref{def_new_anisotropic}. To prove \emph{c)} notice that 
$|d\omega|_{p,q,l+1}=\norm{d\omega}_{p,q,l+1}\leq\norm{\omega}_{p,q,l}+\norm{d\omega}_{p,q,l+1}=|\omega|_{p,q,l},$ hence $d\colon\mathcal{C}^{p,q,l}\rightarrow\mathcal{C}^{p,q,l+1}$ is bounded. $d\circ d=0$ easily follows by the inclusions \emph{a)} and by Lemma \ref{lemma_ex_derivative}. Finally, let us prove \emph{d)}. The statement about the spectral radius is a consequence of the inclusion $\mathcal{C}^{p,q,l}\subseteq\mathcal{B}^{p,q,l}$ and Corollary \ref{cor_spectral_radius}. For the second part, we want to apply again Hennion's Theorem \ref{thm_hennion}.  It is not difficult to show that the inclusion $\mathcal{C}^{p,q,l}\hookrightarrow \mathcal{C}^{p-1,q+1,l}$ is compact. In fact, let $(\omega_n)_{n\in\N}$ be a sequence in $\mathcal{C}^{p,q,l}$ such that $|\omega_n|_{p,q,l}\leq 1. $ Then, $\norm{\omega_n}_{p,q,l}\leq 1$ and $\norm{d\omega_n}_{p,q,l+1}\leq 1.$ Since $\mathcal{B}^{p,q,l}\hookrightarrow\mathcal{B}^{p-1,q+1,l}$ and $\mathcal{B}^{p,q,l+1}\hookrightarrow\mathcal{B}^{p-1,q+1,l+1}$ are compact, there exists a subsequence $(\omega_{n_m})_{m\in\N}$ such that $\omega_{n_m}\rightarrow\tilde{\omega} $ in $\mathcal{B}^{p-1,q+1,l}$ and $d\omega_{n_m}\rightarrow{}d\tilde{\omega}$ in $\mathcal{B}^{p-1,q+1,l+1},$ as $m$ goes to $\infty$.
Moreover, by using the Lasota-Yorke inequality \eqref{Lasota-Yorke_4}, we obtain
\begin{equation*}
	\begin{split}
		&|f_*^n\omega|_{p,q,l}=\norm{f_*^n\omega}_{p,q,l}+\norm{f_*^nd\omega}_{p,q,l+1}\leq\\\leq& C\lambda^{-n(|d_s-l|+\min\{p,q\})}e^{nh_\top }\norm{\omega}_{p,q,l}+C\lambda^{-n|d_s-l|}e^{nh_\top }\norm{\omega}_{p-1,q+1,l}+\\&+C\lambda^{-n(|d_s-l-1|+\min\{p,q\})}e^{nh_\top }\norm{d\omega}_{p,q,l+1}+C\lambda^{-n|d_s-l-1|}e^{nh_\top }\norm{d\omega}_{p-1,q+1,l+1}\leq\\\leq& C\max\{\lambda^{-n|d_s-l|},\lambda^{-n|d_s-l-1|}\}\lambda^{-n\min\{p,q\}}e^{nh_\top }|\omega|_{p,q,l}+\\&+C\max\{\lambda^{-n|d_s-l|},\lambda^{-n|d_s-l-1|}\}e^{nh_\top }|\omega|_{p-1,q+1,l},
	\end{split}
\end{equation*}
which is a Lasota-Yorke inequality for the norm $|\cdot|_{p,q,l}.$
We conclude that the essential spectral radius is bounded by $\max\{\lambda^{-|d_s-l|},\lambda^{-|d_s-l-1|}\}\lambda^{-\min\{p,q\}}e^{h_\top },$ which is smaller than $\nu$ for $p$ and $q$ large enough.

\qed

Property \emph{c)} of Proposition \ref{prop_new_anisotropic} gives the following cochain complex
$$0\xrightarrow{d}\mathcal{C}^{p,q,0}\xrightarrow{d}\mathcal{C}^{p,q,1}\xrightarrow{d}\mathcal{C}^{p,q,2}\rightarrow\dots\xrightarrow{d}\mathcal{C}^{p,q,\dim(M)-1}\xrightarrow{d}\mathcal{C}^{p,q,\dim(M)}\xrightarrow{d}0,$$
hence we can define $\bar{H}^{p,q,l}_{dR}(M)$ as the quotient of closed currents w.r.t.\ exact currents of $\mathcal{C}^{p,q,l}.$ Collecting information about the spectrum of $f_*$ acting on different versions of anisotropic Banach spaces, we obtain the following result.\footnote{See again Figure \ref{pic_spectrum} for a graphic representation of the spectra of $f_*$ and $f_\#$. The following Corollary \ref{cor_spectrum2} tells us that there is no substantial difference on the spectra considering $\mathcal{C}^{p,q,l}$ in place of $\mathcal{B}^{p,q,l}.$}
\begin{corollary}
	\label{cor_spectrum2}
	Let $p,q\in\N$ be large enough. Then
	\begin{enumerate}
		\item $\sigma(f_*|_{\mathcal{C}^{p,q,l}})\cap\{|z|>\nu\}=\sigma(f_*|_{\mathcal{B}^{p\pm i,q\mp i,l}})\cap\{|z|>\nu\}$
		for any $i=0,\dots,\dim(M);$
		\item $\sigma(f_*|_{\mathcal{C}^{p,q,l}})\cap\{|z|>\nu\}\subseteq\left[\sigma(f_*|_{\mathcal{C}^{p,q,l-1}})\cup\sigma(f_\#|_{\bar{H}_{dR}^{p,q,l}})\cup\sigma(f_*|_{\mathcal{C}^{p,q,l+1}})\right]\cap\{|z|>\nu\};$
		\item $\sigma(f_\#|_{\bar{H}_{dR}^{p,q,l}})\cap\{|z|>\nu\}\subseteq \sigma(f_*|_{\mathcal{C}^{p,q,l}})\cap\{|z|>\nu\};$
		\item $\sigma(f_*|_{\mathcal{C}^{p,q,d_s}})\cap\{z\in\mathbb{C}\mid|z|>\lambda^{-1}e^{h_\top }\}=\sigma(f_{\#}|_{\bar{H}_{dR}^{p,q,d_s}})\cap\{z\in\mathbb{C}\mid|z|>\lambda^{-1}e^{h_\top }\}$
	\end{enumerate} 
\end{corollary}
\proof
Equality \emph{1.} is again a consequence of the independence of the spectrum by the choice of the anisotropic space, i.e., Lemma \ref{lemma_baladi_tsujii} with  $\mathcal{B}_0=\Omega^l(M),$ $\mathcal{B}_1=\mathcal{C}^{p,q,l}$, $\mathcal{B}_2=\mathcal{B}^{p\pm i,q\mp i,l}$ and $\mathcal{B}=\mathcal{B}^{\max\{p,p\pm i\},\min\{q,q\mp i\},l}.$ \emph{2.}, resp.\ \emph{3.} can be proved by repeating the proof of Proposition \ref{prop_spectrum_1}, resp.\ Proposition \ref{prop_spectrum_2}, with $\mathcal{B}^{p,q,l}$ replaced by $\mathcal{C}^{p,q,l}$ and $\widetilde{H}_{dR}^{p,q,d_s}$ replaced by $\bar{H}_{dR}^{p,q,d_s}.$ Finally, \emph{2.}, \emph{3.} and \emph{d)} of Proposition \ref{prop_new_anisotropic} imply \emph{4.}

\qed

\begin{remark}
	Notice that one may directly study the action of $f_*$  on $\mathcal{C}^{p,q,l},$ without considering the original anisotropic Banach spaces $\mathcal{B}^{p,q,l}.$ There are several reasons that have led us to our choice. In fact, proofs of Lasota-Yorke inequalities (Lemma \ref{lemma_Lasota_Yorke}) and compact inclusion (Lemma \ref{inclusione_compatta}), as well as the inclusion into currents (Lemma \ref{lemma_correnti}) for $\mathcal{B}^{p,q,l}$ turn out to be simpler from a technical point of view. Secondly, these spaces have been largely investigated in recent years \cite{gouezel-liverani,gouezel-liverani08,giulietti-liverani-pollicott} and we have picked up some ideas from the literature. Lastly, we found the issue when we started treating anisotropic  cohomology for $d_s>1$ and we discovered that without our trick the proof of the isomorphism Theorem  \ref{iso-anisotropic} does not work with $\widetilde{H}_{dR}^{p,q,l},$ that is the cohomology obtained with $\mathcal{B}^{p,q,l}$ (see Remark \ref{rem_lemma_Poincare} below).
\end{remark}

We are now ready to prove the isomorphism between the anisotropic de Rham cohomology $\bar{H}^{p,q,l}_{dR}$ and the standard de Rham cohomology. It is well known, in the fields of algebraic topology and differential geometry, that the standard de Rham cohomology is isomorphic to the \v{C}ech cohomology. Indeed, denoting with $H_{\check{C}}^*(M)$ the \v{C}ech cohomology with complex coefficients, whose definition is recalled just after Theorem \ref{iso-anisotropic}, the following theorem holds.
\begin{theorem}[De Rham isomorphism theorem]
	\label{iso_derham}
	There exists a natural isomorphism between the standard de Rham cohomology and the \v{C}ech cohomology
	$$H_{dR}^*(M)\cong H_{\check{C}}^*(M)$$
\end{theorem}
There are several different proofs of above result. The most elegant one, due to André Weil \cite{weil} (see also \cite{morita}, \cite{bott}), inspires our proof of the following isomorphism theorem. In fact, the careful reader may also be able to reconstruct the proof of Theorem \ref{iso_derham} from the one of Theorem \ref{iso-anisotropic}.
\begin{theorem}
	\label{iso-anisotropic}
	Let $p,q\in\N$ be large enough. There exists a natural isomorphism between the anisotropic de Rham cohomology and the \v{C}ech cohomology
	$$\bar{H}_{dR}^{p,q,*}(M)\cong H_{\check{C}}^*(M).$$
	By Theorem \ref{iso_derham}, we obtain 
	$$\bar{H}_{dR}^{p,q,*}(M)\cong H_{dR}^*(M).$$
\end{theorem}

Before proving Theorem \ref{iso-anisotropic}, we recall the basic facts about \v{C}ech cohomology (for a complete discussion of the topic see \cite{hatcher,morita}). 

Let $\mathcal{U}=\{U_a\}_{a\in\mathcal{A}}$ be a contractible open covering of the manifold $M,$ i.e., we suppose that every finite nonempty intersection
$$U_{a_1}\cap U_{a_2}\cap\dots\cap U_{a_n}\neq\emptyset$$ 
is contractible, i.e., homotopic to a point.
We denote by $(a_0,\dots,a_k):=U_{a_0}\cap\dots\cap U_{a_k}.$
Let $\check{C}_k(M,\mathcal{U})$ be the complex vector space generated by elements $(a_0,\dots,a_k)\neq\emptyset;$ the elements in $\check{C}_k(M,\mathcal{U})$ are called (\v{C}ech) $k$-chains. 
A (\v{C}ech) $k-$cochain $c$ is an element of the dual of $\check{C}_k(M,\mathcal{U})$ such that, for every permutation $\sigma$ of the indexes $\{0,\dots,k\},$ 
$$c(a_0,\dots,a_k)=\text{sgn}(\sigma)c(a_{\sigma(0)},\dots,a_{\sigma(k)}).$$
Let $\check{C}^k(M,\mathcal{U})$ be the complex vector space of all $k-$cochains.
We define the coboundary operator 
$$\delta\colon\check{C}^k(M,\mathcal{U})\rightarrow \check{C}^{k+1}(M,\mathcal{U})$$
such that
$$(\delta c)(a_0,\dots,a_{k+1})=\sum_{j=0}^{k+1} (-1)^j c(a_0,\dots, a_{j-1},a_{j+1},\dots,a_{k+1}).$$
A straightforward computation (see Lemma \ref{lem_equalities}) shows that $\delta^2=\delta\circ\delta=0,$ hence the couple $(\check{C}^*(M,\mathcal{U}),\delta)$ is a cochain complex. Recall that a $k$-cochain $c\in\check{C}^k(M,\mathcal{U})$ is a $k$-cocycle if $\delta(c)=0$, and it is a $k$-coboundary if there exists $c'\in\check{C}^{k-1}(M,\mathcal{U})$ such that $\delta(c')=c.$ As always, every $k$-coboundary is a $k$-cocycle, because of $\delta^2=0$. Therefore, we define the \v{C}ech $k$-cohomology group $H^k_{\check{C}}(M,\mathcal{U})$ as the quotient of $k$-cocycles with respect to $k$-coboundaries. A priori, this new cohomology depends on the covering $\mathcal{U}$ of the manifold $M.$ It will be clear from the proof of Theorem \ref{iso-anisotropic} that it is actually independent of that choice and it make sense to write 
$H_{\check{C}}(M)$ in place of $H_{\check{C}}(M,\mathcal{U}).$

Without loss of generality we can assume that the sets $\{V_i=\psi_i(U_i)\}_{i=1}^m$ is a contractible covering of the manifold $M.$ Let $\phi_i$ be a smooth partition of unity subordinated to the covering $\{V_i\}_{i=1}^m.$ We can also suppose that $\{\text{int}(\text{supp}(\phi_k))\}_{k=1}^m$, where $\text{int}(\text{supp}(\phi_k))$ is the interior of the support of $\phi_k$, is a contractible open covering.
If $\omega$ is a differential form on $M$, then the restriction to $V_i$ is well-defined. On the other hand, in our case $\omega\in\mathcal{C}^{p,q,l}$ is a current and the restriction of $\omega$ to the subset $V_i$ can be defined as $\omega\phi_i.$ A straightforward computation shows that $\omega\phi_i\in\mathcal{C}^{p,q,l}.$

We are now ready to present the proof of Theorem \ref{iso-anisotropic}.\\
\proofof{Theorem \ref{iso-anisotropic}}
Let us outline the main steps of the proof. We have introduced above the two cochain complexes $(\mathcal{C}^{p,q,*}(M),d)$ and $(\check{C}^*(M,\mathcal{U}),\delta)$, that define the anisotropic de Rham cohomology $\bar{H}_{dR}^{p,q,*}(M)$ and the \v{C}ech $k$-cohomology $H^k_{\check{C}}(M,\mathcal{U}),$ respectively. In  the first part of the proof, we build a bridge between $(\mathcal{C}^{p,q,*}(M),d)$ and $(\check{C}^*(M,\mathcal{U}),\delta)$, made of others cochain complexes, that can be visualized as the commutative diagram in Figure \ref{fig_comm_diagram1}. The second part is a zigzag argument on the commutative diagram of Figure \ref{fig_comm_diagram1} that connects every closed non-exact current $\omega\in\mathcal{C}^{p,q,l}$ (in the first column) to a closed non-exact cochain $c\in\check{C}^l(M,\mathcal{U})$ (in the bottom row), and viceversa. Accordingly, we define $\Phi([\omega])=[c].$ In the last part, we show that $\Phi$ is well-defined and invertible, hence it is the desired isomorphism.

\noindent\textit{STEP 1: New cochain complexes and the commutative diagram.}\\
Let us introduce the following notations. $\mathcal{C}^{p,q,l}(M)$ denotes the above anisotropic Banach space $\mathcal{C}^{p,q,l}$. Let $\mathcal{C}^{p,q,l,k}(\mathcal{U})$ be the vector space of linear functions $$\omega\colon\check{C}_k(M,\mathcal{U})\rightarrow\mathcal{C}^{p,q,l},$$ such that, for every permutation $\sigma$ of the set $\{0,\dots,k\}$, 
\begin{equation}
	\label{eq_permutation}
	\omega(a_{\sigma(0)},\dots,a_{\sigma(k)})=\text{sgn}(\sigma) \omega(a_0,\dots,a_k).
\end{equation} 
We underline that $\omega\in\mathcal{C}^{p,q,l,k}(\mathcal{U})$ can be interpreted as a \v{C}ech cochain with values in $\mathcal{C}^{p,q,l}$. In addition, the operator $\delta$ and $d$ can be extended to $\mathcal{C}^{p,q,l,k}$ by defining:
\begin{align*}
	\bar{\delta}\colon \mathcal{C}^{p,q,l,k}(\mathcal{U})&\rightarrow &&\mathcal{C}^{p,q,l,k+1}(\mathcal{U})\\
	\omega&\mapsto&&\bar{\delta}\omega\colon
	\check{C}_{k+1}(M,\mathcal{U})\longrightarrow\mathcal{C}^{p,q,l}\\
	& &&(\bar{\delta}\omega)(a_0,a_1\dots,a_{k+1})=\sum_{j=0}^{k+1}(-1)^j\omega (a_0,\dots,a_{j-1},a_{j+1},\dots,a_{k+1}),\\
	\bar{d}\colon\mathcal{C}^{p,q,l,k}(\mathcal{U})&\rightarrow &&\mathcal{C}^{p,q,l+1,k}(\mathcal{U})\\\omega&\mapsto&&\bar{d}\omega\colon
	\check{C}_k(M,\mathcal{U})\longrightarrow\mathcal{C}^{p,q,l+1}\\
	& &&(\bar d\omega)(a_0,\dots,a_k)=d(\omega(a_0,\dots,a_k)).
\end{align*}
Finally, we define the injections
\begin{align*}
	i\colon\mathcal{C}^{p,q,l}\rightarrow&\mathcal{C}^{p,q,l,0}\\
	\omega\mapsto& i(\omega)\colon \check{C}_0(M,\mathcal{U})\rightarrow\mathcal{C}^{p,q,l}\\
	& i(\omega)(a_0)=\omega,\\
	j\colon \check{C}^{k}(M,\mathcal{U})\rightarrow&\mathcal{C}^{p,q,0,k}\\
	c\mapsto& j(c)\colon \check{C}_k(M,\mathcal{U})\rightarrow\mathcal{C}^{p,q,0}\\
	&j(c)(a_0,\dots,a_k)=c(a_0,\dots,a_k).
\end{align*}
The previous linear maps $d,\ \bar{d},\ \delta,\ \bar{\delta},\ i $ and  $j$ define the diagram of Figure \ref{fig_comm_diagram1}. The following Lemma \ref{lem_equalities} shows that the diagram is commutative and that all rows and all columns are cochain complexes.  Notice that the first column is the anisotropic de Rham cochain complex $(\mathcal{C}^{p,q,*}(M),d)$, while the bottom row represents the \v{C}ech cochain complex $(\check{C}^*(M,\mathcal{U}),\delta)$.
\begin{lemma}
	\label{lem_equalities}
	The following equalities hold true.
	$$\delta\circ\delta=\bar\delta\circ\bar\delta=d\circ d=\bar d\circ \bar d=\bar\delta\circ i=\bar d\circ j=0,$$
	$$\bar d\circ\bar\delta=\bar\delta\circ \bar d,\ j\circ\delta=\bar\delta\circ j,\ i\circ d=\bar d\circ i.$$
\end{lemma}
\proofof{Lemma \ref{lem_equalities}}
For every $(a_0,\dots,a_{k+2})\in \check{C}_{k+2}(M,\mathcal{U})$
\begin{align*}
	&\bar\delta^2\omega(a_0,\dots,a_{k+2})=\sum_{t=0}^{k+2}(-1)^t\bar\delta\omega(a_0,\dots,a_{t-1},a_{t+1},\dots,a_{k+2})=\\
	&=\sum_{t=0}^{k+2}(-1)^t\sum_{s=0}^{t-1}(-1)^{s}\omega(a_0,\dots, a_{s-1},a_{s+1},\dots, a_{t-1},a_{t+1},\dots,a_{k+2})+\\
	&+\sum_{t=0}^{k+2}(-1)^t\sum_{s=t+1}^{k+2}(-1)^{s-1}\omega(a_0,\dots, a_{t-1},a_{t+1},\dots, a_{s-1},a_{s+1},\dots,a_{k+2})=0,
\end{align*}
because each term of the first sum appears in the second one with opposite sign. The same computation, with $\omega\in\mathcal{C}^{p,q,l,k}$ replaced by $c\in \check{C}^k(M,\mathcal{U})$, gives $\delta\circ\delta=0.$ The equality $d\circ d=0$ is a trivial consequence of Lemma \ref{lemma_ex_derivative} and, by definition, $(\overline{d}^2\omega)(a_0,\dots,a_k)=d^2(\omega(a_0,\dots,a_k))=0.$ Next, $(\bar\delta\circ i(\omega))(a_0,a_1)=i(\omega)(a_1)-i(\omega)(a_0)=\omega-\omega=0.$ Moreover, $\bar d\circ j(c)(a_0,\dots,a_k)=d(c(a_0,\dots,a_k))=0,$ because $c(a_0,\dots,a_k)$ is a constant smooth function on $M$. This proves the first line of equalities. Let us show commutation properties. 
\begin{equation*}
	\begin{split}
		&(\bar d\circ\bar\delta\omega)(a_0,\dots,a_{k+1})=d(\delta\omega(a_0,\dots,a_{k+1}))=\\=&\sum_{t=0}^{k+1}(-1)^td(\omega(a_0,\dots,a_{t-1},a_{t+1},\dots,a_{k+1}))=\bar\delta\circ\bar d\omega(a_0,\dots,a_{k+1}). 
	\end{split}
\end{equation*}
Similarly,
\begin{equation*}
	\begin{split}
		&(j\circ\delta(c))(a_0,\dots,a_{k+1})=\delta(c)(a_0,\dots,a_{k+1})=\sum_{t=0}^{k+1}(-1)^tc(a_0,\dots,a_{t-1},a_{t+1},\dots,a_{k+1})=\\
		&=\sum_{t=0}^{k+1}(-1)^tj(c)(a_0,\dots,a_{t-1},a_{t+1},\dots,a_{k+1})=\bar\delta\circ j(c)(a_0,\dots,a_{k+1}).
	\end{split}
\end{equation*}
Finally, 
\begin{equation*}
	\begin{split}
		(\bar d\circ i)(\omega)(a_0)=d(i(\omega)(a_0))=d\omega=(i\circ d)(\omega).
	\end{split}
\end{equation*}
This concludes the proof of Lemma \ref{lem_equalities}.

\qed

\begin{figure}
\begin{center}
	\begin{tikzcd}[column sep=1em]	
		\vdots&\vdots&\vdots&&\vdots&\\
		\mathcal{C}^{p,q,l}(M)\arrow{r}{i}\arrow[u,"d"]&\mathcal{C}^{p,q,l,0}(\mathcal{U})\arrow{r}{\bar\delta}\arrow[u,"\bar d"]&\mathcal{C}^{p,q,l,1}(\mathcal{U})\arrow[r,"\bar\delta"]\arrow[u,"d"]&\arrow[r,"\bar\delta"]\cdots  
		&\mathcal{C}^{p,q,l,k}(\mathcal{U})\arrow[u,"\bar d"]\arrow[r,"\bar\delta"]&\cdots\\
		\vdots\arrow[u,"d"]&\vdots\arrow[u,"\bar d"]&\vdots\arrow[u,"\bar d"]&&\vdots\arrow[u,"\bar d"]&\\
		\mathcal{C}^{p,q,1}(M)\arrow{r}{i}\arrow[u,"d"]
		&\mathcal{C}^{p,q,1,0}(\mathcal{U})\arrow{r}{\bar\delta}\arrow[u,"\bar d"]&\mathcal{C}^{p,q,1,1}(\mathcal{U})\arrow[r,"\bar \delta"]\arrow[u,"\bar d"]&\arrow[r,"\bar\delta"]\cdots&\mathcal{C}^{p,q,1,k}(\mathcal{U})\arrow[u,"\bar d"]\arrow[r,"\bar\delta"]&\cdots  \\
		\mathcal{C}^{p,q,0}(M) \arrow{u}{d}\arrow{r}{i}  
		&\mathcal{C}^{p,q,0,0}(\mathcal{U})\arrow{u}{\bar d}\arrow{r}{\bar\delta}&\mathcal{C}^{p,q,0,1}(\mathcal{U})\arrow{u}{\bar d}\arrow[r,"\bar\delta"]&\arrow[r,"\bar\delta"]\cdots&\mathcal{C}^{p,q,0,k}(\mathcal{U})\arrow{u}{\bar d}\arrow[r,"\bar\delta"]&\cdots\\
		&\check{C}^0(M,\mathcal{U})\arrow{u}{j}\arrow{r}{\delta}&\check{C}^1(M,\mathcal{U})\arrow{u}{j}\arrow[r,"\delta"]&\arrow[r,"\delta"]\cdots&\check{C}^k(M,\mathcal{U})\arrow{u}{j}\arrow[r,"\delta"]&\cdots 
	\end{tikzcd}
\end{center}
\caption{The commutative diagram that relates the cochain complex $(\mathcal{C}^{p,q,*}(M),d)$ to the cochain complex$(\check{C}^*(M,\mathcal{U}),\delta)$}
\label{fig_comm_diagram1}
\end{figure}
\noindent\textit{STEP 2: The zigzag argument.}\\
We can now apply a zigzag procedure to link closed non-exact elements of the first column to a  closed non-exact elements of the bottom row, and viceversa. We emphasize that at some point of the proof we need to go through the vertical arrows $\bar{d}$ in opposite direction, which morally represents an inversion of the exterior derivative $d$. The next Lemma \ref{lemma_Poincare} is the tool we need. In fact, we recall that, by classical results of differential geometry, every closed differential form is locally exact. This is an informal way of writing Poincaré's lemma. 
The following statement generalizes Poincaré's lemma to the context of our anisotropic Banach spaces.  In particular, this is analogous to \cite[Lemma 5.15]{butterley-kiamari-liverani} where the authors proved that every closed 1-currents is, in a way, locally exact. 

\begin{lemma}
	\label{lemma_Poincare}
	Let $h\in\mathcal{C}^{p,q,0}$ be a closed 0-current. Then, for any $k=1,\dots,m,$ there exists $c_k\in\C,$ such that $\phi_k\cdot(h-c_k)=0.$
	Let $\omega\in\mathcal{C}^{p,q,l}(M)$ be an $l$-current, for $l>0$. If $(d\omega)\phi_k=0$, then there exists $u_k\in\mathcal{C}^{p,q,l-1}(M)$ such that 
	$$d(u_k\phi_k)=\omega\phi_k+(-1)^{l-1}u_k\wedge d\phi_k.$$
\end{lemma}

\begin{remark}
	\label{rem_lemma_Poincare}
	Lemma \ref{lemma_Poincare} represents a key difference between this paper and \cite[Section 5]{butterley-kiamari-liverani} and it justifies the introduction of the spaces $\mathcal{C}^{p,q,l}$. Indeed, the authors of \cite{butterley-kiamari-liverani} only considered closed elements of $\mathcal{B}^{p,q,1}$ and such currents are locally exact with potential in $\mathcal{B}^{p+1,q-1,0}.$ In general, if $\omega\in\mathcal{B}^{p,q,l}$ is a closed current for $l>1$, then this is locally exact, but we cannot expect that the potential is more regular than $\omega$. Notice that the same issue holds for differential forms. In fact, if $\omega$ is a closed $C^r$ 1-form on a star-shaped domain $U$, then $\omega=du$ with $u\in C^{r+1}(U).$ On the other hand, if $\omega$ is a closed $C^r$ $l$-form on $U$ for some $l>1$, then $\omega=du$ for some $(l-1)$-form $u$ of class $C^r$ and coefficients of $u$ are $C^{r+1}$ exclusively in some directions. 
\end{remark}

We postpone the proof of Lemma \ref{lemma_Poincare} to Appendix \ref{technical_results} and we describe its application to the zigzag argument. In particular, we use the diagram in Figure \ref{fig_comm_diagram1} and Lemma \ref{lemma_Poincare} to define an isomorphism $$\Phi\colon \bar{H}_{dR}^{p,q,l}(M)\rightarrow H_{\check{C}}(M,\mathcal{U}).$$ We suggest looking at Figure \ref{fig_comm_diagram2} to follow the strategy of the zigzag procedure.
Let $[\omega]\in\bar{H}_{dR}^{p,q,l}(M)$ be an element of the anisotropic de Rham cohomology, where $\omega\in\mathcal{C}^{p,q,l}$ is closed. By definition, $i(\omega)\in\mathcal{C}^{p,q,l,0}(\mathcal{U})$, $(\bar d\circ i(\omega))=d\omega=0$ and $(\bar\delta\circ i(\omega))=0$, by Lemma \ref{lem_equalities}.  Therefore, by Lemma \ref{lemma_Poincare}, for any $a_0=1,\dots,m$, there exists $u_{a_0}^{(l-1)}\in\mathcal{C}^{p,q,l-1}$ such that $d(u_{a_0}^{(l-1)}\phi_{a_0})=\omega\phi_{a_0}+(-1)^{l-1}u_{a_0}^{(l-1)}\wedge d\phi_{a_0}.$ Let us define $u^{(l-1)}\in\mathcal{C}^{p,q,l-1,0}(\mathcal{U}),$ such that $u^{(l-1)}(a_0)=u^{(l-1)}_{a_0}.$ Consequently, $(\bar d u^{(l-1)})(a_0)\phi_{a_0}=(du_{a_0}^{(l-1)})\phi_{a_0}=\omega\phi_{a_0}.$  Next, we consider $\bar\delta u^{(l-1)}\in\mathcal{C}^{p,q,l-1,1}(\mathcal{U}).$
Since the diagram commutes,
\begin{equation*}
	\begin{split}
		&(\bar d\circ\bar\delta u^{(l-1)})(a_0,a_1)\phi_{a_0}\phi_{a_1}=(\bar\delta\circ\bar d u^{(l-1)})(a_0,a_1)\phi_{a_0}\phi_{a_1}=\\=&(du^{(l-1)}_{a_1})\phi_{a_0}\phi_{a_1}-(du^{(l-1)}_{a_0})\phi_{a_0}\phi_{a_1}=\omega\phi_{a_0}\phi_{a_1}-\omega\phi_{a_0}\phi_{a_1}=0.
	\end{split}
\end{equation*}
As a consequence, $d(\bar{\delta}u^{(l-1)}(a_0,a_1))\phi_{a_0}\phi_{a_1}=0,$ hence, by Lemma \ref{lemma_Poincare} there exists $u^{(l-2)}_{a_0,a_1}\in\mathcal{C}^{p,q,l-2}$ such that 
$(du^{(l-2)}_{a_0,a_1})\phi_{a_0}\phi_{a_1}=\bar\delta u^{(l-1)}(a_0,a_1)\phi_{a_0}\phi_{a_1}$ and $$d(u^{(l-2)}_{a_0,a_1}\phi_{a_0}\phi_{a_1})=\bar{\delta}u^{(l-1)}(a_0,a_1)\phi_{a_0}\phi_{a_1}+(-1)^{l-2}  u^{(l-2)}_{a_0,a_1}\wedge d(\phi_{a_0}\phi_{a_1}).$$
As above, we define $u^{(l-2)}\in\mathcal{C}^{p,q,l-2,2}(\mathcal{U})$ such that $u^{(l-2)}(a_0,a_1)=u^{(l-2)}_{a_0,a_1}.$ We point out that the choice of $u^{(l-2)}(a_0,a_1)$ is not unique, but it is unique up to elements $v\in\mathcal{C}^{p,q,l-2}$ such that $(dv)\phi_{a_0}\phi_{a_1}=0.$ Therefore, we set $u^{(l-2)}(a_0,a_1)=u^{(l-2)}_{a_0,a_1}$ for $a_0<a_1,$ while we impose that $u^{(l-2)}(a_1,a_0)=-u^{(l-2)}_{a_0,a_1}.$ In fact, this agrees with \eqref{eq_permutation} and
\begin{equation*}
	\begin{split}
		&du^{(l-2)}(a_0,a_1)\phi_{a_0}\phi_{a_1}=\bar\delta u^{(l-1)}(a_0,a_1)\phi_{a_0}\phi_{a_1}=(u^{(l-1)}(a_1)-u^{(l-1)}(a_0))\phi_{a_0}\phi_{a_1}=\\
		=&-(u^{(l-1)}(a_0)-u^{(l-1)}(a_1))\phi_{a_0}\phi_{a_1}=-\bar\delta u^{(l-1)}(a_1,a_0)\phi_{a_0}\phi_{a_1}=-du^{(l-2)}(a_1,a_0)\phi_{a_0}\phi_{a_1}.
	\end{split}
\end{equation*}
Repeating the argument, we consider $\bar\delta u^{(l-2)}.$ Then, for any $a_0,a_1,a_2\in\{1,\dots,m\},$
\begin{equation*}
	\begin{split}
		&d(\bar\delta u^{(l-2)}(a_0,a_1,a_2))\phi_{a_0}\phi_{a_1}\phi_{a_2}=(\bar d\circ\bar\delta u^{(l-2)})(a_0,a_1,a_2)\phi_{a_0}\phi_{a_1}\phi_{a_2}=\\=&(\bar\delta\circ\bar d u^{(l-2)})(a_0,a_1,a_2)\phi_{a_0}\phi_{a_1}\phi_{a_2}=(du^{l-2}_{a_1,a_2}-du^{l-2}_{a_0,a_2}+du^{l-2}_{a_0,a_1})\phi_{a_0}\phi_{a_1}\phi_{a_2}=\\=&(\bar\delta u^{(l-1)}(a_1,a_2)-\bar\delta u^{(l-1)}(a_0,a_2)+\bar\delta u^{(l-1)}(a_0,a_1))\phi_{a_0}\phi_{a_1}\phi_{a_2}=\\=&(\bar\delta\circ\bar\delta u^{(l-1)})(a_0,a_1,a_2)\phi_{a_0,a_1,a_2}=0,
	\end{split}
\end{equation*}
where, in the last line, we used $\bar\delta\circ\bar\delta=0.$ Consequently, there exists $u^{(l-3)}_{a_0,a_1,a_2}\in\mathcal{C}^{p,q,l-3}$ such that $(du^{(l-3)}_{a_0,a_1,a_2})\phi_{a_0}\phi_{a_1}\phi_{a_2}=\bar\delta u^{(l-2)}(a_0,a_1,a_2)\phi_{a_0}\phi_{a_1}\phi_{a_2}$ and $$d(u^{(l-3)}_{a_0,a_1,a_2}\phi_{a_0}\phi_{a_1}\phi_{a_2})=\bar{\delta}u^{(l-1)}(a_0,a_1)\phi_{a_0}\phi_{a_1}\phi_{a_2}+(-1)^{l-3} u^{(l-3)}_{a_0,a_1,a_2}\wedge d(\phi_{a_0}\phi_{a_1}\phi_{a_2}).$$
We fix a representative $u^{(l-3)}(a_0,a_1,a_2)=u^{(l-3)}_{a_0,a_1,a_2}$, for $a_0<a_1<a_2$, and we define $u^{(l-3)}(a_{\sigma(0)},a_{\sigma(1)},a_{\sigma(2)})=\text{sgn}(\sigma)u^{(l-3)}_{a_0,a_1,a_2},$ for any permutation $\sigma$ of $\{0,1,2\}.$
After $l$ steps, using the same procedure, we obtain $u^{(0)}\in\mathcal{C}^{p,q,0,l-1}(\mathcal{U})$ such that, for any $a_0,\dots,a_{l-1}\in\{1,\dots,m\},$
\begin{align*}
	d(u^{(0)}(a_0,\dots,a_{l-1})\phi_{a_0}\dots\phi_{a_{l-1}})=&\bar\delta u^{(1)}(a_0,\dots,a_{l-1})\phi_{a_0}\dots\phi_{a_{l-1}}+\\+&u^{(0)}(a_0,\dots,a_{l-1})\wedge d(\phi_{a_0}\dots\phi_{a_{l-1}}),\\
	d(u^{(0)}(a_0,\dots,a_{l-1}))\phi_{a_0}\dots\phi_{a_{l-1}}=&\bar\delta u^{(1)}(a_0,\dots,a_{l-1})\phi_{a_0}\dots\phi_{a_{l-1}}
\end{align*}
and 
$$(u^{(0)}(a_{\sigma(0)},\dots,a_{\sigma(l-1)}))=\text{sgn}(\sigma)(u^{(0)}(a_0,\dots,a_{l-1}))$$
for any permutation $\sigma$ of $\{0,\dots,l-1\}.$
Using again the commutation property of the diagram, we obtain that 
$$d(\bar\delta u^{(0)}(a_0,\dots,a_l))\phi_{a_0}\dots\phi_{a_l}=0,$$ for any $a_0,\dots,a_l\in\{1,\dots,m\}.$
Thus, by Lemma \ref{lemma_Poincare}, there exists $c_{a_0,\dots,a_l}\in\C$ such that
$$(\bar\delta u^{(0)}(a_0,\dots,a_l)-c_{a_0,\dots,a_l})\phi_{a_0}\dots\phi_{a_l}=0.$$
Consequently, we can choose $c\colon \check{C}_l(M,\mathcal{U})\rightarrow\C$ fixing $c(a_0,\dots,a_l)=c_{a_0,\dots,a_l}$, whenever $a_0<a_1<\dots<a_l$, and define $c(a_{\sigma(0)},\dots,a_{\sigma(l)})=\text{sgn}(\sigma)c_{a_{\sigma(0)},\dots,a_{\sigma(l)}},$ for each permutation $\sigma$ of $\{0,\dots,l\}.$ Thus, $c\in \check{C}^l(M,\mathcal{U})$ and $j(c)(a_0,\dots,a_l)\phi_{a_0}\dots\phi_{a_l}=\bar\delta u^{(0)}(a_0,\dots,a_l)\phi_{a_0}\dots\phi_{a_l}.$ Since 
\begin{equation*}
	\begin{split}
		&j\circ\delta(c)(a_0,\dots,a_{l+1})\phi_{a_0}\dots\phi_{a_{l+1}}=\bar\delta\circ j(c)(a_0,\dots,a_{l+1})\phi_{a_0}\dots\phi_{a_{l+1}}=\\=&\bar\delta\circ\bar\delta u^{(0)}(a_0,\dots,a_{l+1})\phi_{a_0}\dots\phi_{a_{l+1}}=0
	\end{split}
\end{equation*}
and $j$ is injective, $\delta(c)(a_0,\dots,a_{l+1})\phi_{a_0}\dots\phi_{a_{l+1}}=0,$ hence $\delta(c)(a_0,\dots,a_{l+1})=0.$ Accordingly, $c$ is a coboundary and we set $\Phi([\omega])=[c].$ 
\begin{figure}
\begin{footnotesize}
	\begin{center}
		\begin{tikzcd}[column sep=1em]	
			\vdots&\vdots&\vdots&&\vdots&\\
			\textcolor{blue}{\omega}\in\mathcal{C}^{p,q,l}(M)\arrow{r}{i}\arrow[u,"d"]&\textcolor{blue}{i(\omega)}\in\mathcal{C}^{p,q,l,0}(\mathcal{U})\arrow{r}{\bar\delta}\arrow[u,"\bar{d}"]&\mathcal{C}^{p,q,l,1}(\mathcal{U})\arrow[u,"\bar{d}"]\arrow[r,"\bar\delta"]&\cdots \arrow[r,"\bar{\delta}"]&\mathcal{C}^{p,q,l,l}\arrow[u,"\bar{d}"]&\\
			\mathcal{C}^{p,q,l-1}(M)\arrow{r}{i}\arrow[u,"d"]&\textcolor{blue}{u^{(l-1)}}\in\mathcal{C}^{p,q,l-1,0}(\mathcal{U})\arrow{r}{\bar\delta}\arrow[u,"\bar d"]&\textcolor{blue}{\bar\delta u^{(l-1)}}\in\mathcal{C}^{p,q,l-1,1}(\mathcal{U})\arrow[u,"\bar{d}"]\arrow[r,"\bar\delta"]&\cdots\arrow[r,"\bar\delta"]
			&\mathcal{C}^{p,q,l-1,l}(\mathcal{U})\arrow[u,"\bar{d}"]&\\ %%
			\mathcal{C}^{p,q,l-2}(M)\arrow{r}{i}\arrow[u,"d"]&\mathcal{C}^{p,q,l-2,0}(\mathcal{U})\arrow[r,"\bar\delta"]\arrow[u,"\bar d"]&\textcolor{blue}{u^{(l-2)}}\in\mathcal{C}^{p,q,l-2,1}(\mathcal{U})\arrow[u,"\bar d"]\arrow[r,"\bar\delta"]&\cdots\arrow[r,"\bar\delta"]
			&\mathcal{C}^{p,q,l-2,l}(\mathcal{U})\arrow[u,"\bar{d}"]&\\
			\vdots\arrow{u}{ d}&\vdots\arrow{u}{\bar d}&\vdots\arrow{u}{\bar d}&&\vdots\arrow{u}{\bar d}&\\
			\mathcal{C}^{p,q,0}(M)\arrow[u,"d"]\arrow{r}{i}&\mathcal{C}^{p,q,0,0}(\mathcal{U})\arrow{u}{\bar d}\arrow{r}{\bar\delta}&\mathcal{C}^{p,q,0,1}(\mathcal{U})\arrow{u}{\bar d}\arrow[r,"\bar\delta"]&\arrow[r,"\bar\delta"]\cdots&\textcolor{blue}{\bar\delta u^{(0)}}\in\mathcal{C}^{p,q,0,l}(\mathcal{U})\arrow{u}{\bar d}&\\
			&\check{C}^0(M,\mathcal{U})\arrow{u}{j}\arrow{r}{\delta}&\check{C}^1(M,\mathcal{U})\arrow{u}{j}\arrow[r,"\delta"]&\arrow[r,"\delta"]\cdots&\textcolor{blue}{c}\in\check{C}^l(M,\mathcal{U})\arrow{u}{j}&
		\end{tikzcd}
	\end{center}
\end{footnotesize}
\caption{The zigzag argument on the commutative diagram}
\label{fig_comm_diagram2}
\end{figure}

\noindent \textit{STEP 3: $\Phi$ is well-defined and invertible.}\\
We firstly prove that $\Phi$ is well-defined, that is, if $[\widetilde\omega]=[\omega]\in\bar{H}_{dR}^{p,q,l}(M),$ then it must be true that $[c_\omega]=\Phi([\omega])=\Phi([\widetilde\omega])=[c_{\widetilde\omega}].$ In fact, if $[\omega]=[\widetilde{\omega}]$, then there exists $q\in\mathcal{C}^{p,q,l-1}$ such that $dq=\omega-\widetilde\omega.$
Clearly, $i(\omega)$ and $i(\widetilde\omega)$ are cohomologous, i.e., $i(\omega)-i(\widetilde\omega)=i(\omega-\widetilde\omega)=i(dq)=\bar d(i(q)).$ Let $u^{(l-1)},\widetilde{u}^{(l-1)}\in\mathcal{C}^{p,q,l-1,0}(\mathcal{U})$ such that, for any $a_0\in\{1,\dots,m\},$ $du^{(l-1)}(a_0)\phi_{a_0}=i(\omega)\phi_{a_0}$ and $d\widetilde{u}^{(l-1)}(a_0)\phi_{a_0}=i(\widetilde\omega)\phi_{a_0},$ as determined above. It holds true that $$\bar d(\widetilde{u}^{(l-1)}-u^{(l-1)})(a_0)\phi_{a_0}=\bar d(i(\widetilde{\omega}-\omega))(a_0)\phi_{a_0}=\bar d\circ \bar d(i(q))(a_0)\phi_{a_0}=0,$$ hence, by Lemma \ref{lemma_Poincare}, there is $q^{(l-1)}\in\mathcal{C}^{p,q,l-2,0}(\mathcal{U})$ such that $$\bar dq^{(l-1)}(a_0)\phi_{a_0}=(\widetilde{u}^{(l-1)}-u^{(l-1)})(a_0)\phi_{a_0},$$ for each $a_0\in\{1,\dots,m\}.$ Iterating the previous procedure, we obtain that $$(\widetilde{u}^{(0)}-u^{(0)})(a_0,\dots,a_{l-1})\phi_{a_0}\dots\phi_{a_{l-1}}=j(q^{(0)})(a_0,\dots,a_{l-1})\phi_{a_0}\dots\phi_{a_{l-1}},$$ for some $q^{(0)}\in \check{C}^{l-1}(M,\mathcal{U}).$ Moreover,  $$\delta\widetilde{u}^{(0)}(a_0,\dots,a_l)\phi_{a_0}\dots\phi_{a_l}-j(c_{\widetilde{\omega}})(a_0,\dots,a_l)\phi_{a_0}\dots\phi_{a_l}=0$$ and $$\delta u^{(0)}(a_0,\dots,a_l)\phi_{a_0}\dots\phi_{a_l}-j(c_\omega)(a_0,\dots,a_l)\phi_{a_0}\dots\phi_{a_l}=0.$$
Accordingly, 
\begin{equation*}
	\begin{split}
		&(j(c_{\widetilde{\omega}})-j(c_\omega))(a_0,\dots,a_l)\phi_{a_0}\dots\phi_{a_l}=\bar\delta(\widetilde{u}^{(0)}-u^{(0)})(a_0,\dots,a_l)\phi_{a_0}\dots\phi_{a_l}=\\=&\bar\delta\circ j(q^{(0)})(a_0,\dots,a_l)\phi_{a_0}\dots\phi_{a_l}=j\circ\delta(q^{(0)})(a_0,\dots,a_l)\phi_{a_0}\dots\phi_{a_l},
	\end{split}
\end{equation*}
hence, by injectivity of $j$, $$(c_{\widetilde{\omega}}-c_\omega)(a_0,\dots,a_l)=\delta q^{(0)}(a_0,\dots,a_l),$$ i.e., $c_\omega$ and $c_{\widetilde{\omega}}$ are cohomologous.

We finally show that $\Phi$ is invertible. In particular, we construct the inverse map following the strategy of the direct map. Next result plays the role of Lemma \ref{lemma_Poincare} for the rows of the abelian diagram, i.e., it is a Poincaré's lemma for $\bar\delta.$
\begin{lemma}
	\label{lemma_Poincare2}
	Let $h\in\mathcal{C}^{p,q,l,0}(\mathcal{U})$, such that $\bar\delta h=0$. Then there exists $\omega\in\mathcal{C}^{p,q,l}$ such that $i(\omega)=h.$
	Let $\omega\in\mathcal{C}^{p,q,l,k}(\mathcal{U})$  for some $k>0$.If $\bar\delta\omega=0$, there exists $u\in\mathcal{C}^{p,q,l,k-1}(\mathcal{U})$ such that 
	$\bar\delta u=\omega.$
\end{lemma}
As above, we postpone the proof of Lemma \ref{lemma_Poincare2} to Appendix \ref{technical_results} and we apply it to conclude this proof.
Let $[c]\in H_{\check{C}}^l(M,\mathcal{U})$ be a \v{C}ech cohomology class. Since $\delta(c)=0$, it holds $\bar\delta\circ j(c)=j\circ\delta(c)=0.$ Consequently, Lemma \ref{lemma_Poincare2} implies that there exists $v^{(0)}\in\mathcal{C}^{p,q,0,l-1}(\mathcal{U})$ such that $\bar\delta v^{(0)}=j(c).$ Considering $\bar dv^{(0)},$ we get
$\bar\delta\circ\bar dv^{(0)}=\bar d\circ \bar\delta v^{(0)}=\bar d\circ j(c)=0.$ Thus, we can apply again Lemma \ref{lemma_Poincare2} to find a $v^{(1)}\in\mathcal{C}^{p,q,1,l-2}(\mathcal{U})$ such that $\bar\delta v^{(1)}=\bar dv^{(0)}.$ Iterating this argument, after $l$ steps, we obtain a  $v^{(l)}\in\mathcal{C}^{p,q,l-1,0}(\mathcal{U})$ such that $\bar\delta v^{(l)}=\bar d v^{(l-1)}.$ Since $\bar\delta\circ\bar dv^{(l)}=\bar d\circ \bar\delta v^{(l)}=\bar d\circ\bar dv^{(l-1)}=0,$ we can apply again Lemma \ref{lemma_Poincare2} and we conclude that there exists $\widetilde{\omega}\in\mathcal{C}^{p,q,l}$ such that $i(\widetilde{\omega})=\bar dv^{(l)}.$ One can easily check, using the same method described in the definition of $\Phi$, that $\widetilde{\omega}$ is a closed current and it is unique up to exact currents. As a consequence, we can define $\Psi\colon H_{\check{C}}^l(M,\mathcal{U})\rightarrow \bar{H}_{dR}^{p,q,l}(M)$ such that $\Psi([c])=[\widetilde{\omega}].$

It remains to show that $\Phi$ and $\Psi$ are inverse to each other. Let $[\omega]\in\bar{H}_{dR}^{p,q,l}(M),$ where $\omega\in\mathcal{C}^{p,q,l}.$ Then we have associated to $\omega$ a sequence $(u^{(s)})_{s=0,\dots,l-1}$ of elements $u^{(s)}\in\mathcal{C}^{p,q,s,l-1-s}(\mathcal{U})$ such that
$$d(u^{(s)}(a_0,\dots, a_{l-1-s}))\phi_{a_0}\dots\phi_{a_{l-1-s}}=\bar{\delta}u^{(s+1)}(a_0,\dots, a_{l-1-s})\phi_{a_0}\dots\phi_{a_{l-1-s}}.$$
On the other hand, the construction of $\Psi$ produces a sequence $(v^{(s)})_{s=0,\dots,l-1}$ such that $v^{(s)}\in\mathcal{C}^{p,q,s,l-1-s}$ such that $\bar{\delta}(v^{(s+1)})=\bar{d}v^{(s)},$ that is 
$$\bar{\delta}v^{(s+1)}(a_0,\dots, a_{l-1-s})\phi_{a_0}\dots\phi_{a_{l-1-s}}=d(v^{(s)}(a_0,\dots, a_{l-1-s}))\phi_{a_0}\dots\phi_{a_{l-1-s}}.$$
Without loss of generality we can assume $v^{(s)}=u^{(s)},$ so that $\Phi\circ\Psi=\Psi\circ\Phi=\id.$
This concludes the proof of Theorem \ref{iso-anisotropic}.

\qed

\begin{remark}
	\label{rem_finite_regularity}
	With some more effort and additional technicalities, it is possible to extend the main Theorem \ref{thm_strong} to $C^r$ Anosov diffeomorphisms and $C^r$ observables with $r\geq2$, and possibly to noninteger regularity. We highlight the main differences for the convenience of the reader. First, since the dynamics is only $C^r$, the image of $\Omega^l(M)$ by the pushforward operator $f_*$ is not anymore $C^\infty$. This requires to restrict to the space $\Omega^{l}_r(M)$ of $C^r$ differential forms. Moreover, in the definition of the anisotropic Banach spaces of currents we must assume $p+q\leq r$, limiting the possibility to decrease the essential spectral radius. More importantly, the differential operator $d\colon\mathcal{B}^{p,q,d_s}\rightarrow\mathcal{B}^{p-1,q+1,d_s+1}.$  
	Accordingly, if we want to define a cochain complex $(\mathcal{B}^*,d),$ we must assume $q-d_s\geq 0$ and $p-d_u\geq0$, which forces $r\geq\max\{d_s,d_u\}$. This shows that $\Omega^{l}_r(M)$ is not again the right space to work with. In fact, one should consider the space $\widetilde{\Omega}^l_r(M)$ of $l$-differential forms $\omega\in\Omega^l_r(M)$ such that $d\omega\in\Omega^{l+1}_r(M).$ The idea is very similar to the definition of anisotropic Banach spaces $\mathcal{C}^{p,q,l},$ but it requires to prove the Lasota-Yorke inequalities of Theorem \ref{thm_lasota-yorke} for elements in $\widetilde{\Omega}^l_r(M)$ with respect to the norm $|\omega|_{p,q,l}=\norm{\omega}_{p,q,l}+\norm{\omega}_{p,q,l+1}.$ This is the main technical step that can be avoided by considering $C^\infty$ Anosov diffeomorphisms. The rest of the proof is unchanged.

\end{remark}

\section{Proofs of Theorem \ref{thm_strong} and Corollary \ref{cor_main}}
\label{sec_proofs}
We are ready to prove our main theorem and its corollary.

\proofof{Theorem \ref{thm_strong}} As a consequence of Theorem \ref{iso-anisotropic} and Corollary \ref{cor_spectrum2}
\begin{equation*}
	\begin{split}
		\sigma(f_*|_{\mathcal{C}^{p,q,d_s}})\cap\{z\in\mathbb{C}: |z|>\lambda^{-1}e^{h_\top }\}&=\sigma(f_{\#}|_{\bar{H}_{dR}^{p,q,d_s}})\cap\{z\in\mathbb{C}: |z|>\lambda^{-1}e^{h_\top }\}\\&=\sigma(f_{\#}|_{H_{dR}^{d_s}})\cap\{z\in\mathbb{C}: |z|>\lambda^{-1}e^{h_\top }\}.
	\end{split}
\end{equation*}
Since $e^{h_\top }$ is a simple maximal eigenvalue of $f_*|_{\mathcal{B}^{p,q,d_s}},$ it holds that $\Lambda_1,$ the maximal eigenvalue of $f_\#|_{H_{dR}^{d_s}(M)},$ is $e^{h_\top }$, while the second one fulfills $|\Lambda_2|<e^{h_\top }.$ In addition, by Proposition \ref{prop_measure_of_maximal_entropy}, the eigenvector $\hat\omega$ and the dual eigenvector $\hat t$ related to $e^{h_\top }$ defines the measure of maximal entropy $\mu_{BM}.$ Furthermore, for any other eigenvalue $\Lambda_i$, with $|\Lambda_i|>\lambda^{-1}e^{h_\top },$ we set a Jordan basis $\{\hat\omega_{i,k}\}_{k=1}^{N_i},$  such that $f_*(\hat\omega_{i,1})=\Lambda_i\hat\omega_{i,1}$ and $f_*(\hat\omega_{i,k})=\Lambda_i\hat\omega_{i,k}+\hat\omega_{i,k-1},$ for $k=2,\dots, N_i.$
Let $\{\hat t_{i,k}\}_{k=1}^{N_i}$ be the dual Jordan basis, such that  $\hat t_{i,k}(\hat\omega_{i,j})=\delta_{k,j}.$ Notice that $f_*'\hat t_{i,N_i}=\Lambda_i\hat t_{i,N_i}$ and $f_*'\hat t_{i,k}=\Lambda_i\hat t_{i,k}+\hat t_{i,k+1},$ for $k=1,\dots,N_i-1$.
We point out that every $\Lambda_i$ corresponds to a single Jordan block, because eigenvalues are counted according to their geometric multiplicity. In particular, some of the eigenvalues $\Lambda_2,\dots, \Lambda_m$ such that $ |\Lambda_i|>\lambda^{-1}e^{h_\top}$, may coincide. We obtain that
\begin{equation*}
	\begin{split}
		f_*=e^{h_\top }\hat\omega\otimes\hat t+\sum_{i=2}^m\left[\Lambda_i\left(\sum_{j=1}^{N_i}\hat\omega_{i,j}\otimes \hat t_{i,j}\right)+\sum_{j=1}^{N_i-1}\hat\omega_{i,j}\otimes \hat t_{i,j+1}\right]+\mathcal{Q},
	\end{split}
\end{equation*}
where $\mathcal{Q}$ is a linear operator such that $\norm{\mathcal{Q}}_{(\mathcal{C}^{p,q,d_s})'}\leq\lambda^{-1}e^{h_\top }$. As a consequence,
\begin{equation*}
	\begin{split}
		f_*^n=e^{h_\top }\hat\omega\otimes\hat t+\sum_{i=2}^m\sum_{k=0}^{N_i-1}\binom{n}{k}\Lambda_i^{n-k}\left(\sum_{j=1}^{N_i-k}\hat\omega_{i,j}\otimes \hat t_{i,j+k}\right)+\mathcal{Q}^n,
	\end{split}
\end{equation*}
where $\binom{n}{k}=0$ for $k>n$. Then, for any $\phi,\psi\in C^\infty(M),$
\begin{equation*}
	\begin{split}
		&\int_M\!\!\!\phi(\psi\circ f^n) d\mu_{BM}=\mu_{BM}(\phi(\psi\circ f^n))=\hat t(\phi(\psi\circ f^n)\hat\omega)=e^{-nh_\top }\hat t(f_*^n(\phi\psi\circ f^n\hat\omega))=\\=&e^{-nh_\top }\hat t(\psi(f_*^n(\phi \bar \omega))=\hat t(\phi \hat\omega)\hat t(\psi \hat\omega)+\\&+\sum_{i=2}^m\sum_{k=0}^{N_i-1}\binom{n}{k}\Lambda_i^{n-k}e^{-nh_\top }\left(\sum_{j=1}^{N_i-k} \hat t_{i,j+k}(\phi\hat\omega)\hat t(\psi\hat\omega_{i,j})\right)+
		e^{-nh_\top }\hat t(\psi \mathcal{Q}^n(\phi\hat\omega)).
	\end{split}
\end{equation*}
By defining a finite number of bilinear forms $\{c_{\Lambda_i,k}(\cdot,\cdot)\}_{\substack{i=2,\dots,m\\k=0,\dots,N_i-1}}$ such that 
$$\sum_{k=0}^{N_i-1}\binom{n}{k}\Lambda_i^{-k}\left(\sum_{j=1}^{N_i-k} \hat t_{i,j+k}(\phi\hat\omega)\hat t(\psi\hat\omega_{i,j})\right)=\sum_{k=0}^{N_i-1}n^kc_{\Lambda_i,k}(\phi,\psi),$$
we get
\begin{equation*}
	\begin{split}
		\int_M\phi(\psi\circ f^n) d\mu_{BM}=&\mu_{BM}(\phi)\mu_{BM}(\psi)+\\&+\sum_{i=2}^m\sum_{k=0}^{N_i-1}(\Lambda_ie^{-h_\top })^n n^kc_{\Lambda_i,k}(\phi,\psi)+
		e^{-nh_\top }\hat t(\psi \mathcal{Q}^n(\phi\hat\omega)),
	\end{split}
\end{equation*}

hence
\begin{equation*}
	\begin{split}
		\left|\int_{M}\phi(\psi\circ f^n) d\mu_{BM}-\right.&\int_{M}\phi d\mu_{BM}\int_{M}\psi d\mu_{BM}-\\-&\left.\sum_{i=2}^{m}\sum_{k=0}^{N_i-1}(\Lambda_ie^{-h_\top })^nn^k c_{\Lambda_i,k}(\phi,\psi)\right|\leq C\lambda^{-n} \norm{\phi}_{C^r}\norm{\psi}_{C^r}.
	\end{split}
\end{equation*}
\qed

\proofof{Corollary \ref{cor_main}}
By classical results of Franks \cite{franks,franks_anosov_diffeonorphisms_of_tori} and Newhouse \cite{newhouse},  $f$ is topologically conjugate to a hyperbolic automorphism of the torus $F\colon\mathbb{T}^{\dim(M)}\rightarrow\mathbb{T}^{\dim(M)}.$ Accordingly, since $F$ is topologically transitive \cite{micena}, the same holds for $f$ and we can apply Theorem \ref{thm_strong}. It is enough to show that the second highest eigenvalue $\Lambda_2$ of $F_\#|_{H_{dR}^{d_s}(M)}$ satisfies $|\Lambda_2|<\lambda^{-1}e^{h_\top }$. Since the action of the dynamics on de Rham cohomology is invariant under topological conjugacy, it holds $f_{\#}=F_{\#}.$ But $F$ is linear, hence, with a slight abuse of notation we can write $F^{-1}=F_{\#}|_{H_{dR}^1(M)}.$ Notice that this equality is not properly correct, because $F^{-1}$ acts on $\mathbb{T}^{\dim(M)},$ while $F_{\#}$ acts on $H_{dR}^1(\mathbb{T}^{\dim(M)}).$ On the other hand, fixing a basis of $\R^{\dim(M)}$, $F^{-1}$ is induced by a $\dim(M)\times\dim(M)$ matrix $A\in GL_{\dim(M)}(\Z)$ with $\det(A)=\pm1,$ and the same $A$ is the matrix associated to $F_{\#}$ w.r.t.\ the canonical basis of $H_{dR}^1(\mathbb{T}^{\dim(M)}),$ that is $\{[dx_1],[dx_2],\dots,[dx_{\dim(M)}]\}.$
Next, assume that $\sigma_1=\{\nu_1,\nu_2,\dots,\nu_{\dim(M)}\}$ is the spectrum of $f_\#$ acting on $H_{dR}^1(\mathbb{T}^{\dim(M)}).$ Then, we can assume  that
$$|\nu_1|\leq|\nu_2|\leq\dots\leq|\nu_{d_u}|<1<|\nu_{d_u+1}|\leq|\nu_{d_u+2}|\leq\dots\leq|\nu_{\dim(M)}|$$
and it holds true that $\left|\prod_{i=1}^{\dim(M)}\nu_i\right|=1.$ Moreover, the spectrum $\sigma_l$ of $F_{\#}$ acting on $H_{dR}^l(\mathbb{T}^{\dim(M)})$ can be determined by multiplying $l$ eigenvalues of $\sigma_1,$ that is, $\sigma_l=\{\prod_{i\in I}\nu_i\ |\ I \subseteq\{1,\dots,\dim(M)\},\ |I|=l\}.$ 
Since $e^{h_\top }$ is the maximal eigenvalue of $F_\#|_{H_{dR}^{d_s}(M)}$, we have
$$e^{h_\top }=\prod_{i=d_u+1}^{\dim(M)}\nu_i=\prod_{i=1}^{d_s}\nu_i^{-1}.$$  
Notice that this equality agrees with Ledrappier-Young entropy formula 
\cite[Theorem D.3.1]{brown}
$$h_\top =\sum_{i=d_s}^{\dim(M)}\ln\nu_i=\left(\sum_{i=1}^{d_u}\ln\nu_i\right)^{-1}.$$
In addition, we deduce that 
$$\Lambda_2=\nu_{d_u}\prod_{i=d_u+2}^{\dim(M)}\nu_i,$$
while the maximal eigenvalues of $F_\#$ acting on $H_{dR}^{d_s-1}(M)$, resp.\ $H_{dR}^{d_s+1}(M),$ is 
$$\tau_{d_s-1}=\prod_{i=d_u}^{\dim(M)}\nu_i,\text{ resp.\ }\tau_{d_s+1}=\prod_{i=d_u+2}^{\dim(M)}\nu_i.$$
Furthermore, $|\tau_{d_s-1}|\leq\lambda^{-1}e^{h_\top }$ and $|\tau_{d_{s+1}}|\leq\lambda^{-1}e^{h_\top }.$ In fact, by \ref{cor_spectrum2}, $\tau_{d_s-1},$ resp.\ $\tau_{d_s+1},$ is an eigenvalue of $f_*|_{\mathcal{C}^{p,q,d_s-1}},$ resp.\  $f_*|_{\mathcal{C}^{p,q,d_s+1}},$ and $\rho(f_*|_{\mathcal{C}^{p,q,d_s\pm 1}})\leq\lambda^{-1}e^{h_\top }$.
Finally, noticing that 
$$|\Lambda_2|<\min\{|\tau_{d_s-1}|,|\tau_{d_s+1}|\}\leq\lambda^{-1}e^{h_\top },$$ we conclude that 
$$|\Lambda_2|<\lambda^{-1}e^{h_\top }.$$

\qed

\begin{appendices}
	\section{Technical results}
	\label{technical_results}
	In this section we collect some technical proofs and results that we used in this paper.
	
	The following distortion lemma is a classical result in the field of dynamical systems. On the other hand, we are not aware of any proof adapted to our setting. Thus, we also give the proof.   
	\begin{lemma}[Distortion Lemma]
		\label{lemma_determinant}
		Let $W\in\Sigma$ be an admissible leaf. Let $W_i$ be an admissible leaf, as in Lemma \ref{lemma_foglie}, such that $W_i\subseteq f^{-n}(W).$ Define $\lambda_n^s(x)=\left|\det (d_xf^n|_{T_xW_i})\right|$, that is the contraction of $f^n$ along the stable manifold $W_i.$ Let $\lambda_{n,i}^s=\min_{x\in W_i}\lambda_n^s(x).$ Then there exists a constant $C>0$ such that, for any $n\in\N$ and for each $x\in W_i,$
		$$\lambda_{n,i}^s\leq\lambda_n^s(x)\leq C\lambda_{n,i}^s.$$
		In addition, $\norm{\lambda_n^s}_{C^q(W_i)}\leq |f^n(W_i)|.$ 
	\end{lemma}
	\proof
	We claim that there exists a constant $D>0$ such that for any $n\in\N$ and for any $x,y\in W_i$
	$$D^{-1}\leq\frac{\lambda_n^s(x)}{\lambda_n^s(y)}\leq D.$$
	In fact, assuming that the claim is true, we get
	$\lambda_n^s(x)\leq D\lambda_n^s(y),$ for all $x,y\in W_i,$ hence $\lambda_n^s(x)\leq D\min_{y\in W_i}\lambda_n^s(y)=D\lambda_{n,i}^s.$
	
	Let us prove the claim. Denoting by $W_i^t=f^t(W_i),$ we obtain 
	\begin{equation*}
		\lambda_n^s(x)=\left|\det(d_xf^n|_{T_xW_i})\right|=\left|\det\left(\prod_{t=0}^{n-1}d_{f^t(x)}f|_{T_{f^t(x)}W_i^t}\right)\right|=\prod_{t=0}^{n-1}\left|\det\left(d_{f^t(x)}f|_{T_{f^t(x)}W_i^t}\right)\right|.
	\end{equation*}
	Consequently,
	\begin{equation*}
		\begin{split}
			&\ln\frac{\lambda_n^s(x)}{\lambda_n^s(y)}=\sum_{t=0}^{n-1}\ln\frac{\left|\det (d_{f^t(x)}f|_{T_{f^t(x)}W_i^t})\right|}{\det\left| d_{f^t(y)}f|_{T_{f^t(y)}W_i^t}\right|}=\sum_{t=0}^{n-1}\langle\nabla\ln\left|\det(d_{f^t(z)}f|_{T_{f^t(z)}W_i^t})\right|, x-y\rangle=\\
			=&\sum_{t=0}^{n-1}\langle(d_zf^t)^T\nabla\ln\left|\det(d_{f^t(z)}f|_{T_{f^t(z)}W_i^t})\right|,x-y\rangle=\\=&\sum_{t=0}^{n-1}\langle\nabla\ln\left|\det(d_{f^t(z)}f|_{T_{f^t(z)}W_i^t})\right|,d_zf^t(x-y)\rangle\leq\\\leq& \sum_{t=0}^{n-1}\underbrace{\max_{l\in W_i^t}\norm{\nabla\ln\left|\det(d_{l}f|_{T_{l}W_i^t})\right|}\max_{x,y\in W_i}\norm{x-y}}_{\leq C}\underbrace{\max_{z\in W_i}\norm{d_zf^t|_{T_xW_i}}}_{C\lambda^{-t}}\leq\\\leq&\sum_{t=0}^{n-1} C\lambda^{-t}\leq C\frac{\lambda^{-1}}{1-\lambda^{-1}}\leq D.
		\end{split}
	\end{equation*}
	We have proved $\norm{C^0(W_i)}\leq C\lambda_{n,i}^s.$ Notice that $$|f^n(W_i)|=\int_{f^n(W_i)}\omega_W=\int_{W_i}\lambda_n^s(x)\omega_{W_i}(x)\geq\lambda_{n,i}^s|W_i|\geq C\lambda_{n,i}^s,$$
	hence $\norm{\lambda_n^s}_{C^0(W_i)}\leq C|f^n(W_i)|.$ 
	Let us consider derivative of $\lambda_n^s.$ We compute 
	\begin{equation*}
		\begin{split}
			&\de_{x_j}\lambda_n^s(x)=\de_{x_j}\left(\prod_{t=0}^{n-1}\left|\det (d_{f^t(x)}f|_{T_{f^t}(x)W_i^t})\right|\right)=\\=&\sum_{h=0}^{n-1}\prod_{t=0}^{n-1}\left|\det (d_{f^t(x)}f|_{T_{f^t}(x)W_i^t})\right|\frac{\de_{x_j}\left|\det (d_{f^h(x)}f|_{T_{f^h}(x)W_i^h})\right|}{\left|\det (d_{f^h(x)}f|_{T_{f^h}(x)W_i^h})\right|}=\\=&\lambda_n^s(x)\sum_{h=0}^{n-1}\de_{x_j}\ln\left(\left|\det (d_{f^h(x)}f|_{T_{f^h(x)}W_i^h})\right|\right)=\\=&\lambda_n^s(x)\sum_{h=0}^{n-1}\left\langle\nabla\ln\left(\left|\det (d_{f^h(x)}f|_{T_{f^h(x)}W_i^h})\right|\right),\de_{x_j}f^h(x)|_{T_xW_i}\right\rangle\leq\\\leq&\lambda_n^s(x)\sum_{h=0}^{n-1}\max_{y\in W_i^h}\norm{\nabla\ln|\det (d_y f|_{T_yW_i^h})|}\norm{d_xf^h|_{T_xW_i}}\leq C\lambda_n^s(x)\sum_{h=0}^{n-1}\lambda^{-h}\leq C\lambda_n^s(x).
		\end{split}
	\end{equation*}
	Accordingly, $\norm{\lambda_n^s}_{C^1(W_i)}\leq C\norm{\lambda_n^s}_{C^0(W_i)}\leq C|f^n(W_i)|$ and, iterating the procedure, $\norm{\lambda_n^s}_{C^q(W_i)}\leq C|f^n(W_i)|.$
	
	\qed
	What follows is a generalization of the distortion Lemma \ref{lemma_determinant} to differential forms.
	\begin{lemma}
		\label{lemma_determinant2}
		Let $W\in\Sigma$ be an admissible leaf in $\psi_z(U_z)$. Let $W_i$ be an admissible leaf, as in Lemma \ref{lemma_foglie}, such that $W_i\subseteq f^{-n}(W)$ in $\psi_i(U_i).$ Define $\lambda_n^s(x)=\left|\det (d_xf^n|_{T_xW_i})\right|$, that is the contraction of $f^n$ along the stable manifold $W_i.$ 
		Let $dx_{\overline{k}}\in\Omega^l(\psi_i(U_i))$ and $dx_{\overline{j}}\in\Omega^l(\psi_z(U_z))$, as defined  in Section \ref{Section_Norms_and_Banach_spaces}. Then, 
		$$\norm{\langle dx_{\overline{j}},f_*^ndx_{\overline{k}}\rangle\circ f^n\lambda_n^s}_{C^q(W_i)}\leq C\lambda^{-|d_s-l|n}.$$ 
	\end{lemma}
	\proof
	We firstly consider the case $l\leq d_s.$ Let $\mathcal{S}_{x}^l$ be the family of $l$-dimensional vector subspaces of $\mathcal{C}_x^s.$ For any $V_x\in\mathcal{S}_x^l$ and for any $n\in\N$, let us denote by $V_{f^n(x)}=d_xf^n V_x.$
	Then 
	\begin{equation*}
		\begin{split}
			&|\langle dx_{\overline{j}},f_*^ndx_{\overline{k}}\rangle\circ f^n(x)\lambda_n^s(x)|=	|\langle dx_{\overline{j}},f_*^ndx_{\overline{k}}\rangle\circ f^n(x)\det(d_xf^n|_{T_xW_i})|\leq\\
			\leq&\max_{V_x\in\mathcal{S}_x^l} |\det(d_{f^n(x)}f^{-n}|_{V_{f^n(x)}})\det(d_xf^n|_{T_xW_i})|=\max_{V_x\in\mathcal{S}_x^l}\left|\frac{\det(d_{f^n(x)}f^{-n}|_{V_{f^n(x)}})}{\det(d_{f^n(x)}f^{-n}|_{T_{f^n(x)W}})}\right|
		\end{split}
	\end{equation*}
	Notice that 
	$$\max_{V_x\in\mathcal{S}_x^l}|\det(d_{f^n(x)}f^{-n}|_{V_{f^n(x)}})|\leq \max_{V_x\in\mathcal{S}_x^{d_s}}|\det(d_{f^n(x)}f^{-n}|_{V_{f^n(x)}})|\lambda^{-(d_s-l)n},$$
	hence 
	$$\max_{V_x\in\mathcal{S}_x^l}\left|\frac{\det(d_{f^n(x)}f^{-n}|_{V_{f^n(x)}})}{\det(d_{f^n(x)}f^{-n}|_{T_{f^n(x))W}})}\right|\leq \max_{V_x\in\mathcal{S}_x^{d_s}}\left|\frac{\det(d_{f^n(x)}f^{-n}|_{V_{f^n(x)}})}{\det(d_{f^n(x)}f^{-n}|_{T_{f^n(x)W}})}\right|\lambda^{-(d_s-l)n}$$
	Since both $d_{f^n(x)}f^{-n}(V_{f^n(x)})$ and $d_{f^n(x)}f^{-n}(T_{f^n(x)W})$ converge to the stable subbundle as $n\rightarrow +\infty,$ by continuity of the differential, there exist $\nu\in(0,1)$ and $\bar n\in\N$ such that, for $m>\bar n,$
	$$|\det(d_{f^m(x)}f^{-1}|_{V_{f^m(x)}})|-|\det(d_{f^m(x)}f^{-1}|_{T_{f^{m}(x)}f^m(W_i)})|<\nu<1,$$
	hence
	$$|\det(d_{f^n(x)}f^{-n}|_{V_{f^n(x)}})|-|\det(d_{f^n(x)}f^{-n}|_{T_{f^n(x)}W})|< C \nu^{n-\bar n}\leq C$$
	We obtain that 
	$$\max_{V_x\in\mathcal{S}_x^{d_s}}\left|\frac{\det(d_{f^n(x)}f^{-n}|_{V_{f^n(x)}})}{\det(d_{f^n(x)}f^{-n}|_{T_{f^n(x)W}})}\right|\leq 1+\frac{1}{|\det(d_{f^n(x)}f^{-n}|_{T_{f^n(x)W}})|}\leq 1+\frac{C}{\lambda^{nd_s}}\leq C$$
	We conclude that 
	$$\norm{\langle dx_{\overline{j}},f_*^ndx_{\overline{k}}\rangle\circ f^n(x)\lambda_n^s(x)}_{C^0(W_i)}\leq C\lambda^{-(d_s-l)n}$$
	Let us compute derivatives. Without loss of generality, we can assume that the derivative $\de_{x_j}\in T_xW_i.$ Then
	\begin{equation*}
		\begin{split}
			\de_{x_j}(|\det(d_{f^n(x)}f^{-n}|_{V_{f^n(x)}})\lambda_n^s(x)|)=&\de_{x_j}(|\det(d_{f^n(x)}f^{-n}|_{V_{f^n(x)}})\|)\lambda_n^s(x)+\\
			&+|\det(d_{f^n(x)}f^{-n}|_{V_{f^n(x)}})\|\de_{x_j}\lambda_n^s(x)
		\end{split}
	\end{equation*}
	Since $\de_{x_j}\lambda_n^s(x)$ (see the proof of Lemma \ref{lemma_determinant}), the second term of the sum can be estimated by the $C^0$-norm. 
	Similarly, one can repeat the argument of the proof of Lemma \ref{lemma_determinant} to prove that $$\de_{x_j}|\det(d_{f^n(x)}f^{-n}|_{V_{f^n(x)}})|\leq C|\det(d_{f^n(x)f^{-n}|_{V_{f^{n}(x)}}})|$$
	Thus, also the first term of the sum can be estimated by the $C^0$-norm. By iterating that procedure one obtains the estimate for the $C^q$-norm. 
	It remains to prove the case $l>d_s.$ Let $V_x$ be an $l$-dimensional subspace of $T_xM.$ The worst estimate for $|\det(d_{f^n(x)}f^{-n}|_{V_{f^n(x)}})|$ is given by the case for which $V_x=S_x\oplus U_x,$ where $S_x$ is a $d_s$-dimensional subspace in the stable cone and $U_x$ is a $(l-d_s)$-dimensional unstable subspace. Since $|\det (d_{f^n(x)}f^{-n}|_{U_{f^n(x)}})|\leq\lambda^{-(l-d_s)n},$ we obtain again 
	\begin{equation*}
		|\langle dx_{\overline{j}},f_*^n dx_{\overline{k}}\rangle\circ f^n(x)\lambda_n^s(x)|\leq \max_{S_x\in\mathcal{S}_x^{d_s}}\left|\frac{\det(d_{f^n(x)}f^{-n}|_{S_{f^n(x)}})}{\det(d_{f^n(x)}f^{-n}|_{T_{f^n(x)}W})}\right|\lambda^{-(l-d_s)n}
	\end{equation*}
	and we can conclude as above.
	\qed
	
	We now prove two generalizations of Poincaré's lemma that we used in Section \ref{section_connection_de_Rham} to get the isomorphism between the anisotropic de Rham cohomology and the \v{C}ech cohomology.
	
	\proofof{Lemma \ref{lemma_Poincare}}
	The first part of the lemma informally says that closed $0$-currents are constant.
	Let $h\in\mathcal{C}^{p,q,0}$ be a closed current of degree 0 (and dimension $\dim(M)$). By Lemma \ref{lemma_correnti}, $\mathcal{B}^{p,q,0}$ can be identified with a subspace the dual space of $C^{p+q}$ function on $M,$ i.e., a space of distribution, and the behavior of $h$, acting on $\phi\in C^{p+q}$ as a current, is obtained disintegrating the integral 
	$$i(h)(\phi)=\int_M  h\phi\omega_0$$
	along leaves of $\Sigma.$ The same holds true for $\mathcal{C}^{p,q,0},$ which is a subspace of $\mathcal{B}^{p,q,0}.$ Let $\phi\in C^{p+q}(M)$ such that $\text{supp}(\phi)\subseteq \text{int}(\text{supp}(\phi_k)),$ where we recall that $\{\phi_k\}_{k=1}^m$ is the partition of unity subordinated to the contractible open covering $\{V_k=\psi_k(U_k)\}$. Then
	\begin{equation*}
		\begin{split}
			i(\de_{x_i}h)(\phi)=\int_M\phi\de_{x_i}h\omega_0=\!\int_{a_1}^{b_1}\!\dots\!\int_{a_{d_u}}^{b_{d_u}}\!\int_{W_{t_1,\dots t_{d_u}}}\!\phi\de_{x_i}h f_{t_1,\dots,t_{d_u}}\omega_{t_1,\dots,t_{d_u}}dt_1\dots dt_{d_u}
		\end{split}
	\end{equation*}
	where we disintegrated the integral along $d_s$-dimensional admissible leaves $W_{t_1,\dots t_{d_u}}\in\Sigma$ depending on $d_u$ parameters $t_1,\dots t_{d_u}$. Notice that, since leaves are smooth, the Jacobian $f_{t_1,\dots,t_{d_u}}$ is a smooth function. Consequently, 
	\begin{equation*}
		\begin{split}
			i(\de_{x_i}h)(\phi)=\!\int_{a_1}^{b_1}\!\dots\!\int_{a_{d_u}}^{b_{d_u}}\!\int_{W_{t_1,\dots t_{d_u}}}\!\langle f_{t_1,\dots,t_{d_u}}\phi_k^{-1}\phi dx_i,\phi_kdh\rangle\omega_{t_1,\dots,t_{d_u}}dt_1\dots dt_{d_u}=0.
		\end{split}
	\end{equation*}
	This implies that $\de_{x_i}h=0$ as a distribution on the interior of $\text{supp}(\phi_k).$ Since, by assumption, $\text{supp}(\phi_k)$ is simply connected, hence connected, there exists $c_k\in\C$ such that $\phi_k\cdot(h-c_k)=0$ in $(C^{p+q})^*.$ Thus, $\phi_k\cdot(h-c_k)=0$ in $\mathcal{C}^{p,q,0},$ because the inclusion $i'\colon\mathcal{C}^{p,q,0}\rightarrow\mathcal{B}^{p,q,0}$ and the map $\iota\colon\mathcal{B}^{p,q,0}\rightarrow(C^{p+q})^*$ are injective (see Lemma \ref{lemma_correnti}).
	
	Let us prove the second part of the lemma. We firstly introduce an enlarged partition of unity $\{\bar{\phi}_k\}_{k=1}^m$ such that $\text{supp}(\phi_k)\subseteq \text{supp}(\phi_k^+)\subseteq V_k$ and $\phi_k^+=1$ on $\text{supp}(\phi_k).$ Next
	consider a differential form $\omega\in\Omega^{l}(M)$ such that $d\omega\phi_k=0$ for some $k\in\{1,\dots,m\}.$ We fix $x_k\in V_k\setminus \text{supp}(\phi_k^+)$ and we define the linear operator $\alpha_k\colon\Omega^l(M)\rightarrow\Omega^{l-1}(M)$ such that
	$$\alpha_k(\omega)_x=\psi_{k*}\left(\int_0^1t^{l-1}(\iota_{\psi_k^{-1}(x)-\psi_k^{-1}(x_k)}\psi_k^*\omega)_{\psi^{-1}(x_k)(1-t)+\psi^{-1}_k(x)t}dt\right),$$
	where $\iota_v\omega$  denotes the interior product of $\omega$ with a vector field $v.$ Next, denoting by $y_k=(r_1^{(k)},\dots,r_{\dim(M)}^{(k)})=\psi_k^{-1}(x_k)$ and $y=(r_1,\dots,r_{\dim(M)})=\psi_k^{-1}(x),$ we can compute
	\begin{equation}
		\label{eq_ap1}
		\begin{split}
			&((d\alpha_k+\alpha_kd)\omega)_x=\psi_{k*}\left(\int_{0}^1t^{l-1}((d\iota_{y-y_k}+\iota_{y-y_k}d)\psi_k^*\omega)_{y_k(1-t)+yt}dt\right)=\\
			=&\psi_{k*}\left(\int_{0}^1t^{l-1}(L_{y-y_k}\psi_k^*\omega)_{y_k(1-t)+yt}dt\right)
		\end{split}
	\end{equation}
	where, in the second line, we used  Cartan's magic formula \cite[Theorem 2.11]{morita} which states that, given a smooth vector field $X$,  $L_X=d\iota_X+\iota_Xd.$ Let $\Phi_s(y)=y_k+e^s(y-y_k)$, so that $\Phi_0(y)=y$ and $\frac{d}{ds}|_{s=0}\Phi_s(y)=y-y_k.$ Writing $v_k$ in place of $\psi_k^*\omega$ and assuming that $(v_k)_z=\sum_{\overline{j}\in\mathcal{J}^l}f_{\overline{j}}(z)dr_{\overline{j}}$ in coordinates, we get $(\Phi_s^*v_k)_z=e^{ls}(v_k)_{\Phi_s(z)}$ and
	$$(L_{y-y_k}v_k)_z=\left.\frac{d}{ds}\right|_{s=0}(\Phi_s^*v_k)_z=l(v_k)_z+\sum_{i=1}^{\dim(M)}(z-y_k)_i\frac{\de f_{\overline{j}}}{\de r_i}(z)dr_{\overline{j}},$$
	where the first equality is the definition of the Lie derivative.
	Consequently, 
	\begin{equation*}
		\begin{split}
			&t^{l-1}(L_{y-y_k}v_k)_{y_k(1-t)+yt}=\\
			&=t^{l-1}\left(l(v_k)_{y_k(1-t)+yt}+\sum_{\overline{j}\in\mathcal{J}^l}\sum_{i=1}^{\dim(M)}t(y-y_k)_i\frac{\de f_{\overline{j}}}{\de r_i}(y_k(1-t)+yt)dr_{\overline{j}}\right)=\\
			&=\left(lt^{l-1}(v_k)_{y_k(1-t)+yt}+t^l\sum_{\overline{j}\in\mathcal{J}^l}\langle\nabla f_{\overline{j}}(y_k(1-t)+yt),y-y_k\rangle dr_{\overline{j}}\right)\\
			&=\frac{d}{dt}\left(t^{l}(v_k)_{y_k(1-t)+yt}\right)
		\end{split}
	\end{equation*}
	hence, from \eqref{eq_ap1}, we get
	$$((d\alpha_k+\alpha_kd)\omega)_x=\psi_{k*}\int_0^1 \frac{d}{dt}\left(t^l\psi_k^*\omega_{y_k(1-t)+yt}\right)dt=\omega_x.$$
	We have thus proved that $(d\alpha_k+\alpha_kd)\omega=\omega.$ Accordingly, if $\omega$ is a $l$-form such that $d\omega\phi_k=0,$ then we get $d\alpha_k(\omega)\phi_k=\omega\phi_k-\alpha_k(d\omega)\phi_k=\omega\phi_k,$ that is we have just proved  Poincaré's lemma for differential forms. Once we have the result for forms, we could try to extend $\alpha_k$ to a bounded linear operator from $\mathcal{C}^{p,q,l}$ to $\mathcal{C}^{p,q,l-1}.$ Unfortunately, $\alpha_k(\omega),$ for a differential form $\omega,$ is not defined on every admissible stable leaf of $W\in\Sigma$. In fact, we can just consider leaves inside $\psi_k(U_k).$ On the other hand, the product $\alpha_k(\omega)\phi_k^+$ is well defined on the full manifold, because it is null out of $\text{supp}(\phi_k^+).$ Accordingly, we consider the operator $\beta_k\colon\Omega^l(M)\rightarrow\Omega^{l-1}(M)$, such that $\beta_k(\omega)=\alpha_k(\omega)\phi_k^+$, and we want to prove that it can be extended to a bounded linear operator from $\mathcal{C}^{p,q,l}$ to $\mathcal{C}^{p,q,l-1}.$ 
	Let $\omega\in\Omega^l(M),$ and $\psi_k\circ G_{p,F}(\mathcal{B}_{d_s}(0,\rho))=W\in\Sigma,$ with $p\in\mathcal{B}(0,\rho)$ and $F\in\mathcal{F},$ as defined in Section \ref{section_leaves}. Let  $\phi\in\Gamma_0^{q,l-1}(W)$ be a test form. Then, denoting again by $y=\psi_k^{-1}(x)$, $y_k=\psi_k^{-1}(x_k)$, $v_k=\psi_k^*\omega,$ and $\chi_k=\psi_k^*\phi$
	\begin{equation}
		\label{eq_ap2}
		\begin{split}
			&\int_W\!\!\!\langle\phi_x,\alpha_k(\omega)_x\rangle\phi_k^+(x)\omega_W(x)\!\!=\!\!\!\int_W\!\!\!\langle\phi_x,\psi_{k*}\left(\int_0^1\!\!\!t^{l-1}(\iota_{y-y_k}v_k)_{y_k(1-t)+yt}dt\right)_x\!\!\rangle\phi_k^+(x)\omega_W(x)=\\
			=&\psi_{k*}\left(\int_{\psi_k^{-1}W}\langle(\chi_k)_y,\int_0^1t^{l-1}(\iota_{y-y_k}v_k)_{y_k(1-t)+yt}dt\rangle\phi_k^+(\psi_k(y))\psi_k^*\omega_W(y)\right).
		\end{split}
	\end{equation}
	Setting $l_t(y)=y_k(1-t)+yt,$ we have $t^{l-1}(i_{y-y_k}v_k)_{l_t(y)}=(l_{t*}i_{y-y_k}v_k)_y$ and  $(\chi_k)_y=(l_t^*l_{t*}\chi_k)_y=t^{l-1}(l_{t*}\chi_k)_{l_t(y)}.$ Thus, from \eqref{eq_ap2}, we obtain
	\begin{equation}
		\label{eq_ap3}
		\begin{split}
		&\int_W\!\!\!\langle\phi_x,\alpha_k(\omega)_x\rangle\phi_k^+(x)\omega_W(x)\!\!=\\=&\psi_{k*}\left(\int_{\psi_k^{-1}W}\int_0^1\langle(\chi_k)_y,(l_t^*\iota_{y-y_k}v_k)_y\rangle \phi_k^+(\psi_k(y))dt\psi_k^*\omega_W(y)\right)=\\
			=&\psi_{k*}\left(\int_{\psi_k^{-1}W}\int_0^1t^{l-1}\langle l_{t*}(\chi_k)_{l_t(y)},(\iota_{y-y_k}v_k)_{l_t(y)}\rangle \phi_k^+(\psi_k(y)) dt\psi_k^*\omega_W(y)\right)=\\
			=&\!\!\!\sum_{i=1}^{\dim(M)}\psi_{k*}\left(\int_{\psi_k^{-1}W}\int_0^1t^{l-1}\langle l_{t*}\chi_k\wedge dr_i,v_k\rangle_{l_t(y)}(y-y_k)_i\phi_k^+(\psi_k(y))dt\psi_k^*\omega_W(y)\!\!\right)
		\end{split}
	\end{equation} 
	Next, we compute
	\begin{equation*}
		\begin{split}
			&\int_{\psi_k^{-1}W}\int_0^1\!\!\!\!t^{l-1}\langle l_{t*}\chi_k\!\wedge\! dr_i,v_k\rangle_{l_t(y)}(y-y_k)_i\phi_k^+(\psi_k(y))dt\psi_k^*\omega_W(y)=\\
			=&\!\!\!\int_{\mathcal{B}_{d_s}(0,\rho)}\!\int_0^1\!\!\!\!t^{l-1}\langle l_{t*}\chi_k\!\wedge\! dr_i,v_k\rangle_{l_t(G_{p,F}(s))}(G_{p,F}(s)\!-\!y_k)_iG_{p,F}^*\psi_k^*\phi_k^+(s)dtG_{p,F}^*\psi_k^*\omega_W(s)=\\
			=&\!\!\int_0^1\!\!\!\!\int_{\mathcal{B}_{d_s}(0,t\rho)}\!\!\!\!\!\!\!\!\!\!\!\!\!\!\!\!\!t^{l-1-d_s}\langle l_{t*}\chi_k\!\wedge\! dr_i,v_k\rangle_{l_t((G_{p,F}(t^{-1}\!z)))}(\!G_{p,F}(t^{-1}\!z)\!-\!y_k\!)_iG_{p,F}^*\psi_k^*\phi_k^+(\!t^{-1}\!z\!)G_{p,F}^*\psi_k^*\omega_W(\!z\!)dt
		\end{split}
	\end{equation*}
	Let $q_t=y_k(1-t)+tp$ and $F_t\colon\mathcal{B}_{d_s}(0,t\rho)\rightarrow\mathcal{B}_{d_u}(0,\rho)$ such that $F_t(z)=F(t^{-1}z).$ Then $q_t\in\mathcal{B}(0,\rho)$ and $F_t$ can be extended to a function $\bar{F_t}\in\mathcal{F}.$ As a consequence, setting $W_t=\psi_k\circ G_{q_t,F_t}(\mathcal{B}_{d_s}(0,t\rho))\subseteq \widehat{W_t}=\psi_k\circ G_{q_t,F_t}(\mathcal{B}_{d_s}(0,\rho))\in\Sigma$ and
	\begin{equation}
		\label{eq_bar_phi}
		(\bar\phi_t)_x=\sum_{i=1}^{\dim(M)}(\psi_{k*}l_{t*}\chi_k\wedge dx_i)_x(\psi_k^{-1}(x)-y_k)_i\phi_k^+(x),
	\end{equation}
	we obtain, continuing equation \eqref{eq_ap3},
	\begin{equation*}
		\int_W\!\!\!\langle\phi_x,\alpha_k(\omega)_x\rangle\phi_k^+(x)\omega_W(x)\!\!=\!\!\int_0^1t^{l-1-d_s}\!\!\!\int_{W_t}\langle(\bar{\phi}_t)_x,\omega_x\rangle\omega_{W_t}(x)dt.
	\end{equation*}
	In summary, we have shown that 
	\begin{equation*}
		\begin{split}
			\left|\int_W\!\!\!\langle\phi_x,\alpha_k(\omega)_x\rangle\phi_k^+(x)\omega_W(x)\right|\leq C\int_0^1 t^{l-1-d_s} \norm{\bar{\phi}_t}_{\Gamma_0^{q,l}}\norm{\omega}_{0,q,l}|W_t| dt
		\end{split}
	\end{equation*}
	One can easily check, using \eqref{eq_bar_phi}, that $\norm{\bar{\phi}_t}_{\Gamma_0^{q,l}}\leq C\norm{\phi}_{\Gamma_0^{q,l-1}}.$ In addition, $|W_t|\leq Ct^{d_s}|\widehat{W_t}|\leq Ct^{d_s},$ hence we conclude that
	\begin{equation}
		\label{eq_first_norm}
		\begin{split}
			\left|\int_W\!\!\!\langle\phi_x,\alpha_k(\omega)_x\rangle\phi_k^+(x)\omega_W(x)\right|&\leq C\int_0^1 t^{l-1} \norm{\phi}_{\Gamma_0^{q,l-1}}\norm{\omega}_{0,q,l} dt\\&\leq C\norm{\phi}_{\Gamma_0^{q,l-1}}\norm{\omega}_{0,q,l},
		\end{split}
	\end{equation}
	i.e., $\norm{\beta_k(\omega)}_{0,q,l-1}\leq C\norm{\omega}_{0,q,l}$ and $\beta_k\colon\mathcal{B}^{0,q,l}\rightarrow\mathcal{B}^{0,q,l-1}$ is a bounded linear operator. 
	Similarly, let $\gamma_k\colon\Omega^l(M)\rightarrow\Omega^l(M),$ such that 
	$\gamma_k(\omega)=\alpha_k(\omega)\wedge d\phi_k^+. $ Then, by the same proof, one can easily check that it can be extended to a bounded linear operator $\gamma_k\colon\mathcal{B}^{0,q,l}\rightarrow\mathcal{B}^{0,q,l}.$ 
	Next, we compute
	\begin{equation}
		\label{eq_third_norm}
		\begin{split}
			&\int_W\langle\phi,d\beta_k(\omega)\rangle\omega_W=\int_W\langle\phi,d\alpha_k(\omega)\rangle\phi_k^+\omega_W+(-1)^l\int_W\langle\phi,\alpha_k(\omega)\wedge d\phi_k^+\rangle \omega_W=\\
			=&\int_W\langle\phi,\omega\rangle\phi_k^+\omega_W-\int_W\langle\phi,\alpha_k(d\omega)\rangle\phi_k^+\omega_W+(-1)^l\int_W\langle\phi,\alpha_k(\omega)\wedge d\phi_k^+\rangle \omega_W,
		\end{split}
	\end{equation}
	hence
	\begin{equation}
		\label{eq_second_norm}
		\begin{split}
			\left|\int_W\!\!\!\!\langle\phi,d\beta_k(\omega)\rangle\omega_W\right|&\leq C\norm{\phi}_{\Gamma_0^{q,l}}\norm{\omega}_{0,q,l}+C\norm{\phi}_{\Gamma_0^{q,l}}\norm{\beta_k(d\omega)}_{0,q,l}+C\norm{\phi}_{\Gamma_0^{q,l}}\norm{\gamma_k(d\omega)}_{0,q,l+1}\\&\leq  C\norm{\phi}_{\Gamma_0^{q,l}}\norm{\omega}_{0,q,l}+C\norm{\phi}_{\Gamma_0^{q,l}}\norm{d\omega}_{0,q,l+1},
		\end{split}
	\end{equation}
	where, in the last inequality, we used the continuity of $\beta_k\colon\mathcal{B}^{0,q,l+1}\rightarrow\mathcal{B}^{0,q,l}$ and  $\gamma_k\colon\mathcal{B}^{0,q,l+1}\rightarrow\mathcal{B}^{0,q,l+1}.$
	From \eqref{eq_first_norm} and \eqref{eq_second_norm} we get 
	\begin{equation*}
		\begin{split}
			|\beta_k(\omega)|_{0,q,l-1}=\norm{\beta_k(\omega)}_{0,q,l-1}+\norm{d\beta_k(\omega)}_{0,q,l}\leq C\norm{\omega}_{0,q,l}+C\norm{d\omega}_{0,q,l+1}\leq C|\omega|_{0,q,l}
		\end{split},
	\end{equation*}
	that is $\beta_k\colon\mathcal{C}^{0,q,l}\rightarrow\mathcal{C}^{0,q,l-1}$ is continuous. Since $\omega$ is closed, $d\gamma_k(\omega)=\omega\wedge d\phi_k^+$ and $\norm{\gamma_k(\omega)}_{0,q,l}\leq C\norm{\omega}_{0,q,l}.$ Therefore, $|\gamma_k(\omega)|_{0,q,l}\leq C|\omega|_{0,q,l}$ and also $\gamma_k\colon\mathcal{C}^{0,q,l}\rightarrow\mathcal{C}^{0,q,l}$ is bounded. 
	
	We now want to prove that $\beta_k\colon\mathcal{C}^{1,q,l}\rightarrow\mathcal{C}^{1,q,l-1}$ is still a bounded linear operator. 
	Therefore, let $\omega\in\Omega^l(M),$  $W\in\Sigma,$ $\phi\in\Gamma_0^{q+1,l-1}(W)$ and $v\in\mathcal{V}^{q+1}(U(W)).$ 
	We firstly show that $\iota_v\alpha_k(\omega)=-\alpha_k(\iota_v\omega).$ In fact, once we have this preliminary result, we get $L_v\alpha_k(\omega)=(d\iota_v+\iota_vd)\alpha_k(\omega)=-d\alpha_k(\iota_v\omega)+\iota_v\omega-\iota_v\alpha_k(d\omega)=-\iota_v\omega+\alpha_k(d\iota_v\omega)+\iota_v\omega+\alpha_k(\iota_vd\omega)=\alpha_k(L_v\omega).$
	With a slight abuse of notation, in order to simplify notations,  we confuse $\omega$, $x$, $x_k$ and $v$, with $\psi_{k*}\omega$, $\psi_k(x)$, $\psi_k(x_k)$ and $\psi_{k*}v$, respectively. Accordingly, we compute
	\begin{equation*}
		\begin{split}
			\iota_v\alpha_k(\omega)=&\iota_v\left(\int_0^1t^{l-1}(\iota_{x-x_k}\omega)_{x_k(1-t)+xt}dt\right)=\int_0^1\iota_v(l_t^*\iota_{x-x_k}\omega)_xdt\\=&\int_0^1(l_t^*\iota_{l_{t*}v}\iota_{x-x_k}\omega)_xdt=\int_0^1t^{-1}(l_t^*\iota_{v}\iota_{x-x_k}\omega)_xdt=\\=&-\int_0^1t^{-1}(l_t^*\iota_{x-x_k}\iota_v\omega)_xdt=-\int_0^1t^{l-2}\iota_{x-x_k}(\iota_v\omega)_{x_k(1-t)+tx}dt=-\alpha_k(\iota_v\omega).\\
		\end{split}
	\end{equation*}
	As a consequence,
	\begin{equation*}
		\begin{split}
			\int_W\langle\phi, L_v\beta_k(\omega)\rangle\omega_W=&\int_W\langle\phi, L_v\alpha_k(\omega)\rangle\phi_k^+\omega_W+\int_W\langle\phi, \alpha_k(\omega)v\phi_k^+\rangle\omega_W=\\
			=&\int_W\langle\phi, \alpha_k(L_v\omega)\rangle\phi_k^+\omega_W+\int_W\langle\phi, \alpha_k(\omega)v\phi_k^+\rangle\omega_W\\=&\int_W\langle\phi, \beta_k(L_v\omega)\rangle\omega_W+\int_W\langle\phi, \alpha_k(\omega)v\phi_k^+\rangle\omega_W,,
		\end{split}
	\end{equation*}
	hence
	\begin{equation*}
		\begin{split}
			\left|\int_W\!\!\!\!\langle\phi, L_v\beta_k(\omega)\rangle\omega_W\right|\leq& C\norm{\phi}_{\Gamma_{0}^{q+1,l-1}}\norm{\beta_k(L_v\omega)}_{0,q+1,l-1}+C\norm{\phi}_{\Gamma_0^{q+1,l-1}}\norm{v}_{C^{q+1}}\norm{\omega}_{0,q+1,l}\leq\\
			\leq& C\norm{\phi}_{\Gamma_{0}^{q+1,l-1}}\norm{v}_{C^{q+1}}\norm{\omega}_{1,q,l}+C\norm{\phi}_{\Gamma_0^{q+1,l-1}}\norm{v}_{C^{q+1}}\norm{\omega}_{0,q+1,l},
		\end{split}
	\end{equation*}
	which proves that $\norm{\beta_k(\omega)}_{1,q,l-1}\leq C\norm{\omega}_{1,q,l}.$
	By the same procedure of \eqref{eq_third_norm} we obtain
	\begin{equation*}
		\begin{split}
			\int_W\!\!\!\!\langle\phi, L_vd\beta_k(\omega)\rangle\omega_W=&\int_W\!\!\!\!\langle\phi,L_v\omega\rangle\phi_k^+\omega_W-\int_W\!\!\!\!\langle\phi,\alpha_k(L_vd\omega)\rangle\phi_k^+\omega_W+\\+&(-1)^l\int_W\!\!\!\!\langle\phi,\alpha_k(L_v\omega)\wedge d\phi_k^+\rangle \omega_W,
		\end{split}
	\end{equation*}
	hence
	\begin{equation*}
		\begin{split}
			\left|\int_W\!\!\!\!\langle\phi, L_vd\beta_k(\omega)\rangle\omega_W\right|\leq& C\norm{\phi}_{\Gamma_0^{q+1,l-1}}\norm{L_v\omega}_{0,q+1,l}+C\norm{\phi}_{\Gamma_0^{q+1,l-1}}\norm{\beta_k(L_vd\omega)}_{0,q+1,l}+\\+&C\norm{\phi}_{\Gamma_0^{q+1,,l-1}}\norm{\gamma_k(L_v\omega)}_{0,q+1,l}\leq C\norm{\phi}_{\Gamma_0^{q+1,l-1}}\norm{v}_{C^{q+1}}\norm{\omega}_{1,q,l}+\\+&C\norm{\phi}_{\Gamma_0^{q+1,l-1}}\norm{v}_{C^{q+1}}\norm{d\omega}_{1,q,l+1}+C\norm{\phi}_{\Gamma_0^{q+1,l-1}}\norm{v}_{C^{q+1}}\norm{\omega}_{1,q,l},
		\end{split}
	\end{equation*}
	i.e $\norm{d\beta_k(\omega)}_{1,q,l}\leq C\norm{\omega}_{1,q,l}+C\norm{d\omega}_{1,q,l+1}.$ We conclude that $$|\beta_k(\omega)|_{1,q,l-1}=\norm{\beta_k(\omega)}_{1,q,l-1}+\norm{d\beta_k(\omega)}_{1,q,l}\leq C\norm{\omega}_{1,q,l}+C\norm{d\omega}_{1,q,l+1}\leq C|\omega|_{1,q,l},$$ 
	so that $\beta_k\colon\mathcal{C}^{1,q,l}\rightarrow\mathcal{C}^{1,q,l-1}$ is a bounded linear operator. A similar computation shows that $\gamma_k\colon\mathcal{C}^{1,q,l}\rightarrow\mathcal{C}^{1,q,l}$ is a continuous operator. By the same argument $\beta_k$ and $\gamma_k$ extend to bounded linear operators $\beta_k\colon\mathcal{C}^{p,q,l}\rightarrow\mathcal{C}^{p,q,l-1}$ and $\gamma_k\colon\mathcal{C}^{p,q,l}\rightarrow\mathcal{C}^{p,q,l}.$ Finally, since $\phi_k^+\phi_k=\phi_k,$ we obtain 
	\begin{equation*}
		\begin{split}
			d(\beta_k(\omega)\phi_k)=&d\beta_k(\omega)\phi_k+(-1)^{l-1}\beta_k(\omega)\wedge d\phi_k=\\=&\omega\phi_k^+\phi_k+(-1)^{l-1}\alpha_k(\omega)\wedge d\phi_k^+\phi_k+(-1)^{l-1}\beta_k(\omega)\wedge d\phi_k=\\=&\omega\phi_k+(-1)^{l-1}\beta_k(\omega)\wedge d\phi_k,
		\end{split}
	\end{equation*}
	hence, setting $u_k=\beta_k(\omega)$ we get the lemma.
	
	\qed
	
	\proofof{Lemma \ref{lemma_Poincare2}}
	Let $h\in\mathcal{C}^{p,q,l,0}$ such that $\bar{\delta}h=0.$ Then $0=\bar{\delta}h(a_0,a_1)=h(a_1)-h(a_0),$ i.e., $h(a_0)=h(a_1)$, for all $a_0,a_1\in\{1,\dots,m\}.$ Accordingly, defining $\omega=h(a_0)\in\mathcal{C}^{p,q,l},$ we get $i(\omega)=h.$ 
	
	Next, we consider $\omega\in\mathcal{C}^{p,q,l,k},$ for some $k>0$, such that $\bar{\delta}\omega=0.$ We define the operator $\Xi\colon\mathcal{C}^{p,q,l,k}\rightarrow\mathcal{C}^{p,q,l,k-1}$ such that 
	$$\Xi(\omega)(a_0,\dots,a_{k-1})=\sum_{j=1}^m\phi_j\omega(j,a_0,\dots,a_{k-1}),$$ 
	where $\{\phi_j\}_{j=1}^m$ is the partition of unity subordinated to the cover $\{V_j=\psi_j(U_j)\}_{j=1}^m$ that we used to define \v{C}ech cohomology.
	We show that $\bar{\delta}\circ\Xi+\Xi\circ\bar{\delta}=\id.$ Indeed,
	\begin{equation*}
		\begin{split}
			\bar\delta\circ\Xi\omega(a_0,\dots,a_k)=&\sum_{i=0}^k(-1)^i\Xi\omega(a_0,\dots,a_{i-1},a_{i+1},\dots,a_k)=\\=&\sum_{i=0}^k(-1)^i\sum_{j=1}^m\phi_j\omega(j,a_0,\dots,a_{i-1},a_{i+1},\dots,a_k),
		\end{split}
	\end{equation*}
	while
	\begin{equation*}
		\begin{split}
			\Xi\circ\bar\delta\omega(a_0,\dots,a_k)=&\sum_{j=1}^m\phi_j\bar\delta\omega(j,a_0,\dots,a_k)=\sum_{j=1}^m\phi_j\omega(a_0,\dots,a_k)+\\+&\sum_{i=0}^k(-1)^{i+1}\sum_{j=1}^m\phi_j\omega(j,a_0,\dots,a_{i-1},a_{i+1},\dots,a_k)=\\=&(\omega-\bar{\delta}\circ\Xi\omega)(a_0,\dots,a_k).
		\end{split}
	\end{equation*}
	Accordingly, if $\bar{\delta}(\omega)=0,$ defining $u=\Xi\omega,$ we get 
	$$\bar{\delta}u=\bar{\delta}\circ\Xi\omega=\omega-\Xi\circ\bar\delta\omega=\omega.$$
	
	\qed
\end{appendices}

\section*{Funding and/or Conflicts of interests/Competing interests}
	Partial financial support was received from the Marie Sklodowska-Curie grant agreement No 777822 (‘GHAIA’), and the PRIN Grant 2017S35EHN, MUR, Italy. The author has no relevant financial or non-financial interests to disclose.

%%% Section references
\bibliographystyle{plain}
\bibliography{res_revised}

\end{document}